\documentclass{article}

\usepackage{arxiv}

\usepackage[utf8]{inputenc} 
\usepackage[T1]{fontenc}    
\usepackage{hyperref}       
\usepackage{url}            
\usepackage{booktabs}       
\usepackage{amsfonts}       
\usepackage{nicefrac}       
\usepackage{microtype}      
\usepackage{lipsum}		
\usepackage{graphicx}
\usepackage{natbib}
\usepackage{doi}
\usepackage{caption}
\usepackage{subcaption}
\usepackage{algorithm}
\usepackage{algorithmic}

\usepackage{amsmath}
\usepackage{amssymb}

\newtheorem{theorem}{Theorem}[section]
\newtheorem{lemma}[theorem]{Lemma}

\newcommand{\blackslug}{\penalty 1000\hbox{
    \vrule height 8pt width .4pt\hskip -.4pt
    \vbox{\hrule width 8pt height .4pt\vskip -.4pt
          \vskip 8pt
          \vskip -.4pt\hrule width 8pt height .4pt}
    \hskip -3.9pt
    \vrule height 8pt width .4pt}}

\newenvironment{proof}{$\;$\newline \noindent {\sc Proof.}$\;\;\;$\rm}{\qed}
\newcommand{\qed}{\hspace*{\fill}\blackslug}
\newenvironment{definition}{$\;$\newline \noindent {\bf Definition}$\;$}{$\;$\newline}
\def\boxit#1{\vbox{\hrule\hbox{\vrule\kern4pt
  \vbox{\kern1pt#1\kern1pt}
\kern2pt\vrule}\hrule}}

\title{Projective H\"{o}lder-Minkowski Colors: \\A Generalized Set of Commutative \& Associative Operations with Inverse Elements \\for Representing and Manipulating Colors}

\author{ 
\href{https://orcid.org/0000-0003-3618-4166}{\includegraphics[scale=0.06]{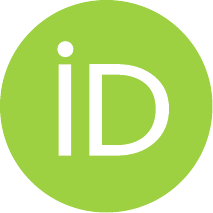}\hspace{1mm}Ergun Akleman}\thanks{Joint with Computer Science and Engineering Department.} \\
	Visual Computing \& Computational Media, \\School of PVFA\\
 Texas A\&M University, \\College Station, TX, 77831\\
	\texttt{ergun@tamu.edu} \\
 \And
        Somyung (David) Oh\\
Computer Science and Engineering\\ 
College of Engineering \\
Texas A\&M University, \\College Station, TX, 77831\\
	\texttt{somyungoh@gmail.com } \\
 	\And 
       Youyou Wang\\
Computer Science and Engineering\\ 
Texas A\&M University, \\
College Station, TX, 77831\\
	\texttt{kingyy2010@gmail.com} \\
     \And    
   Bekir Tevfik Akgun\\
    College of Computer and Information Sciences\\
	Yeditepe University, \\
 Istanbul, Turkey\\
	\texttt{bekirtevfik.akgun@yeditepe.edu.tr} \\
   \And
Jianer Chen\\
 Computer Science and Engineering, \\
 College of Engineering\\
 Texas A\&M University, \\College Station, TX, 77831\\
	\texttt{chen@cs.tamu.edu} \\
}

\renewcommand{\shorttitle}{H\"{o}lder-Minkowski Colors}

\hypersetup{
pdftitle={A template for the arxiv style},
pdfsubject={q-bio.NC, q-bio.QM},
pdfauthor={David S.~Hippocampus, Elias D.~Striatum},
pdfkeywords={First keyword, Second keyword, More},
}

\begin{document}
\maketitle

\begin{abstract}
One of the key problems in dealing with color in rendering, shading, compositing, or image manipulation is that we do not have algebraic structures that support operations over colors. In this paper, we present an all-encompassing framework that can support a set of algebraic structures with associativity, commutativity, and inverse properties. To provide these three properties, we build our algebraic structures on an extension of projective space by allowing for negative and complex numbers. These properties are important for (1) manipulating colors as periodic functions, (2) solving inverse problems dealing with colors, and (3) being consistent with the wave representation of the color. Allowance of negative and complex numbers is not a problem for practical applications, since we can always convert the results into desired range for display purposes as we do in High Dynamic Range imaging. This set of algebraic structures can be considered as a generalization of the Minkowski norm $L_p$ in projective space. These structures also provide a new version of the generalized H\"{o}lder average with associativity property. Our structures provide inverses of any operation by allowing for negative and complex numbers. These structures provide all properties of generalized H\"{o}lder average by providing a continuous bridge between the classical weighted average, harmonic mean, maximum, and minimum operations using a single parameter $p$. The parameter $p=1$ corresponds to a classical weighted average operation, and $p=-1$ provides a harmonic mean. The $p \rightarrow \infty$ produces the maximum operator and the $p \rightarrow -\infty$ minimum operator. Our main contribution is the demonstration of the associative property in a projective space. We have also shown that these operations are scale invariant. Our version is slightly different than the classical H\"{o}lder average, but it still guarantees that the results of the operations are always between the maximum and minimum numbers. We also show that these operations are qualitatively similar to a weighted average. Any of these structures can be directly implemented using current shader languages by treating three color channels as spatial positions. For more color channels or a general version of H\"{o}lder-Minkowski colors, there is a need for an extension that directly allows vectors of complex numbers in shading languages. 
\end{abstract}

 \begin{figure*}
    \centering  
        \begin{subfigure}[t]{0.32\textwidth}
        \includegraphics[width=1.0\textwidth]{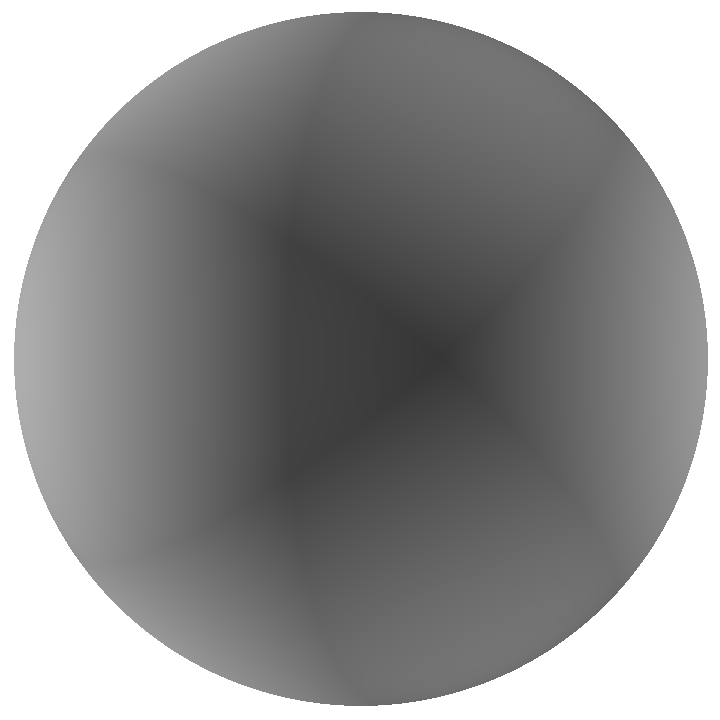}
        \caption{$p=1.0$}
        \label{fig_images/hedron/2}
    \end{subfigure}
      \hfill 
        \begin{subfigure}[t]{0.32\textwidth}
        \includegraphics[width=1.0\textwidth]{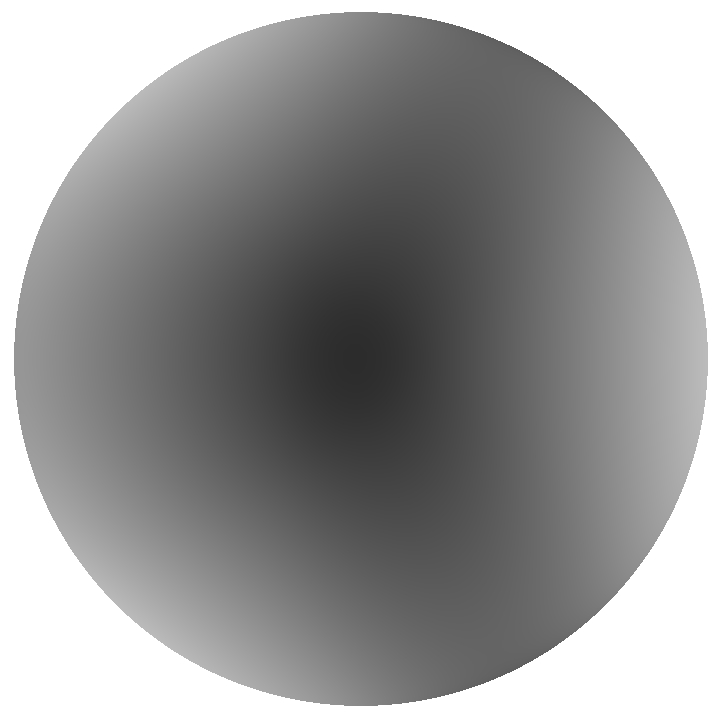}
        \caption{$p=2.0$}
        \label{fig_images/hedron/4}
    \end{subfigure}
      \hfill
        \begin{subfigure}[t]{0.32\textwidth}
        \includegraphics[width=1.0\textwidth]{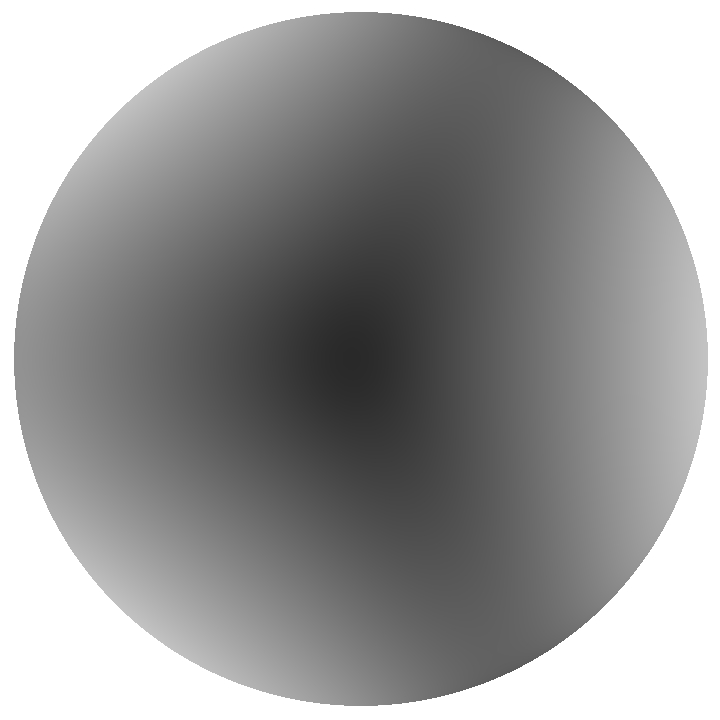}
        \caption{$p=2.5$}
        \label{fig_images/hedron/5}
         \end{subfigure}
        \hfill
       \begin{subfigure}[t]{0.32\textwidth}
        \includegraphics[width=1.0\textwidth]{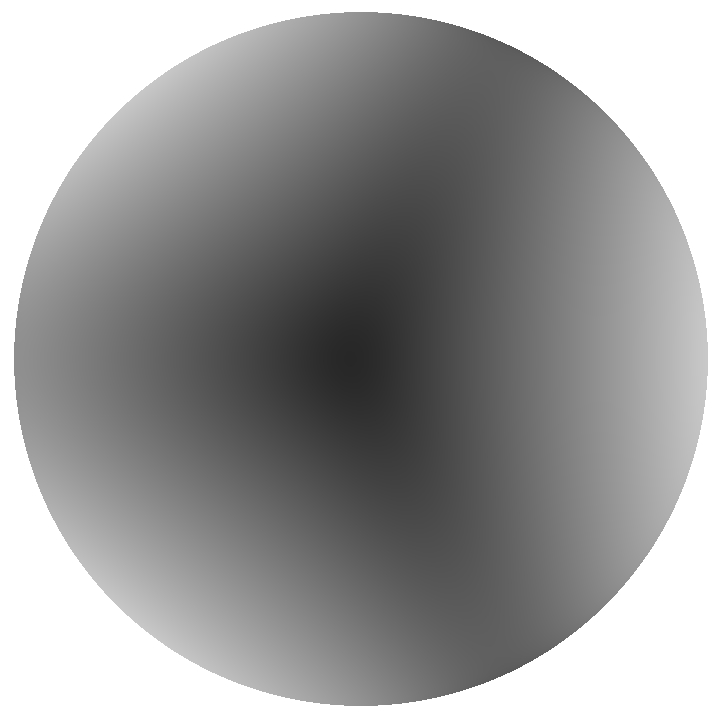}
        \caption{$p=3.0$}
        \label{fig_images/hedron/6}
         \end{subfigure}
        \hfill
        \begin{subfigure}[t]{0.32\textwidth}
        \includegraphics[width=1.0\textwidth]{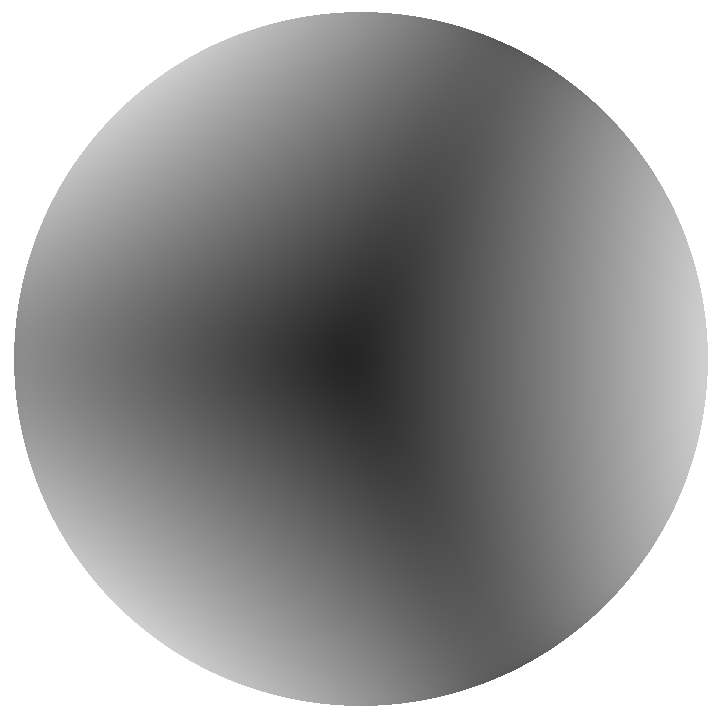}
        \caption{$p=4.0$}
        \label{fig_images/hedron/8}
         \end{subfigure}
        \hfill
        \begin{subfigure}[t]{0.32\textwidth}
        \includegraphics[width=1.0\textwidth]{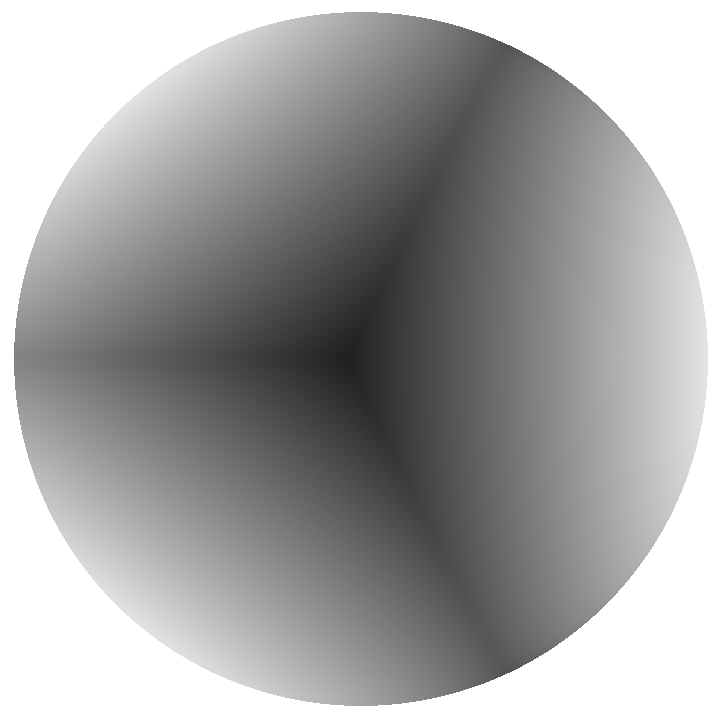}
        \caption{$p \rightarrow \infty$}
        \label{fig_images/hedron/9}
         \end{subfigure}
      \hfill
\caption{A sphere illuminated by four point-lights located in the vertices of a tetrahedron. The images show the impact of different parameters that are based on different weighted averages of the direct illumination term $max(\cos \theta,0)$, which is a discontinuous derivative. The case $p=1$ corresponds to the classical weighted average (or addition), visual artifacts caused by the derivative discontinuity of the function $max(\cos \theta,0)$. Parameter values greater than $p=1$ create smoother results by removing these visual artifacts caused by the discontinuity of the derivative at zero. Very high values $p$ approach the maximum operator and create another visual artifact. The $p$ values smaller than $1$ are not meaningful; therefore, we did not include them. }
\label{fig_teaser}
\end{figure*}

 \begin{figure*}
    \centering  
        \begin{subfigure}[t]{0.32\textwidth}
        \includegraphics[width=1.0\textwidth]{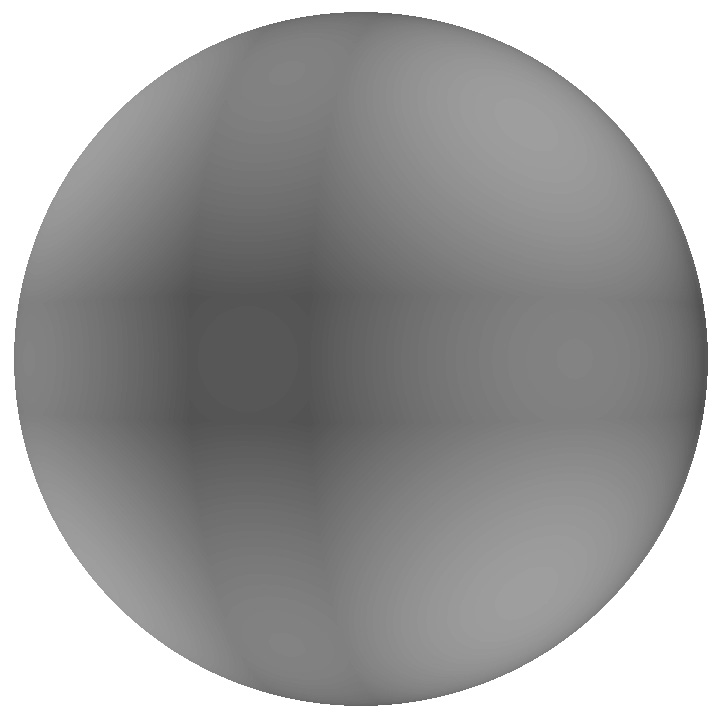}
        \caption{$p=1.0$}
        \label{fig_images/octahedron/2}
    \end{subfigure}
      \hfill 
        \begin{subfigure}[t]{0.32\textwidth}
        \includegraphics[width=1.0\textwidth]{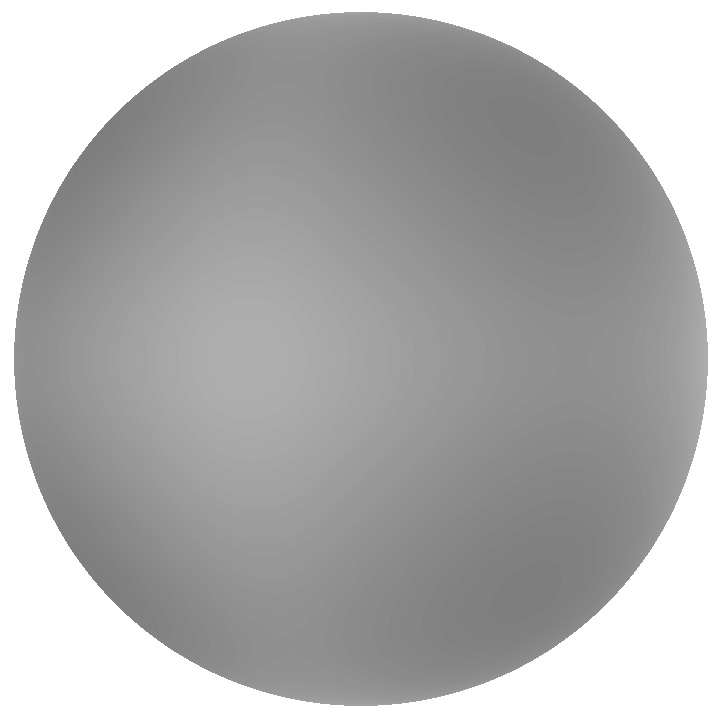}
        \caption{$p=2.0$}
        \label{fig_images/octahedron/4}
    \end{subfigure}
      \hfill
        \begin{subfigure}[t]{0.32\textwidth}
        \includegraphics[width=1.0\textwidth]{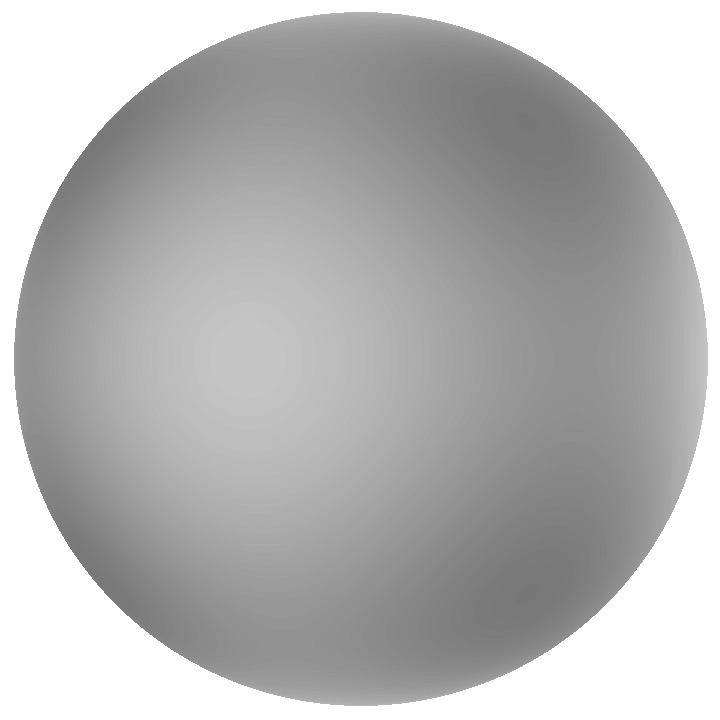}
        \caption{$p=2.5$}
        \label{fig_images/octahedron/5}
         \end{subfigure}
        \hfill
       \begin{subfigure}[t]{0.32\textwidth}
        \includegraphics[width=1.0\textwidth]{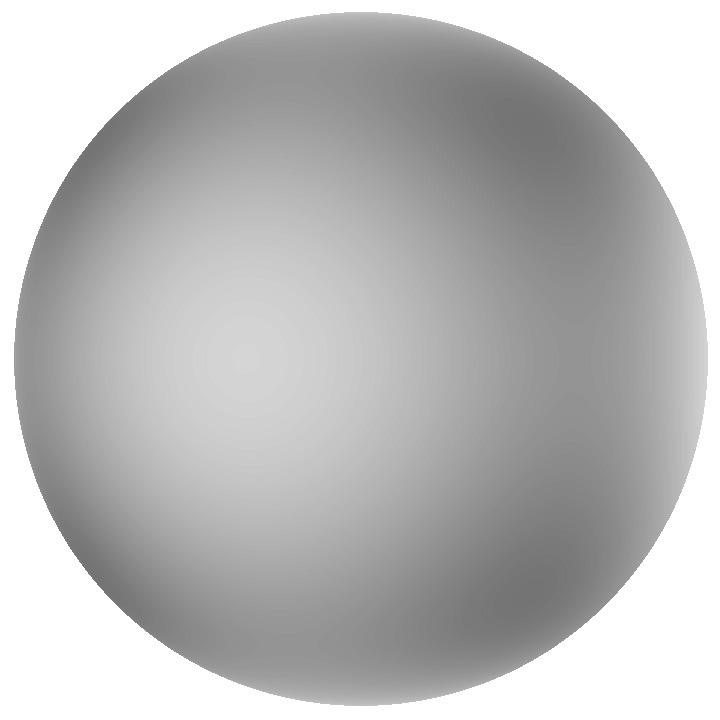}
        \caption{$p=3.0$}
        \label{fig_images/octahedron/6}
         \end{subfigure}
        \hfill
        \begin{subfigure}[t]{0.32\textwidth}
        \includegraphics[width=1.0\textwidth]{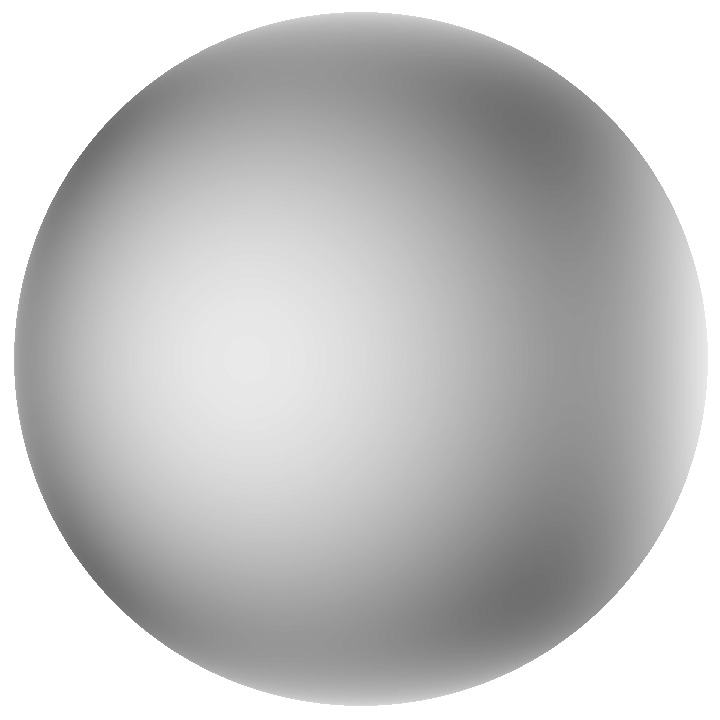}
        \caption{$p=4.0$}
        \label{fig_images/octahedron/8}
         \end{subfigure}
        \hfill
        \begin{subfigure}[t]{0.32\textwidth}
        \includegraphics[width=1.0\textwidth]{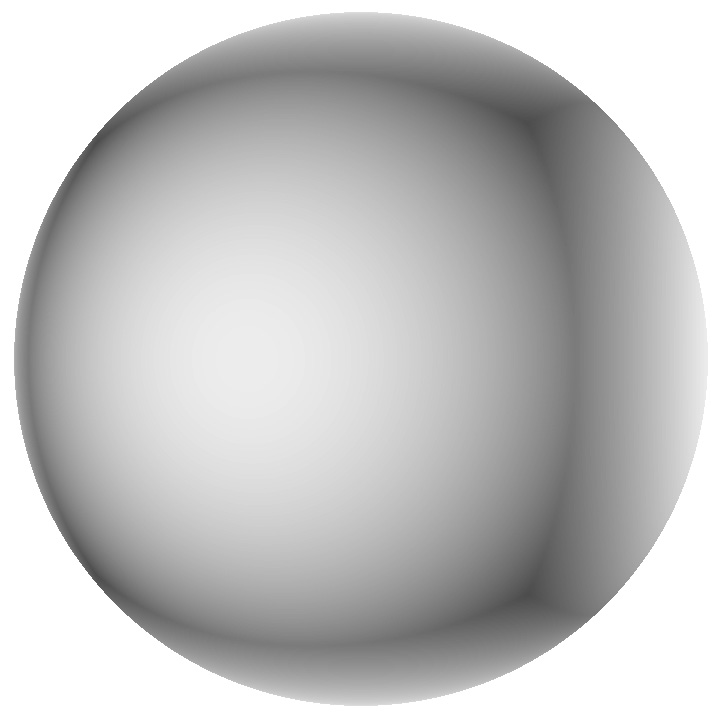}
        \caption{$p \rightarrow \infty$}
        \label{fig_images/octahedron/9}
         \end{subfigure}
      \hfill
\caption{A sphere illuminated by six point-lights located in the vertices of an octahedron. The images show the impact of different parameters that are based on different weighted averages of the direct illumination term $max(\cos \theta,0)$, which is a discontinuous derivative. The case $p=1$ corresponds to the classical weighted average (or addition), visual artifacts caused by the derivative discontinuity of the function $max(\cos \theta,0)$. Parameter values greater than $p=1$ create smoother results by removing these visual artifacts caused by the discontinuity of the derivative at zero. Very high values $p$ approach the maximum operator and create another visual artifact. The $p$ values smaller than $1$ are not meaningful; therefore, we did not include them. }
\label{fig_teaser2}
\end{figure*}

\section{Introduction}\label{sec_introduction} 

Operations over colors are one of the key elements of computer graphics, from rendering to compositing and image manipulations. There is currently a wide range of operations, from addition and multiplication to maximum, minimum, and average, that are used to manipulate colors. However, these operations are not part of any algebraic framework, which can cause difficulties in solving problems with colors. In this paper, we propose to represent colors as complex functions of time over projective space, i.e. in homogeneous coordinates. The projective space representation will allow us to formally differentiate ``point (or H\"{o}lder) functions'' and ``vector (or Minkowski) functions''. This formalization also allows us to differentiate between light and material colors. We further allow Minkowski operations, which naturally provide a conceptually consistent framework to operate on point functions and vector functions on projective space. This representation along with Minkowski operations, for instance, naturally provides addition over vector functions and weighted average over point functions.

The weighted average is among the most common operations in a wide variety of engineering and scientific applications, including color processing \cite{hardy1952}. Most low-pass filters are weighted average operations \cite{oppenheim2014,mersereau1984}. In curves and surfaces, we often use weighted average forms, such as B-splines \cite{bartels1995,farin2014,cohen2001}. The illumination of multiple lights is implemented as a sum or a weighted average of functions. Another set of operations commonly used in engineering and scientific applications includes maximum and minimum, which turn into Boolean operations on binary variables \cite{givant2008}. They are used as nonlinear filters in the form of dilation and erosion \cite{maragos1987a,maragos1987b}, and in implicit surfaces to obtain union and intersection \cite{bloomenthal1997}. 

{\it H\"{o}lder average} (or {\it generalized average} or {\it power average}) is a family of functions for aggregating sets of numbers, including the common arithmetic mean, geometric mean, harmonic mean, maximum, and minimum. It is the one that unifies the common linear and nonlinear averaging operations as special cases of a single-parameter family of nonlinear operators that behave like linear operations. Unfortunately, the H\"{o}lder average is not naturally associative. On the other hand, associativity is critical to obtain consistent results and to develop computationally efficient algorithms. For example, the popularity of $4\times4$ matrices comes from the fact that they provide an associative algebra to handle linear, affine, and projective operations. Similarly, associative over operation is more useful than non-associative over operation \cite{porter1984}. Associativity guarantees to make the computation of results independent of the order of the operations. 

In this paper, we show that a generalized version of the Minkowski $L_p$ norm over complex time functions in projective space provides an associative version H\"{o}lder average,  which we, therefore,  call {\it H\"{o}lder-Minkowski Colors}. We focus mainly on the theoretical properties of this family of operations and the structure of the proposed H\"{o}lder-Minkowski operations. Furthermore, we demonstrate that H\"{o}lder-Minkowski operations are particularly useful for color applications, which involve the combination of colors. Finally, we point out that having such a nonlinear family of operations is useful for a wide variety of applications, even beyond colors. Using the parameter as a search space, it can be possible to identify operations that are better for given applications.  

\subsection{Inspiration \& Rationale}

Our inspiration for this paper came from the emergence of complex algebra as a consequence of problem-solving. It is helpful to present a historical perspective \cite{crowe1994}. As is well known, the concept of zero first appeared in India at the end of the seventh century. It reached western Europe around the 12th century by Fibonacci's Liber Abaci \cite{mcclenon1919}.  The inclusion of zero also enforces the introduction of negative numbers as inverse elements. This requirement was first observed by Cardano in the 16th century \cite{heeffer2007}. Cardano and Bombelli also realized the need for complex numbers to have a complete algebra as a by-product of Tartaglia's method for computing cubic roots \cite{heeffer2008}. The full realization of the need took another 200 years until the popularization of complex algebra by Argand and Gauss \cite{green1976}. This popularization also led to the sudden emergence of other algebras, such as Quaternion, vector, matrix, and Clifford (or Geometric) algebras, in the 19th century. See \cite{crowe1994} for details. Using only positive numbers appears natural, but it always complicates dealing with inverse problems. 

This work was initially motivated by an obscure problem that can cause aliasing in rendering. The derivative of the direct illumination term, $max(\cos \theta, 0)$,  which we frequently use in rendering processes, is discontinuous at zero. This discontinuity of the derivative causes visual artifacts, as shown in Figure~\ref{fig_images/hedron/2}. One can argue that this problem can be solved by using a gamma correction such as $max(\cos \theta^\gamma,0)$, where $\gamma \geq 2$. However, this approach is a hack that does not provide a formal solution. For example, it is not possible to update the result once the computation is complete. To provide a solution, we reformulated the $gamma$ correction approach into a set of operations that can be considered as H\"{o}lder and Minkowski operations. This extension eliminates the discontinuity of the derivative and provides the desired results by removing visual artifacts as shown in Figures~\ref{fig_teaser} and ~\ref{fig_teaser2}.
However, this does not solve the problem completely. Graphics rendering works with positive numbers. Therefore, once the rendering is completed, it will still not be possible to update the results, such as removing or moving some of the lights. To achieve such a result, we need to recompute the whole scene again. This practical requirement forces us to introduce the concept of ``negative lights,'', which is usually done in computer graphics practice. As a result, if we treat lights as elements, then operations on light elements demand an \textit{algebraic structure} over the light elements, which operations can also be applied on ``negative elements.'' We point out that such applications appear not only in the practice of rendering but also in other areas such as image processing and digital compositing. 

This extension is not as straightforward as simply introducing negative lights. Since our approach is based on the H\"{o}lder and Minkowski operations, which are nonlinear, we also need to deal with the inverse of these nonlinear operations, which require the domain to cover the entire gamut of complex numbers. However, this extension turned out to be consistent with the interpretation of colors as waves instead of energy. The extension can even be handled using any existing image file format, since complex numbers can be created using two positive numbers, amplitude $r$ and phase $\theta$ in two forms of $re^{i\theta}$. Current color representations only provide a single positive number for each channel, which corresponds to the amplitude, i.e. $r$ term. In our new representation, the phase term $\theta$ can be provided simply as an additional term using an additional image. For instance, instead of a single png or jpg, we need two png or jpg images. For the projective term, we can simply use the existing $\alpha$ channels. It is also possible to add two more images to provide a separate projective term per channel.

To handle colors using our new representation, operations on colors need to have efficient implementations. Two properties of the efficiency of operations are the commutativity and the associativity of the operations. The original H\"{o}lder and Minkowski operations are already commutative. Our extensions preserve the original commutativity property. On the other hand, associativity does not hold in the original H\"{o}lder operations. To provide associative versions of original H\"{o}lder operations, there is also a need to represent colors using projective geometry. These extensions are sufficient to provide an algebraic structure, which removes derivative discontinuity, provides inverse operations on lights, and has efficient implementations of its operations. 

\section{An Overview of The H\"{o}lder-Minkowski Colors}
\label{sec_Overview}

In this section, we will provide a quick overview of H\"{o}lder-Minkowski Colors. In these color structures, for every color channel, we have a distinct set of elements. We do not limit the number of channels. However, without loss of generality, these channels can be considered red, green, and blue channels. The following is the list of elements H\"{o}lder-Minkowski colors for a given channel: 

\begin{itemize}
\item $x$, $y$, $z$ or $a$, $b$, $c$: scalars. These can be real or complex numbers.
\item $\mathbf{v}=(x,a)$: the general H\"{o}lder-Minkowski element that can represent ``lights''.  
\item $(x,a)$:  a point (or H\"{o}lder element), where $a \neq 0 $ with $(x,a)=(x/a,1)$. 
\item $(x,0)$: a vector (or Minkowski element). It is used to represent singularities, such as point light, at its center.  
\item $(0,a)$: a special type of point that can represent ``materials''.
\item $(0,0)$: the identity element of  H\"{o}lder-Minkowski Colors. 
\end{itemize}

The following is a list of operations that can manipulate these elements. 

\begin{itemize}
\item $b (x,a) = (bx, ba)$: Multiplication with a scalar, where $b$ can represent ``materials''.
\item $H_p(\mathbf{v}_0,\mathbf{v}_1) = ((x_0^p+x_1^p)^{1/p}, (a_0^p+a_1^p)^{1/p})$: H\"{o}lder-Minkowski operation. This becomes a generalized Minkowski over vector elements and a generalized H\"{o}lder type average over point elements. 
\item $(x,a)^{-1,p} = ((-1)^{1/p}x, (-1)^{1/p}a ) $: H\"{o}lder-Minkowski inverse.
\end{itemize}

Using these operations and elements, there are a wide number of ways to create individual structures. The following is the list of ways to create individual structures: 
\begin{itemize}
\item Scalar Types: Real Numbers or Complex Numbers. 
\item Element Types: Only vector types of elements or both vector and point types of elements. \item Operation type: The $p$ value. 
\end{itemize}

We provide examples of structures and discuss the differences in their representational powers in Section~\ref{sec_Algebras}. 
A significant part of all these structures is the existence of inverse elements. The inverse elements depend on $p$ and are not necessarily unique. We discuss it in Section~\ref{sec_minkinverse}. As an example, we list all types of operations related to t$p=1$ case, since it is widely used in computer graphics. 

\begin{itemize}
\item $H_1((x_0,0),(x_1,0)) = ((x_0+x_1), 0)$: vector addition. The result is a vector. 
\item $H_1((x_0,a_0),(x_1,0)) = ((x_0+x_1), a_0)$: Adding a point to the vector. The result is a point. 
\item $H_1((x_0,a_0),(x_1,a_1)) = ((x_0+x_1), (a_0+a_1))$: weighted average. The result is a point. 
\item $H_1(a_0(x_0,1),a_1(x_1,1)) = ((a_0x_0+a_0x_1), (a_0+a_1))$: another weighted average. The result is again a point. 
\end{itemize}

We want to point out that only using $H_1$ on vectors by constraining scalars into real numbers can we obtain a formal algebraic structure with inverse and identity elements that are also consistent with most common rendering models. Note that in most rendering types of applications, there are forward problems, in which there is no need to make updates or solve problems. In such cases, the four operations above can also be used with only positive real numbers. This creates the illusion that positive numbers are sufficient. However, once we need to solve a problem, we need inverse operations. For the case $p=1$, having inverse elements requires only negative numbers. Complex numbers, on the other hand, are needed for most cases where $p \neq 1$. 

\section{Our Contributions}

This introduction of the paper already provides a basic overview of how to use real or complex numbers in projective space using addition and average operators. We think that even the brief overview in Section~\ref{sec_Overview} could be sufficient to extend the color models for many graduate students in computer graphics. In the rest of the paper, we will provide the need for these algebraic structures, the theorems demonstrating that they are associative and commutative, examples of algebraic structures, and a few examples of practical power coming from using different structures for different problems. Our contribution to the rest of the paper includes the following: 
\begin{enumerate} 
\item We have developed the generalized Minkowski operation based on the $L_p$ norm over complex numbers. 
\item Based on our generalized Minkowski operation in the projective space, we have developed a single parameter family of binary operations that becomes an associative version of H\"{o}lder mean, which can be used in a wide variety of graphics applications and beyond. 
\item The new operations are not the same as H\"{o}lder mean in the weighted average case. Therefore, we have mathematically proven that this new family of operations also provides averaging of positive real numbers.  
\item We have mathematically proven that this family of binary operations is commutative and associative in a homogeneous space.
\item We have shown that the inverse always exists. This property makes it easy to update computations. 
\item We have also proven that they are scale-invariant.
\item We have provided a qualitatively similar continuous function for shadows to avoid infinite frequencies caused by discontinuities.
\end{enumerate}

In the rest of the paper, we first provide the mathematical and physical needs in Sections~\ref{sec_motivation} and ~\ref{sec_Physical}. We then provide formal proofs to demonstrate that these concepts can work for a larger set of operations by forming more general algebraic structures. We also provide examples of H\"{o}lder-Minkowski Colors in Section~\ref{sec_Algebras} and some of their applications in Section~~\ref{sec_Applications}. The rest of the paper is almost no visual, which is very unusual for a computer graphics paper. That is a conscious choice for us. We focus only on the theory to demonstrate the mathematical validity of this algebraic framework. Nevertheless, we have provided a few examples that cannot be done with existing methods. 

 \begin{figure*}
    \centering  
    \begin{subfigure}[t]{0.32\textwidth}
        \fbox{\includegraphics[width=1.0\textwidth]{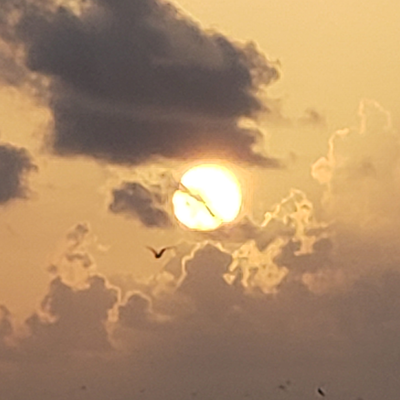}}
        \caption{$r_0$ amplitude image. }
        \label{fig_process2/00}
    \end{subfigure}
          \hfill
    \begin{subfigure}[t]{0.32\textwidth}
        \fbox{\includegraphics[width=1.0\textwidth]{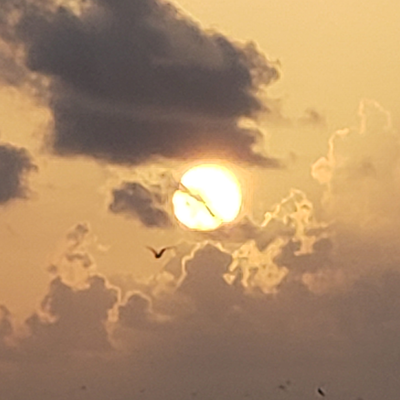}}
        \caption{$\theta_0$ phase image.}
        \label{fig_process2/01}
    \end{subfigure}
          \hfill
    \begin{subfigure}[t]{0.32\textwidth}
        \fbox{\includegraphics[width=1.0\textwidth]{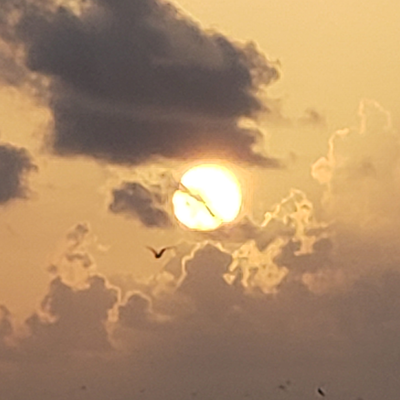}}
        \caption{$r_1$: another amplitude image.}
        \label{fig_process2/10}
    \end{subfigure}
         \hfill
    \begin{subfigure}[t]{0.32\textwidth}
        \fbox{\includegraphics[width=1.0\textwidth]{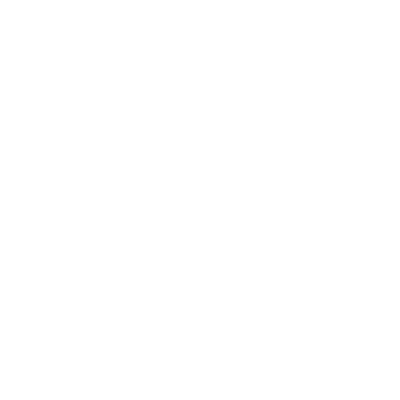}}
        \caption{$\theta_1=(1,1,1)$: white phase image. }
        \label{fig_process2/11}
    \end{subfigure}
      \hfill
    \begin{subfigure}[t]{0.32\textwidth}
        \fbox{\includegraphics[width=1.0\textwidth]{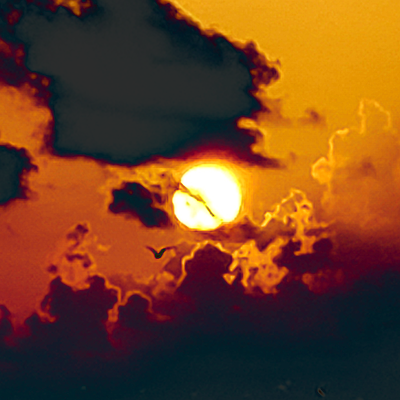}}
        \caption{Amplitude of \\$0.5(r_0e^{i\theta_0}+r_1e^{i\theta_1})$ }
        \label{fig_process2/final}
    \end{subfigure}
\caption{An example that demonstrates how to create complex representations using multiple images. This example also demonstrates the impact of the phase term. Note that each image consists of three channels, red, green, and blue. We assume that in phase images the color value $c$ corresponds to $2c \pi$. For example, the image of the white phase corresponds to $(2\pi, 2\pi, 2\pi )$. }
\label{fig_process1}
\end{figure*}

\section{Mathematical Motivations}\label{sec_motivation}

\subsection{Motivation for Inverse} 

Our motivation for this paper came from purely mathematical constraints that existed in the color domain. In current practice, we restrict the colors to a positive range $[0,max]$ and do not allow inverse elements. Without loss of generality, we can assume that colors are positive real numbers in $\Re^+$. Note that positive real numbers include any Low Dynamic Range (LDR) or High Dynamic Range (HDR) case since they can be represented as a subset of $\Re^+$, e.g. a subset that consists of only rational numbers in the form of $i/2^k$ where $i=0,1,...,2^k-1$. For the Low Dynamic Range case, we usually choose $k=8$. High dynamic range (HDR) approaches improved these color representations, allowing us to significantly increase the range \cite{debevec2008}. However, HDR colors are still essentially positive numbers \cite{reinhard2010,banterle2017}. Constraining our domain of colors into positive numbers creates a major issue in solving problems in the color domain even using only linear operations. 

Let us now demonstrate this issue using some examples. First, consider a simple algebraic problem such as 
$$r+0.5=0.75$$
even though the result is simply  
$$r=0.25$$
we cannot solve this using standard algebraic tools in the color space, since standard algebraic tools require adding the inverse of $+0.5$. Unfortunately, the inverse of $+0.5$ does not exist in color space since the inverse, which is $-0.5$, cannot be considered formally a color. This can be solved with an ad hoc approach, by allowing subtraction but not additional inverse. 

We need to open parentheses here and point out that the inverse of operations, such as subtraction, should not be used at all. These inverse operations are created to demonstrate concepts to elementary school children. Inverse operators are usually not associative or commutative. On the other hand, addition with an inverse element is both associative and commutative. 

It is, therefore, that this ad hoc approach, i.e. allowing the inverse of the addition operator, is not sufficient if the problem requires solving a chain of equations. Now, consider the following set, in which we need to compute the value of $r_1$ using an intermediate auxiliary color $r_0$ as follows,  
$$r_0+1.0=0.5$$
$$r_1=0.5+r_0$$
Since we define $r_0$ as a color, it also must be positive. In most applications such as Photoshop, $r_0=-0.5$ is turned into either $r_0=0$ or $r_0=0.5$ depending on the negation operation.  Then, the value of $r_1$ is calculated from $r_1=0.5$ or $r_1=1.0$ instead of $r_1=0$. Note that this issue is not something that can be solved by using HDR. The domain already includes any dynamic range we want to achieve. 

The problem becomes even bigger when we allow nonlinear operations such as generalized (or H\"{o}lder average) or Minkowski operations. The inverse elements according to such non-linear operations can be complex numbers. As we have demonstrated negative numbers appear during the solution process involving linear equations, and complex numbers appear during the solution process of nonlinear equations. By not allowing them, we guarantee that we will get the wrong results even for simple problems. Unfortunately, we cannot contain this problem, since we cannot control what will happen during the computation process. 

\subsection{Motivation for Associativity \& Commutative} 

Associativity is critical for obtaining consistent results, as it guarantees that the results are independent of the order of the operations. Associativity is also critical to develop computationally efficient algorithms. For example, the popularity of $4\times4$ matrices in homogeneous space comes from the fact that this framework provides an associative algebra to handle linear, affine, and projective operations \cite{roberts1963}. This property speeds up the computation of multiple operations by concatenating all operations into a single $4\times4$ matrix. Undo operation also becomes easy, multiplication of inverse matrix. Similarly, associative over operation is popular and more useful than non-associative over operation \cite{porter1984}. The associative alpha operation is associative, but not commutative. We cannot implement operations in random order. Non-commutativity makes sense for the specific use of over-operation. However, in our case, we also want commutativity, since we do not want order dependency. We observe that the order of the light arrival in rendering to any surface must not be important. Using this analogy, we allow only commutative operations even in compositing operations. The advantage of commutativity is that we do not have to keep a stack of operations to make undo as in the $4\times4$ matrix operations. Undo can be implemented in any order. We can remove and add any effect in any order. 

\section{Physical Connections}\label{sec_Physical}

The need for complex numbers in this paper initially emerged from mathematical requirements, as we have explained in the previous section. Later, we realized that complex numbers are also meaningful in handling physical problems since complex numbers help to reformulate light as a wave. We also realized that current practice already conceptually allows for wave formulation. However, since only positive numbers are used, it cannot produce waves. Our formalization makes this connection clear, and it allows one to produce waves without making significant changes in algorithms except by replacing positive real numbers with complex numbers. 

\subsection{Waves and Complex Functions}

Now, let us consider the current view of color. In current practice, we consider color as an energy. Our rendering algorithms reconstruct a set of positive real numbers for any given position and direction, each representing the ``energy'' of different wavelengths of visible light, which we call ``a color''. 
For electromagnetic radiation, the multiplication of wavelength and frequency is constant, and it is the speed of light. Therefore, we can replace the wavelength with its corresponding frequency. We can then view the current representation of color as a "simplified" Fourier transformation of a periodic function of time that is given for any given point and direction,
$$F_{\mathbf{p},\vec{n}}(t) = \sum_{\omega=\omega_{min}}^{\omega_{max}} a_{\mathbf{p},\vec{n},\omega}e^{i\omega t}$$
where $\omega$ is the frequency that corresponds any particular wavelength and $a_{\mathbf{p},\vec{n},\omega}$ is just a positive real number, $\mathbf{p}$ is a 3D position, and $\vec{n}$ is the direction given by a unit vector. Since the Fourier transform is a linear operator, all common operations such as addition or average over such functions turn into operations over coefficients. Therefore, we can operate on the coefficients $a_{\mathbf{p},\vec{n},\omega}$ only. 

One problem with this formulation is that although all functions $F_{\mathbf{p},\vec{n}}(t)$ together can be considered to form some type of wave equation, they are not waves. All points and directions work in sync as if the speed of the wave were infinite. This is just a theoretical problem, not a computational problem. If we want a representation that can provide waves, we need to introduce a delay term that could be different for each point. Note that the term $a_{\mathbf{p},\vec{n},\omega}$ can be viewed as the magnitude of a complex number. Then, the delay can be introduced by a phase term that replaces the term $a_{\mathbf{p},\vec{n},\omega}$ with complex numbers
$$c_{\mathbf{p},\vec{n},\omega}=a_{\mathbf{p},\vec{n},\omega}e^{i \theta_{\mathbf{p},\vec{n},\omega}}$$
as follows 
$$F_{\mathbf{p},\vec{n}}(t)  = \sum_{\omega=\omega_{min}}^{\omega_{max}} c_{\mathbf{p},\vec{n},\omega}e^{i\omega t}$$
where $\omega$ is still the frequency. 

This is not a significant change in formulation. The addition or average over these functions still turns into operations over coefficients. In this paper, we also demonstrate that H\"{o}lder and Minkowski operations over functions can also be implemented as operations over coefficients. Moreover, there is no need for subtraction. It can be included as an addition with an inverse.  Using a phase term is also meaningful for representing some physical processes. For instance, in iridescence, some wavelengths with different phases can cancel each other, which can only be obtained by allowing complex numbers. 

We also want to point out that this is not just useful for rendering. Even an image or a photograph should (or could) have a phase term. The reason is that the color of any pixel in an image should also be a function representing waves. We can potentially compute the phase term. We can also include such a phase term using an additional image. Assuming that the phase term is computed and stored in the second image, the impact of this extension is shown in the examples provided in Figures~\ref{fig_process1} and~\ref{fig_process2}. Even when the amplitudes of coefficients are the same, having a phase term can significantly change the resulting amplitude image.

The phase term is of paramount importance in our approach. Sensors (such as eyes or cameras) may record only the amplitude since it is difficult to record the phase. Our brains, therefore, may be processing these data using only positive real numbers. However, our computers should not necessarily mimic the process in the brain. We know that interactions outside of our sensors occur using the full flexibility of waves. Our computational framework should mimic this flexibility. Algebraic consistency appears to follow directly from this physical flexibility.  

In conclusion, we believe that the extension to complex numbers is not based only on mathematical or computational needs. Extension to complex numbers also makes sense to represent the wave behavior of light. We claim that by viewing color as an energy, we lose important phase information that will not only be physically more meaningful, but will also be useful for problem-solving. From a mathematical point of view, this restriction is also the reason why we do not have inverse elements.

\subsection{Need for Projective Space for Physical Consistency}

The projective space was originally motivated also by the mathematical need to implement an associative version of H\"{o}lder operations. We noticed later that this also required physical consistency. In current practice, we represent both lights and materials as colors. However, as we have discussed in the previous section, the lights must be functions. If materials are also functions, the multiplication of materials with lights does not make formal sense, since the multiplication changes the frequency, which will be meaningless. Therefore, materials cannot be functions. They should only be a set of real numbers that are given as ratios between $0$ and $1$. 

To make them part of the same algebraic structure, we need a formal device that can allow both lights and materials to be used as elements of this formal algebraic structure. Projective space gives us such a formalism. Note that projective spaces are given by two elements as $(x,a)$. In our case, $x$ can be a periodic function in the form of $x=c e^{i\omega t}$ where $c$ and $a$ can be any complex number for each wavelength. This formulation also allows us to include singularities such as the singularity in the center of point light by using $a=0$ as $(c_{\omega} e^{i\omega t},0)$.

This formulation also gives two ways to define material or loss. The first way is to represent materials as just scalars that are multiplied vector (or Minkowski) or point (or H\"{o}lder) elements. Note that the multiplication of a vector element with a real scalar changes its amplitude, which corresponds to one of the most common rendering operations. The multiplication of a vector element with a complex scalar changes both its amplitude and phase, which can be used in wave applications. It is even more interesting with point (or H\"{o}lder) elements. 
Multiplication of a point or H\"{o}lder entity with any number does not change the amplitude or color. However, if we have multiple terms, it gives different weights to different weights, indirectly affecting H\"{o}lder averaging. The result is easier to grasp in this case.  In other words, the projective space not only helps to transform H\"{o}lder average into an associative operation, but also gives us a general framework in which every type of element can be consistently included. 

There is also another way to define materials that can lead to different types of operation. If the $(0,a_{\omega})$ element is an added vector element, such as the point light $(c_{\omega} e^{i\omega t},0)$, it changes its amplitude and phase. For example, a light leaving point light and travel for a length of $r$ can be calculated using the following formula: 
$$(c_{\omega} e^{i\omega t},0) + (0,r^2) = (c_{\omega} e^{i\omega t},r^2)$$
We can also introduce phase shift by choosing $a_{\omega}$ term appropriately. For materials, the value of $a_{\omega}$ can be provided for the BRDF of the material. In this case, care is needed in the calculation of the term $a_{\omega}$ to obtain the desired results. This method is complicated and requires additional investigation to have more control over the results. 

\section{Background} 

Now we are ready to discuss related work. The concepts we present in this section are already well-known. We have included them here to complete the presentation and for comparing and contrasting with our H\"{o}lder-Minkowski Colors.

\subsection{Distance Functions and Lp Norm}
\label{subsec_Lp}

The concept of a distance function has been developed to provide a formal description of measuring the distance between two points in a vector space. \emph{Any} function that
satisfies the following logical conditions for distance can be
used to measure distance in a vector space.
\begin{definition} 
Let $V$ be a vector space and $\Re^*$ be the set of all nonnegative numbers, a function $f : V \longrightarrow \Re^*$
is a distance function, if it satisfies the following conditions:
\begin{itemize}
\item $f(v) = 0$ if and only if $v = 0$;
\item $f(v) = f(-v)$;
\item $f(v1 + v2) \leq f(v1) + f(v2)$ for any $v1, v2 \in V$.
\end{itemize}
\end{definition} 
There exist various distance functions defined over
different vector spaces, such as the $L_p$ norm, the Hausdorff-Besicovitch distance and the Hamming distance. $L_p$ norm is developed based on the well-known Minkowski’s inequality
$$\left(\sum(a_i^p+b_i^p)\right)^{\frac{1}{p}} \leq \left(\sum(a_i^p) \right) ^{\frac{1}{p}}   + \sum\left( b_i^p\right) ^{\frac{1}{p}} $$
where $a_i \geq 0$, $b_i \geq 0$, $p \geq 1$, and $i=0,1,\ldots,n-1$ where $n$ is any positive integer. 
Interested readers can find proof of this inequality in many classical mathematics
textbooks (e.g., \cite{mitrinovic1970}, page 55).
In an n-dimensional vector space where $v=(x_0,x_1,\ldots,x_{n-1})$, using the Minkowski operator, the norm $L_p$ is given as
$$f(v) = \left(\sum(|x_i|^p)\right)^{\frac{1}{p}}$$
$L_p$ The norm was originally introduced by Hermann Minkowski and is also called the Minkowski distance. It has been widely used in implicit shape modeling for quite a long time to model shapes. Ricci \cite{ricci1973} extended these operators by including negative values of $p$ and providing exact and approximate set operations, which we now call Ricci operators. Barr \cite{barr1981}, independently, has developed superquadrics using Minkowski operators. Hanson \cite{hanson1988} and Akleman \cite{akleman1996interactive} also used Minkowski operators to generalize superquadrics to hyperquadrics. The concept of constant-time updateability is also developed using the associative property of Ricci operators \cite{akleman1999}. Using Ricci operators, Wyvill \cite{wyvill1999} developed Constructive Soft Geometry. Even Rvachev’s exact set operations \cite{pasko1995} are related to the $L_p$ norm. Distance functions are also useful for generating implicit field functions from various types of shape information \cite{blanc1995,crespin1996}. 

Minkowski operations are defined on nonnegative numbers. It is not possible to extend Minkowski operations to all real numbers, since we cannot avoid complex numbers if we allow negative numbers. In our H\"{o}lder-Minkowski Colors, we allow complex numbers to be closed under Minkowski operations. Distance properties such as triangular equality are, of course, meaningful for only positive real numbers, since we can order numbers trivially only in one dimension.  

\subsection{H\"{o}lder Average: Generalized Mean \& Associativity}
\label{subsec_Hm}

The generalized mean is defined as 
\begin{eqnarray*}
M(x_0,...,x_{n-1}) &=&  \left(\sum\limits_{j=0}^{n-1} a_j x_j^p\right)^{1/p}\\
\mbox{where}  && \sum\limits_{j=0}^{n-1}{a_j} = 1 \\
\end{eqnarray*}
This formula provides all the operations mentioned above. However, its binary form is not associative, even though the original Minkowski operations are associative. In this paper, we demonstrate that we can solve this problem by rewriting the equation in a projective space using the $L_p$ norm. The new formula is equivalent to the generalized mean only for equal weights. On the other hand, it still guarantees to provide average over positive numbers. 

Associativity has always been important in computer graphics, as it provides several advantages such as computational efficiency with better data management \cite{foley1996}. For example, $4\times4$ matrices were introduced in 1963 to make translation and perspective associative operations \cite{roberts1963}. Similarly, associative alpha matting has been developed to make the computation of alpha operations order independent \cite{porter1984,smith1995}, which allowed alpha compositing to be used in volume rendering \cite{drebin1988}. We recognize that the basic conceptual idea in these cases is the same as ours: (1) extend the problem to the projective space and (2) use the same operation for both $x$ and $w$.  

Similarly to Minkowski operations, our extension allows complex numbers in generalized mean. For complex numbers, we cannot obtain properties similar to "average". For instance, our H\"{o}lder Average of n-number of complex numbers may not necessarily stay inside the convex hull of the positions of all complex numbers. On the other hand, if all numbers are positive real, we obtain the convex hull property. 

\subsection{Color Representations}

Even though H\"{o}lder-Minkowski Colors can potentially be used in many applications, in this paper we focus on color applications due to their direct impact on color representations. In computer graphics, we currently use special subsets of positive real numbers to represent color as an energy in the form of $$C= [c_{min}, c_{max}] \subset  [0,\infty),$$ where $c_{min} < c_{max}$. The values of $c_{min}$ and $c_{max}$ are imposed by the constraints of the display or sensor device, and never hold the values of $0$ or $\infty$. This is a de facto constraint when dealing with color on display devices. Similarly, we use special subsets of positive real numbers to represent wavelengths 
$$ \omega= \{\omega_{0},\omega_{1},\ldots, \omega_{K-1}\} \subset [0,\infty),$$ where the values of $\omega_{0},\omega_{1},\ldots, \omega_{K-1}$ are also imposed by the constraints of display or sensory devices. For example, detectable wavelengths are different for dogs, humans, birds, or infrared cameras. The colors can then be represented as functions of discrete wavelengths $\mathbf{c}(\omega)$, such as $$ \mathbf{c}: \omega \longrightarrow C $$
This is useful because it allows us to simply ignore all device dependencies and view the color as a positive real number between $[c_{min}, c_{max}]$ for a given frequency. 

\subsection{Dynamic Range}

Two positive real numbers $c_{min}$ and $c_{max}$ directly define the dynamic range of any given display or sensory device \cite{debevec2008} as
$$d=c_{max}/c_{min}.$$ From this equation it is clear why $0$ and $\infty$ are not theoretically allowed - infinite dynamic range is not possible.  To increase the dynamic range of a given device, we have two options: (1) increase $c_{max}$, and (2) decrease $c_{min}$. It also is possible to normalize the numbers by dividing each by $c_{max}$. Using the subset $[c_{min},1]$, we can still obtain an indefinitely high dynamic range, but not an infinite range. 

This discussion also helps to ignore actual device dependencies and to view the dependencies as abstract mathematical constraints. From this discussion, it must be clear that regardless of dynamic range we still have to deal with boundaries imposed by display or sensory devices. To avoid problems caused by producing numbers larger than $c_{max}$ and smaller than $c_{min}$, our operations should produce only the numbers in this subset from the numbers in this subset. This constraint can be satisfied only by using a Barycentric framework \cite{akleman2016} since with a Barycentric algebra we can guarantee to stay in the convex hull of a given set of numbers. 

\subsection{Barycentric Algebra, Partition of Unity and Convex Hull}

Barycentric algebra emerged as a dominating algebra in particular curve and surface modeling \cite{bartels1995,farin2014,cohen2001} since in curve and surface modeling we face a similar requirement.  For instance, consider the problem of interactive curve drawing on a screen. In this case, users expect the curve to stay within the boundaries of the screen. We can guarantee that the curve stays within the boundaries with barycentric algebra by using basis functions that satisfy the \textbf{partition of unity principle}. 

Let a set of points $\mathbf{v}_0, \mathbf{v}_1,\ldots,\mathbf{v}_{m-1}$ be given, then
$$\mathbf{v} = \sum_{i=0}^{n-1} a_i v_i$$
is in convex hull that is constructed by $\mathbf{v}_0,\mathbf{v}_1,\ldots,\mathbf{v}_{m-1}$ if and only if $a_i$'s satisfy the partition of unity principle, namely, 
$$\sum_{i=0}^{n-1} a_i  = 1, \; \; \; \mbox{and} \; \; \; 0 \leq a_i $$

The partition of unity property is useful, as it provides translation, rotation, and scaling invariance. In other words, if all original points are translated, rotated, or scaled in the same way, the result of the operation $\mathbf{v}$, also translates, rotates, and scales in the same way. Barycentric algebra is not only limited to curve and surface modeling. In rendering and compositing, there are already many operations that provide barycentric properties. For instance, all blur kernels, such as Gaussian blur, satisfy the Partition of Unity property. Similarly, associative over operation satisfies partition of unity. Barycentric algebra is also used to develop shaders \cite{akleman2016} to obtain nonphotorealistic styles from charcoal \cite{du2016}, crosshatching \cite{du2017} to Painterly Portraits \cite{castaneda2017paper}. Barycentric algebra is also useful to obtain specific artists such as the contemporary Chinese artist Yang Ming-Yi  \cite{liu2015chinese},  the American Western landscape painter Edgar Payne \cite{justice2018}, wine and glass still life paintings \cite{subramanian2020painterly}, and the American contemporary Fauvist-Expressionist painter Anne Garney \cite{clifford2020}.

\subsection{Projective Alpha Colors}

The projective alpha colors are another barycentric structure for colors. The projective space for colors was first proposed by Jim Blinn \cite{blinn1994} as an extension of the associative over operations. Philip Willis later introduced the projective alpha color formally \cite{willis2006}. One issue is projective alpha colors; it does not have inverse elements. In this paper, we will show how to extend projective alpha colors into a consistent structure under the H\"{o}lder-Minkowski framework (See section~\ref{sec_Algebras}). Projective alpha colors are considered as positive real four-tuples that can be transformed by 4x4 matrices. H\"{o}lder-Minkowski Colors allow matrix transformations only for $p=1$ cases.  In particular, extensions for $p=1$ cases can still allow an associative property but must lose the commutativity property. We discuss this case in Section~\ref{subsec_matrices}. 

\subsection{Overview of Projective Space}

Working in projective space simplifies dealing with the partition of unity property \cite{foley1996} in general. In this section, we will provide an algebraic extension to classical homogeneous coordinates. Let us first introduce the classical formulation. This will help to differentiate with our extension. In classical formulation, any point in projective space using homogeneous coordinates be represented in two different forms that are defined as follows
\begin{eqnarray}
 (a \mathbf{v}, a) &=&(\mathbf{v},1) \label{eq:1}\\
 ( \mathbf{v}, a) &=&(\mathbf{v}/a,1)\label{eq:2}
\end{eqnarray}
where $a$ is a positive real number including zero. In Section~\ref{subsec_vec}, we will see the need to include zero even in the classical formulation. Equation~\ref{eq:1} is called the post-multiplication form and Equation~\ref{eq:2} is called the pre-multiplication form. We will show later in Section~\ref{sec_AoCPS} that there is no need to differentiate between these two forms. 

In the classical formulation, we usually view $\mathbf{v}$ as 2D or 3D spatial information such as $\mathbf{v}=(x,y)$ and $\mathbf{v}=(x,y,z)$. In the projective color space, we view the discrete color information as $\mathbf{v}=(r,g,b)$. The reason behind this formulation is usually to use matrices to manipulate positions. In this way, we can include translation and perspective as matrix multiplication. The formulation is also useful for manipulating color in the projective color space using matrices. For instance, using rotation, we can move the red component into the blue component \cite{willis2006}. The associativity of the matrices also made use of the projective space desirable. 

\subsection{Overview of Vectors and Points in Projective Space \label{subsec_vec}}

A not well-known property of the classical homogeneous coordinates is that they can represent points and vectors in the same form. Now consider the triple $\mathbf{v}=(x,y,z)$ with the following case in which $a$ is approaching zero as 
$$\lim_{t \rightarrow \infty} \Big(\frac{1}{t} t x , \frac{1}{t} t  y , \frac{1}{t} t  z,  \frac{1}{t} \Big) \longrightarrow (x, y, z, 0 )$$
the correspondent position in homogeneous coordinates approaches $(x, y, z, 0)$ in limit. In post-multiplication form, this could be considered meaningless at first impression since its canonical form corresponds to undefined terms of $(0 \times \infty, 0 \times \infty, 0 \times \infty, 1)$ regardless of values of $x,y$ and $z$. 

Fortunately, these types of entities are possible premultiplication forms and correspond to classical vectors. For vectors, we usually ignore the last term $0$ since it does not provide additional information. As we shall see later, this will provide a consistent framework for both positions, relative positions, colors, and relative colors. 

An important implication of this formalization through the classical projective space is that it can successfully differentiate between points and vectors without any explicit awareness. In practice, we ignore the term $a$ in the calculation. However, ignoring the term $a$ is not desirable, since in that case vectors and points appear to be the same. The problem with this ignorance is that we cannot apply the same operations to vectors and points. We can add vectors but we cannot add points. We can only find the barycentric average of points. As we have demonstrated in the following, such an awareness is not necessary in a homogeneous space. The addition is meaningful even without an algebraic structure. 

\subsection{Overview of Addition on Projective Space}

Now assume that we add two vectors $(x_0,y_0,z_0,0)$ and $(x_1,y_1,z_1,0)$, the addition is still the classical vector addition as $(x_0+x_1,y_0+y_1,z_0+z_1,0)$. However, if we add two points, given in post-multiplication form as 
$$(a_0x_0,a_0y_0,a_0z_0,a_0)=(x_0,y_o,z_0,1)\; \; \mbox{and}$$  $$(a_1x_1,a_1y_1,a_1z_1,a_1)=(x_1,y_1,z_1,1);$$ the addition is still legal and turns into a Barycentric average by producing a point between these two points as 
\begin{eqnarray*}
(a_0x_0+a_1x_1,a_0y_0+a_1y_1,a_0z_0+a_1z_1,a_0+a_1)= \\
\left(\frac{a_0x_0+a_1x_1}{a_0+a_1},\frac{a_0y_0+a_1y_1}{a_0+a_1},\frac{a_0z_0+a_1z_1}{a_0+a_1},1 \right).\\
\end{eqnarray*}
In this case, we use the post-multiplication form to demonstrate that the addition of the two points is a weighted average operation. 
We can also add a point and vector, such as $(a_0x_0,a_0y_0,a_0z_0,a_0)$ and $(x_1,y_1,z_1,0)$, the result being a point $(a_0x_0+a_0x_1,a_0y_0+y_1,z_0+z_1,a_0)$. This is expected since vectors are conceptually considered as differences between two points.  
In the last case, we used mixed forms to demonstrate the concept, but all of these examples can be shown in pre-multiplied form. As discussed earlier, we use the post-multiplication form to provide a simple demonstration of how weighted averages work.  

These are very useful properties, but unfortunately, the addition operation does not form an algebraic structure. Consider the two points as
$$(sa_0x_0,sa_0y_0,sa_0z_0,sa_0)=(x_0,y_o,z_0,1)\; \; \mbox{and}$$  $$(a_1x_1,a_1y_1,a_1z_1,a_1)=(x_1,y_1,z_1,1);$$
where $s$ is any positive number. Now if we compute addition we still obtain a weighted average as follows, 
\begin{eqnarray*}
(sa_0x_0+a_1x_1,sa_0y_0s+a_1y_1,sa_0z_0+a_1z_1,sa_0+a_1)= \\
\left(\frac{sa_0x_0+a_1x_1}{sa_0+a_1},\frac{sa_0y_0+a_1y_1}{sa_0+a_1},\frac{sa_0z_0+a_1z_1}{sa_0+a_1},1 \right).\\
\end{eqnarray*}
However, the result depends on $s$. If $s\rightarrow \infty$, it goes to $(x_0,y_0,z_0,1)$. If $1/s\rightarrow \infty$, it goes to $(x_1,y_1,z_1,1)$. Having infinitely many results does not make algebraic sense if we assume the original equality of projection space holds. There is a need for a solution to this problem to turn this addition operation into an algebraic structure. 

The second issue is that all these operations are well defined for any $\mathbf{v}$ as long as we use the same domain for all $\mathbf{v}$. On the other hand, identity and inverse elements are still not defined. In other words, despite its elegance, this formulation in homogeneous coordinates does not provide an algebraic structure. To develop such a structure like an Abelian group we also need to provide identity and inverse elements. The next section presents our solution. 

\section{Algebralization of Projective Space} 
\label{sec_AoCPS}

For algebraization of the classical projective space of 3D points and vectors, we need to reinterpret the projective space. In this section, we first introduce a multiplication operation as follows:  
$$b(x,y,z,a) = (bx,by,bz,ba)$$
for any real number $b$.
We then assume that 
$$(x,y,z,a) \neq (x/a,y/a,z/a,1)$$
for any real number $a$. Note that this also means that 
$$a(x,y,z,1)=(ax,ay,az,a) \neq (x,y,z,1)$$
Note that these new conditions remove the need for post-multiplication and pre-multiplication forms. Moreover, each 4D element is now unique, and the addition operation over two 4D elements creates another 4D element. In addition, we now have an identity element for addition, which is $(0,0,0,0)$ since 
$$(x,y,z,a)+(0,0,0,0) = (x,y,z,a)$$
and we have the inverse element of $(x,y,z,a)$ as $(-x,-y,-z,-a)$ since 
$$(x,y,z,a)+(-x,-y,-z,-a) = (0,0,0,0)$$
Another important property of this framework is that projection is already included as multiplication with multiplicative inverse of $a$ as follows: 
$$(x,y,z,1)=\frac{1}{a}(ax,ay,az,a)$$
We can view this as a projection operator as follows: 
$$(x,y,z,1)=C((ax,ay,az,a))$$
where $C((a\mathbf{v},a)) = (\mathbf{v}, 1)$ is a projection operator. 
This reinterpretation appears to be conceptually similar to homogeneous coordinates as follows: 
$$C((x,y,z,1))=C((ax,ay,az,a))=C((bx,by,bz,b))$$
In other words, the elements are not equivalent, but their projections are equivalent. 
Note that the projection of the addition is not algebraic. This is not a problem, since we will not use it to solve problems. It is used only after completing the computation to obtain the final result. Now, we are ready to extend this framework with a more general form by replacing $\mathbf{v}$ and $a$ with complex numbers and functions. 

\section{Extension to Complex Functions} 

One of our key ideas is to view $\mathbf{v}$ and $a$ as complex functions in the form of $\mathbf{v} : D \rightarrow \mathbb{C}$ where $D$ can be a discrete or continuous domain and $\mathbb{C}$ is set of complex numbers. In the explanation of the idea, we first ignore the complex space and focus on the function property. Note that if the domain is discrete, these functions can still be viewed as 2D or 3D spatial information such as $\mathbf{v}=(x,y)$ and $\mathbf{v}=(x,y,z)$, discrete color information such as $\mathbf{v}=(r,g,b)$. However, it should be noted that these are not 2D or 3D vectors in the classical sense. As we have explained earlier, this formulation consists of three periodic functions
$$r(t)=re^{i\omega_r t} $$
$$g(t)=ge^{i\omega_g t} $$
$$b(t)=be^{i\omega_b t} $$
and $\mathbf{v}$ is actually a single function that can be written as follows
$$\mathbf{v}(t)=re^{i\omega_r t}+ge^{i\omega_g t}+be^{i\omega_b t}$$
where $\omega_r$, $\omega_g$ and $\omega_b$ frequencies that correspond to related wavelengths. We treat them as three separate entities to simplify the operations. 
We do not consider a continuous domain, but note that the number of channels can be arbitrarily large and the color can be represented as a summation $$\mathbf{v}(t)=\sum_{\omega} \mathbf{c}(\omega) e^{i\omega t}$$  where the term $\omega$ is the frequency that corresponds to the related wavelength.
It is important to note that the number of channels does not change much. We still treat each channel separately. Now, let $a$ denote a discrete function of the wavelength as $\mathbf{a} = \mathbf{a}(\omega)$. In classical homogeneous coordinates, the term $a$ is a positive scalar $[0,\infty)$. In this paper, we extend $a$ to all complex numbers as $a \in \mathbb{C}$. Transformations are achieved in this form by simply dividing and multiplying color functions with $a \in \mathbb{C}$ such as $\mathbf{v}(\omega) / \mathbf{a}(\omega)$ and $\mathbf{v}(\omega) \mathbf{a}(\omega)$. 
For example, the classical $\mathbf{v}=(r,g,b)$ term can be given as follows:
$$\mathbf{v}=(r,g,b) = \left(\frac{x_r}{a_r}, \frac{x_g}{a_g}, \frac{x_b}{a_b} \right)$$
The corresponding functions are now 
$$r(t)=\frac{x_r}{a_r}e^{i\omega_r t} $$
$$g(t)=\frac{x_g}{a_g}e^{i\omega_g t} $$
$$b(t)=\frac{x_b}{a_b}e^{i\omega_b t} $$
and $\mathbf{v}$ is also a single function that can be written as follows
$$\mathbf{v}(t)=\frac{x_r}{a_r}e^{i\omega_r t}+\frac{x_g}{a_g}e^{i\omega_g t}+\frac{x_b}{a_b}e^{i\omega_b t}$$
where $\omega_r$, $\omega_g$, and $\omega_b$ are frequencies that correspond to related wavelengths. One problem with these forms is that they do not support associativity. The identity element that requires $x_r=x_g=x_b=a_b=a_g=a_b=0$ is not defined even when we include complex numbers.    

Our goal is to operate effectively on $r(t)$, $g(t)$, and $b(t)$. Our first observation is that we can obtain each term as a projection of positions in projection space as follows 
$$r(t)=C(x_r e^{i\omega_r t}, a_r) $$
$$g(t)=C(x_g e^{i\omega_r t}, a_g) $$
$$b(t)=C(x_b e^{i\omega_r t}, a_b) $$
Note that these are the final color results. Our operations do not have to be in this space. We can simply consider the following element as our basic element for any given $\omega$ as 
$$\mathbf{v}(\omega,t) = (x(\omega) e^{i\omega t}, a(\omega))$$
In the previous section, we have already shown that this space already has an identity and an inverse for addition. Since the addition of elements $(x(\omega) e^{i\omega t}, a(\omega))$ is also associative, this structure already supports all properties of the Abelian group (for a full treatment of this case, see Section~\ref{section_APAPS}). This extended form still captures the properties of the original projective space for the positive real values of $a$. The only problem is to handle the terms of $e^{i\omega t}$. Fortunately, we can ignore the terms of $e^{i\omega t}$ if we deal with only addition. It can also be ignored under the projection that gives the average. Then we can make operations directly on coefficients considering elements in the form of
$$\mathbf{v}(\omega) = (x(\omega), a(\omega)).$$
If we know the channel $\omega$, there is no need to explicitly list terms of $\omega$. As a result, we can further remove the term of $\omega$ and view elements as only two complex numbers as follows
$$\mathbf{v}=(x, a).$$
We can also ignore the terms of $e^{i\omega t}$ for the H\"{o}lder-Minkowski operations as discussed in Section~\ref{subsec_RPFT}. Now we are ready to replace the addition operation with H\"{o}lder-Minkowski operations.

\section{H\"{o}lder-Minkowski Colors}

Our key idea in this paper is the generalization of the $L_p$ norm, which we can call the generalized Minkowski operator. $L_P$ norm, which was originally invented by Minkowski, is outlined in Section~\ref{subsec_Lp}. Our generalized Minkowski operator, or H\"{o}lder-Minkowski operator as we call it, is given as follows: $$H(x_0,x_1)=(x_0^p+ x_1^p )^{\frac{1}{p}}$$
where $x_0$ and $x_1$ are complex numbers and $p \in \Re$. 
There are two differences from the classical $L_p$ norm. 
\begin{enumerate}
\item Unlike $L_p$ norm, $x_0$ and $x_1$  are not just positive real numbers. Note that for the $L_p$ norm, we usually need to apply the absolute value operator to deal with only positive real numbers. 
\item Unlike the $L_p$ norm, we allow negative values of $p$ like H\"{o}lder operators. 
\end{enumerate} 
It appears at first glance that allowing negative real numbers and zero for parameter $p$ is not meaningful from the perspective of distance functions. On the other hand, negative values of $p$ are acceptable and useful for the H\"{o}lder averaging. It is well known that the value of $p=0$ provides a geometric mean in H\"{o}lder averaging as it will be explained later. The operation also provides basic H\"{o}lder Average for equal weights. 

We now build our operations using this family of generalized Minkowski operators. Let $a_0v_0=(a_0 x_0,a_0)$ and $a_1v_1=(a_1 x_1,a_1)$ denote two color points, then our H\"{o}lder-Minkowski addition operator is defined as a binary operator as follows:
\begin{eqnarray*}
H_p(v_0,v_1) &=&  \left( \left( a_0^p x_0^p + a_1^p x_1^p \right)^\frac{1}{p}, \left(a_0^p + a_1^p\right)^\frac{1}{p}\right) \\
\end{eqnarray*}
where $x$ and $w$ are complex numbers or functions, and $p \in \Re$. We will later show that this operator is associative and commutative. Moreover, if we obtain the project result, the projection operation is given as follows.
\begin{eqnarray*}
C(H_p(v_0,v_1)) &=&  \left( \left( \frac{ a_0^p x_0^p + a_1^p x_1^p }{ a_0^p + a_1^p } \right)^\frac{1}{p}, 1 \right)
\end{eqnarray*}
Note that the project operation itself is not associative. It is useful only to obtain generalized, or H\"{o}lder, an average of $v_0$ and $v_1$. Note that for the $C$ operation, we can ignore the term $1$ since this will always be the case. Therefore, we usually write directly as
\begin{eqnarray*}
C(H_p(v_0,v_1)) &=&  \left( \frac{ a_0^p x_0^p + a_1^p x_1^p }{ a_0^p + a_1^p } \right)^\frac{1}{p}
\end{eqnarray*}
Note that in addition to H\"{o}lder-Minkowski operations, we also allow multiplication any element $v=(x,a)$ with any complex scalar $z$ as follows: 
$$z v  = z (x,a) = (x z,a z)  $$
Note that the multiplication does not change the projection. It is useful mainly to control weights for H\"{o}lder averaging after an application of an H\"{o}lder-Minkowski operation. For example, letting $v_0=(x_0,a_0)$ and $v_1=(x_1,a_1)$ denote two color points, we can then change the H\"{o}lder-Minkowski computation by using this multiplication as follows: 
\begin{eqnarray*}
H_p(z_0v_0,z_1v_1) &=&  \left( \left( z_0^p a_0^p x_0^p + z_1^p a_1^p x_1^p \right)^\frac{1}{p}, \left(z_0^p a_0^p + z_1^p a_1^p\right)^\frac{1}{p}\right) \\
C(H_p(z_0 v_0, z_1 v_1)) &=&  \left( \frac{ z_0^p a_0^p x_0^p + z_1^p a_1^p x_1^p }{ z_0^p a_0^p + z_1^p a_1^p } \right)^\frac{1}{p}
\end{eqnarray*}

\subsection{Removing Periodic Function Terms}
\label{subsec_RPFT}

Recall that $x$ terms are not scalar, but functions in the form of $x= x(\omega) e^{i \omega t}$. We will show that we can ignore periodic functions and operate only on coefficients. Note that
the two color points are actually $v_0=(x_0(\omega) e^{i \omega t},a_0)$ and $v_1=(x_1(\omega) e^{i \omega t},a_1)$. We will now show that it is possible to ignore $e^{i \omega t}$ terms. Note the following formulation:  
\begin{eqnarray*}
 \left( \left( a_0^p x_0^p(\omega) e^{i p \omega t} + a_1^p x_1^p(\omega) e^{i p\omega t}\right)^\frac{1}{p}, \left(a_0^p + a_1^p\right)^\frac{1}{p}\right)  &=&  \\
 \left( \left( a_0^p x_0^p(\omega)  + a_1^p x_1^p(\omega) \right)^\frac{1}{p} \; \; e^{\frac{i p\omega t}{p}}, \left(a_0^p + a_1^p\right)^\frac{1}{p}\right) &=&  \\
  \left( \left( a_0^p x_0^p(\omega)  + a_1^p x_1^p(\omega) \right)^\frac{1}{p} \; e^{i \omega t}, \left(a_0^p + a_1^p\right)^\frac{1}{p}\right)\\
\end{eqnarray*}
Note that the operations do not have any effect on $e^{i \omega t}$. It is also the same for the term of $C(H_p(v_0,v_1))$. This property is useful; we can simply ignore the term of $e^{i \omega t}$ and apply the operation on coefficients. In the rest of the paper, we simply ignore $e^{i \omega t}$ in our proofs. We also do not have to speak in terms of $(\omega)$, since each frequency is computed separately. In conclusion, after this point, $(x,a)$ can be considered truly as only two complex numbers. 

\subsection{A Remark on Minkowski Operations over Complex Numbers}
\label{subsec_remark}

Minkowski operations on complex numbers need to be carefully defined. 
In the rest of the paper, for any complex number $z$, we assume that the following is always correct $$(z^{p})^{\frac{1}{p}} =z.$$ 
Note that this is not really correct if we assume that the numbers are real. Consider $p=2$ and $z=-1$. Then $((-1)^{2})^{0.5}$ will most likely be computed as $1$. To avoid that, we assume the numbers are always represented in a general complex form as, for example, $-1 = e^{i (2n+1) \pi}$ and $1 = e^{i (2n) \pi}$. Now, the operation $(z^{p})^{\frac{1}{p}}$ simply returns the original number, since $p$ and $1/p$ are multiplicative inverses and cancel each other. In fact, we use this property in our proofs.

Working in a general form does not only help to obtain a single result with Minkowski operations. It is also important to obtain multiple inverses as we discuss in section~\ref{sec_minkinverse}. 

\subsection{Realizing Classical Operations}

It is known that the H\"{o}lder average provides all classical or well-known averaging operations on positive real numbers. Therefore, it is also called generalized average. Our formulation is related H\"{o}lder average. Among positive real numbers, it can be considered associative H\"{o}lder operator. It is easy to show that the H\"{o}lder-Minkowski Average is capable of providing classical operations if $x$ and $a$ are positive numbers using different values of $p$. For $p=1$, the operation becomes the classical weighted (additive) average as follows: 
$$
C(H_1(v_0,v_1)) =  \frac{ a_0 x_0 + a_1 x_1 }{ a_0 + a_1 } 
$$
For $\lim p \rightarrow 0$, the operation turns into multiplication regardless of the values of $w$ as follows: 
$$
\lim_{p \longrightarrow 0} C(H_{p}(v_0,v_1)) =  \frac{ a_0 x_0 a_1 x_1 }{ a_0 a_1  }  =  x_0 x_1
$$
 
For $p=-1$, the operation turns into a form of  harmonic mean as follows: 
$$
C(H_{-1}(v_0,v_1)) =  \left( \frac{ a_0^{-1} x_0^{-1} + a_1^{-1} x_1^{-1} }{ a_0^{-1} + a_1^{-1} } \right)^{-1} = \frac{ (a_0 + a_1)  } { a_0 x_0 + a_1 x_1 } x_0 x_1
$$
When $p \rightarrow \infty$ and $w=a_0=a_1$, the operation turns into maximum as follows: 
$$
\lim_{p \longrightarrow \infty} C(H_{p}(v_0,v_1)) = \lim_{p \longrightarrow \infty} \left( \frac{ w^{p} x_0^{p}{-1} + w^{p} x_1^{p} }{ w^{p} + w^{p} } \right)^\frac{1}{p} = max( x_0 , x_1)
$$
Similarly, when $p \rightarrow -\infty$ and $w=a_0=a_1$, the operation turns into minimum as follows: 
$$
\lim_{p \longrightarrow -\infty} C(H_{p}(v_0,v_1)) = \lim_{p \longrightarrow -\infty} \left( \frac{ w^{p} x_0^{p}{-1} + w^{p} x_1^{p} }{ w^{p} + w^{p} } \right)^\frac{1}{p} = min( x_0 , x_1)
$$
Unfortunately, H\"{o}lder-Minkowski operations do not provide a general geometric mean for $p=0$ as in H\"{o}lder average. Instead, it provides only a specific geometric mean. To demonstrate how the geometric mean is obtained H\"{o}lder-Minkowski, we provide a version of the classical proof of what is the result of the Minkowski operation at $zero$. 

\begin{theorem}
\label{themMink0}
$$\lim_{p\rightarrow 0} (a_0 x_0^p + a_1 x_1^p )^{\frac{1}{p}} \longrightarrow x_0^{\frac{a_0}{(a_0  + a_1)}} x_1^{\frac{a_1}{(a_0  + a_1)}}$$ for all positive real values of $x_0$, $x_1$, $a_0$, and $a_1$.\\
\textbf{Remark 1.} Note that here we do not take the $p$ power of $a_0$ and $a_1$. Also note that even if $a_0+a_1$ is not $1$, this formula still results in the geometric mean. 
\begin{proof}
we can rewrite the function using exponential and logarithm as follows: 

\begin{eqnarray*}
(a_0 x_0^p + a_1 x_1^p )^{\frac{1}{p}} & = &\exp\left( \ln\left((a_0 x_0^p + a_1 x_1^p )^{\frac{1}{p}}\right)\right) \\
& = &\exp\left( \frac{\ln(a_0 x_0^p + a_1 x_1^p )}{p} \right)
\end{eqnarray*}
In the limit p → 0, we can apply L'Hôpital's rule to the argument of the exponential function. Differentiating the numerator and denominator with respect to $p$, we have
\begin{eqnarray*}
\lim_{p\rightarrow 0}\left( \frac{\ln(a_0 x_0^p + a_1 x_1^p )}{p}\right)& = &
\lim_{p\rightarrow 0}\left( \frac{
\frac{a_0 x_0^p \ln(x_0) + a_1 x_1^p \ln(x_1)}{(a_0 x_0^p + a_1 x_1^p)}
}{1}\right)\\
& = &
\lim_{p\rightarrow 0}\left(
\frac{a_0 x_0^p \ln(x_0) + a_1 x_1^p \ln(x_1)}{(a_0 x_0^p + a_1 x_1^p)}
\right)\\
& = &
\left(
\frac{a_0  \ln(x_0) + a_1  \ln(x_1)}{(a_0  + a_1)}
\right)\\
& = &
\left(
t_0 \ln(x_0) + (1-t_0) \ln(x_1)
\right)\\
\end{eqnarray*}
where $$t_0 = \frac{a_0}{(a_0  + a_1)}. $$
Now, by the continuity of the exponential function, we can substitute back into the above relation to obtain
\begin{eqnarray*}
\lim_{p\rightarrow 0}\left( \frac{\ln(a_0 x_0^p + a_1 x_1^p )}{p}\right)& = &
 \exp \left(
t_0 \ln(x_0) + (1-t_0) \ln(x_1)
\right)\\
& = & x_0^{t_0} x_1^{1-t_0}\\
& = & x_0^{\frac{a_0}{(a_0  + a_1)}} x_1^{\frac{a_1}{(a_0  + a_1)}}\\
\end{eqnarray*}
This concludes the proof. 
\end{proof}
\end{theorem}

In our case, the formula turns into the following form, which is still a geometric mean, but we no longer have control over power: 
$$
\lim_{p\rightarrow 0} C(H_1(v_0,v_1)) \longrightarrow \frac{ (a_0x_0)^{1/2} a_1x_1^{1/2} }{ (a_0)^{1/2} a_1^{1/2} } = x_0^{1/2} x_1^{1/2}
$$

Next, we show that the operation guarantees to yield an average even for complex numbers if they satisfy a special property. 

\subsection{Averaging Properties on Complex Numbers}

We first need to define the averaging operation for complex numbers. The problem is that we allow complex numbers and they cannot simply be ordered. We can only use the convex hull property. Using the fact that the convex hull will be a line segment for two points, we will define the strong averaging properties. 
\begin{definition}
Let $H$ be a binary operation from $S^2$ to $ S$ and $0 \leq x_0 \leq x_1 \in S$. Then the function $H$ is a \textit{strong averaging operator} if
the following condition always holds  
$$ H(x_0,x_1) = \frac{ x_0 +t_1  x_1}{1+t_1} $$ 
for a positive real $t_1$. 
\end{definition} 
Note that addition can provide this property on complex numbers. On the other hand, this property is too strong for other H\"{o}lder-Minkowski operations on complex numbers. However, they provide a less strong averaging property. 
\begin{definition}
Without loss of generality, let $x_0 \leq x_1$. 
Let $H$ be a binary operation from $S^2$ to $ S$ and $0 \leq x_0 \leq x_1 \in S$. Then the function $H$ is a \textit{weak averaging operator} if
the following condition always holds for all $x_0$ and $x_1$ that satisfy the following equality $x_1=t_1 x_0$: where $t_1$ is a positive real number greater than $1$. Then the following inequality holds: 
$$ 1 \leq \frac{H(x_0,x_1)}{x_0} \leq t_1 .$$ 
\end{definition}
Note that a weak averaging operator works as a strong operator if we focus on only positive real numbers. In practice, all forward operations will work mostly on positive real numbers. As we shall see later, we will need negative and complex numbers to implement inverse problems. Now, we are ready to demonstrate the weak averaging properties of H\"{o}lder-Minkowski operations. 

\begin{theorem}
\label{them1}
$C(H_p(a_0 x_0,a_0),(a_1 x_1,a_1))$ is a weak averaging operator on $x_0=C(a_0 x_0,a_0)$ and $x_1=C(a_1 x_1, a_1) $ for any $x_0,x_1,a_0,a_1 \in S$ such that $x_1/x_0$ and $a_1/a_0$ are positive real numbers.

\begin{proof}
Without loss of generality, let $x_0 \leq x_1$ and $t_1=x_1/x_0$ be positive real numbers larger than $1$. We want to show that the following inequality holds
$$1 \leq \frac{ C(H_p(a_0 x_0,a_0),(a_1 x_1,a_1))}{x_0} \leq  t_1 $$ 
for all values $p \in \Re$ and $s_1=a_1/a_0$ is positive real. 
Now, let us replace $$C(H_p(a_0 x_0,a_0),(a_1 x_1,a_1))$$ with actual formulation: 
$$\frac{( a_0^p x_0^p + a_1^p x_1^p )^\frac{1}{p}}{( a_0^p + a_1^p)^\frac{1}{p} }  = 
\frac{a_0x_0(1 + s_1^p t_1^p )^\frac{1}{p}}{a_0( 1 + s_1^p)^\frac{1}{p} } =
\frac{x_0(  1 + s_1^p t_1^p )^\frac{1}{p}}{( 1 + s_1^p)^\frac{1}{p} }
$$ 
Note that all we need to show now is that
$$1 \leq \frac{( 1 + s_1^p 1_1^p )^\frac{1}{p}}{( 1 + s_1^p)^\frac{1}{p} } \leq  t_1 $$ 
Since now all the numbers are positive, we can take $p$'th power of all sides without changing the direction, then we obtain following equation: 
$$1 \leq  \frac{ 1 + s_1^p t_1^p }{ 1 + s_1^p } \leq  t_1^p $$
If we multiply all sides with $1 + 1_1^p$, we obtain the following inequality: 
$$ 1 +  s^p \;\; \leq \;\; 1 + s_1^p t_1^p \;\; \leq \;\;  t_1^p +  t_1^p s_1^p$$
Since we chose, $1 \leq t_1$, the middle part is always larger or equal to left part and it is always smaller or equal to right part regardless of the values of $p, a_0$, and $a_1$. This completes the proof. 
\end{proof}
\end{theorem}

\subsection{Commutativity and Associativity}

In this section, we will evaluate the commutativity and associativity properties of $H_p(v_0,v_1)$. Demonstrating commutativity, that is, $H_p(v_0,v_1)=H_p(v_1,v_0)$ is trivial since the operation is symmetric. On the other hand, associativity needs a little more elaboration. 
\begin{theorem}
\label{them2}
The operator $H_p$ is associative for all values of $p \in \Re$ except zero. 

\begin{proof}
Let $v_0,v_1,v_2$ denote any three color points in the projective space. If $H_p$ is an associative operator, then the following equality must hold: 
$$H_p(H_p(v_0,v_1), v_2) = H_p(v_0, H_p(v_1,v_2) $$

Let us first compute $H_p(H_p(v_0,v_1),v_2)$ as follows: 
$$
H_p\left(\left( \left( a_0^p x_0^p + a_1^p x_1^p \right)^\frac{1}{p}, \left(a_0^p + a_1^p\right)^\frac{1}{p}\right), (a_2 x_2, a_2) \right) $$
We can simply compute the first and second terms for all values of $p \in \Re - \{0\}$ as follows: 
\begin{eqnarray*}
\left(\left(\left( a_0^p x_0^p + a_1^p x_1^p \right)^\frac{1}{p}\right)^p + a_2^p x_2^p \right)^\frac{1}{p} &=&
\left( a_0^p x_0^p + a_1^p x_1^p + a_2^p x_2^p \right)^\frac{1}{p} \\
\left(\left(\left( a_0^p   + a_1^p  \right)^\frac{1}{p}\right)^p + a_1^p   \right)^\frac{1}{p} &=&
\left( a_0^p  + a_1^p   + a_2^p    \right)^\frac{1}{p} \\
\end{eqnarray*}
Therefore, 
$$H_p(H_p(v_0,v_1), v_2) = \left( \left( a_0^p x_0^p + a_1^p x_1^p + a_2^p x_2^p \right)^\frac{1}{p}, \left( a_0^p  + a_1^p   + a_2^p    \right)^\frac{1}{p} \right)  $$

Now, let us compute $H_p(v_0, H_p(v_1,v_2))$ as follows: 
$$
H_p\left(\left( \left( a_0^p x_0^p + a_1^p x_1^p \right)^\frac{1}{p}, \left(a_0^p + a_1^p\right)^\frac{1}{p}\right), (a_1 x_1, a_1) \right) $$
We can again compute the first and second terms for all values of $p \in \Re - \{0\}$ as follows: 
\begin{eqnarray*}
\left( a_0^p x_0^p \right)^\frac{1}{p} + \left(\left( a_1^p x_1^p + a_2^p x_2^p \right)^\frac{1}{p}\right)^p &=&
\left( a_0^p x_0^p + a_1^p x_1^p + a_2^p x_2^p \right)^\frac{1}{p} \\
\left( a_0^p  \right)^\frac{1}{p} + \left(\left( a_1^p  + a_2^p  \right)^\frac{1}{p}\right)^p &=&
\left( a_0^p + a_1^p  + a_2^p  \right)^\frac{1}{p} \\
\end{eqnarray*}
Therefore, 
$$H_p(v_0, H_p(v_1,v_2), ) = \left( \left( a_0^p x_0^p + a_1^p x_1^p + a_1^p x_1^p \right)^\frac{1}{p}, \left( a_0^p  + a_1^p   + a_1^p    \right)^\frac{1}{p} \right)  $$
In conclusion, for all values of $p \in \Re - \{0\}$
$$H_p(H_p(v_0,v_1), v_2) = H_p(v_0, H_p(v_1,v_2) $$

This completes the proof. 
\end{proof}
\end{theorem}

As a direct consequence of associative property, we can remove the parentheses and can write the operation as an n-ary operation where $n$ can be any integer as follows:
$$H_p(H_p(v_0,v_1), v_2) = H_p(v_0, H_p(v_1,v_2) = H_p(v_0,v_1,v_2)$$
Commutative and associative properties are useful for practical applications. For instance, these operators can be used in path tracing to update color values per pixel. It is sufficient to keep cumulative $x$ and $w$ values per pixel to update a color or a pixel. As n-ary operators, they can be formulated to be used as image filters.  

\subsection{Scale Invariant Property of Projection}

\begin{theorem}
\label{them3}
The operation $C(H_p)$ is scale invariant for all $p \in \Re$ except zero, in other words, for all $p \in \Re$ except zero \newline
$C(H_p((a_0 a x_0, a_0),(a_1 a x_1, a_1)) = a C(H_p((a_0 x_0, a_0),(a_1 x_1, a_1)).$
\begin{proof}
Let two color points are denoted by $(a_0 x_0, a_0)$ and $(a_1 x_1, a_1)$. Let $C(a_0 x_0, a_0) = x_0$ and $C(a_1 x_1, a_1) =x_1$ be scaled by $a$. Note that $a$ can be a function of $a$ as $a(a)$. With this scaling, we obtain two new colors $C(a_0 a x_0, a_0) = a x_0$ and $C(a_1 a x_1, a_1) = a x_1$. Now, 
\begin{eqnarray*}
C(H_p((a_0 a x_0, a_0),(a_1 a x_1, a_1)) 
\!&=& \!\left( \frac{ a_0^p a^p x_0^p + a_1^p a^p x_1^p }{ a_0^p + a_1^p } \right)^\frac{1}{p} \\
\!&=& \! (a^p)^\frac{1}{p} \left( \frac{ a_0^p  x_0^p + a_1^p  x_1^p }{ a_0^p + a_1^p } \right)^\frac{1}{p} \\
\!&=& \! a \left( \frac{ a_0^p  x_0^p + a_1^p  x_1^p }{ a_0^p + a_1^p } \right)^\frac{1}{p} \\ \\
\!&=& \! a C(H_p((a_0 x_0, a_0),(a_1 x_1, a_1))
\end{eqnarray*}
This completes the proof. 
\end{proof}
\end{theorem}

This means that if we multiply each color by the same amount, the final computed color will also be multiplied by the same amount. Unfortunately, the operation is not affine invariant; i.e. it is not invariant under translation. As we discussed earlier, the reason we do not consider matrix transformations on colors (as in projective alpha colors) is that the general transformation is not invariant under rotation or other matrix manipulations. On the other hand, for $p=1$, we can obtain the affine and rotation invariance, but that is already known. 

\subsection{Qualitative Similarity} 

We also demonstrate that the family of operations is qualitatively similar to the weighted average operator. 

We base our model on the concept of qualitative similarity introduced by Forbus et al.~\cite{forbus1984} to compare functions. The following definition comes directly from~\cite{forbus1984}.  

\begin{definition}
Let $F_x: \mathfrak{D} \rightarrow \mathbb{R}$ and $F_y: \mathfrak{D} \rightarrow \mathbb{R}$ denote functions $x=F_x(z)$ and $y=F_y(z)$ where $z \in \mathfrak{D}$ denote any given domain and $x, y \in \mathbb{R}$ denote real numbers. We then say that $F_x$ and $F_y$ are qualitatively proportional to each other if there exists a monotonically increasing mapping $y=f(x)$. One advantage of this formulation is that there is no need to explicitly derive the monotonically increasing function $f$. The symbol $\propto_{+}$ is used to relate two qualitatively proportional entities to each other as $$F_y \propto_{+} F_x.$$
\end{definition}

\begin{lemma}
\label{lemma1}
Any power function $F_y(z)=z^\gamma$ where $a$ is any real number except zero is qualitatively similar to
$F_x(z)=z$. 
\begin{proof}
The proof is straightforward. $F_y(z)=(F_x(z))^\gamma$, i.e. $y=f(x)=x^\gamma$, which is a monotonically increasing function for all values of $a$ except zero. Note that even zero is formally monotonically increasing but practically useless. 

\textbf{Remark:} This is the reason why $\gamma$ correction is popular in image manipulation. The $\gamma$ correction operation only creates resulting images that are qualitatively similar to the original image.
\end{proof}
\end{lemma}

It is not always easy to find a monotonically increasing function that relates the two functions. Fortunately, there is another way to demonstrate that $F_x$ and $F_y$ are qualitatively proportional by showing that the following inequality is correct for all $x_0$ and $x_1$:
$$ \frac{F_x(z_1) - F_x(z_0)}{F_y(z_1) - F_y(z_0)} \geq 0$$

If we fix all the other inputs, an operation becomes a function. Assume fixing $v_0=(a_0x_0, a_0)$ and $a_1$ in $C(H_p(v_0.v_1))$ and varying only $x_1$; we can call this function $F_{p}(z)$ where $z=v_1$. 

\begin{theorem}
\label{them4}
The functions 
$$F_{p}(z) =  \left( \frac{ a_0^p x_0^p + a_1^p z^p }{ a_0^p + a_1^p }\right)^\frac{1}{p}$$
and
$$F_{1}(z) =  \frac{ a_0x_0 + a_1 z }{ a_0 + a_1 }$$
are qualitatively similar to each other. 
\begin{proof}
We want to show that all for all $p$, $z_0$ and $z_1$:
$$ \frac{F_{p}(z_1) - F_{p}(z_0)}{F_{1}(z_1) - F_{1}(z_0)} \geq 0$$
Note that
\begin{eqnarray*}
F_{p}(z_0) &=&  \left( \frac{ a_0^p x_0^p + a_1^p z_0^p }{ a_0^p + a_1^p }\right)^\frac{1}{p}\\
F_{p}(z_1) &=&  \left( \frac{ a_0^p x_0^p + a_1^p z_1^p }{ a_0^p + a_1^p }\right)^\frac{1}{p}\\
F_{1}(z_0) &=&   \frac{ a_0 x_0 + a_1 z_0 }{ a_0 + a_1 }\\
F_{1}(z_1) &=&   \frac{ a_0 x_0 + a_1 z_1 }{ a_0 + a_1 }\\
\end{eqnarray*}

Now let us take the $p$th power of $F_p$ functions and compute the difference: 
\begin{eqnarray*}
F_{p}(z_1)^p - F_{p}(z_0)^p &=&  \frac{ a_1^p (z_1^p - z_0^p) }{  a_0^p + a_1^p }\\
F_{1}(z_1) - F_{1}(z_0) &=&  \frac{ a_1  (z_1  -  z_0)  }{ a_0 + a_1 }\\
\end{eqnarray*}

Then if we divide these two into each other, we obtain 

$$ \frac{F_{p}(z_1)^p - F_{p}(z_0)^p}{F_{1}(z_1) - F_{1}(z_0)} =A \frac{ z_1^p - z_0^p }{  z_1 - z_0}$$
where $A$ is a positive number given as
$$A = \frac{(a_0 + a_1)a_1^{p-1}}{a_0^p + a_1^p}.$$
Since $A$ is positive, we know that the rest are positive. This is because power functions are qualitatively similar to each other as already shown in Lemma~\ref{lemma1}. Note that $F_{p}(z)^p$ is also a monotonically increasing function of $F_{p}(z)$. As a result, $F_{p}(z)$ is qualitatively similar to $F_{1}(z)$
This concludes the proof. 
$$F_{p} \propto_{+} F_{1}.$$
\end{proof}
\end{theorem}
This theoretical result is important for practical applications. Assume that we have already computed an average value and new input is fed to the operation; based on this result, we know that we always guarantee that larger inputs yield larger outputs, and vice versa. 

\section{H\"{o}lder-Minkowski Identity and Inverse}
\label{sec_minkinverse}

The H\"{o}lder-Minkowski Colors include inverses by allowing $a$ and $x$ to be complex numbers. It is easy to see that $(0,0)$ is an identity element.

\begin{theorem}
\label{them5}
The term $(0,0)$ is the identity element for all H\"{o}lder-Minkowski operations. 
\begin{proof}
Let $(x,a)$ be any element of an H\"{o}lder-Minkowski Color, let us now combine it with $(0,0)$. The combination is the following:
\begin{eqnarray*}
H_p((x,a),(0,0))& =& \left( (x^p+0^p)^{1/p},(a^p+0^p)^{1/p} \right)\\
&=&\left( (x^p)^{1/p},(a^p)^{1/p} \right)\\
&=&(x,a)\\
\end{eqnarray*}
Since combining $(0,0)$ with any element using the H\"{o}lder-Minkowski operation gives the the element itself;  
this concludes the proof. 
\end{proof}
\end{theorem}

The introduction of the identity element into our set also helps us to formally introduce H\"{o}lder-Minkowski inverse.

\begin{theorem}
\label{them6}
$(x,a)^{-1_p} = \left(xe^{\frac{i (2n+1) \pi}{p}}, ae^{\frac{i(2n+1)\pi}{p}} \right) $ is H\"{o}lder-Minkowski inverse if $n$ is any positive integer, including zero. 
\begin{proof}
Let $(x,a)$ and $(y,b)=\left(xe^{\frac{i (2n+1) \pi}{p}}, ae^{\frac{i(2n+1)\pi}{p}} \right)$ be two elements of a H\"{o}lder-Minkowski Colors. We want to show that H\"{o}lder-Minkowski operations: 
$$H_p((x,a),(y,b))=(0,0)$$
Now, let us plug these two into 
\begin{eqnarray*}
H_p((x,a),(y,b))&=&  ( x^p + y^p ) ^\frac{1}{p}, ( a^p + b^p ) ^\frac{1}{p}\\
&=&  \left( x^p + \left(x e^{\frac{i (2n+1) \pi}{p}} \right)^p \right) ^\frac{1}{p}, \left( a^p + \left(a e^{\frac{i (2n+1) \pi}{p}} \right)^p \right) ^\frac{1}{p}\\
&=&  (x^p + x^p e^{i (2n+1) \pi}) ^\frac{1}{p}, ( a^p + a^p  e^{i (2n+1) \pi}) ^\frac{1}{p}\\
&=&  ( x^p - x^p ) ^\frac{1}{p}, ( a^p - a^p  ) ^\frac{1}{p}\\
&=&  (0,0)
\end{eqnarray*}
This concludes the proof. 
\end{proof}
\end{theorem}
This theorem demonstrates the need for complex numbers. Even if we allow $a$ and $x$ to be real numbers (positive and negative), inverse exists only for odd integer values $p$ where $p=2n+1$. For even integers, there is no inverse with real numbers. The number of inverse elements is $p$ for integers. For irrational numbers, there are infinitely many inverse elements. This result is useful in practice. Assume that we have a dynamic scene where an object is moving. The impact of that object can be removed from some shading points and added to other shading points. The major impact of this property, which is called constant time update-ability in implicit surfaces \cite{akleman1999}, is that we do not have to re-render the whole scene to make an update if we use this during rendering. In fact, negative lights are already being used in current practice to make changes. People even use negative lights to darken some areas. Since these approaches correspond to the $p=1$ case, our algebraic structures turn these practical solutions into formally consistent methods.

\section{Examples of H\"{o}lder-Minkowski Colors}
\label{sec_Algebras}

In this section, we provide a set of examples of H\"{o}lder-Minkowski Colors. They can be used according to the problems to be solved. 

\subsection{Algebralization of the Standard Model}

Our basic algebraic structure is the generalization of the standard model in computer graphics, which is based on the rendering equation \cite{kajiya1986rendering}. Note that the rendering equation is approximated using summation. Therefore, it can be implemented using $H_1$. For generalization, it is sufficient to add zero and negative elements by allowing all real numbers. 

\begin{itemize}
\item $x$, $y$, $z$ or $a$, $b$, $c$: real scalars including negative real numbers and zero.
\item $\mathbf{v}=(x,0)$: The Minkowski element.  
\item $\mathbf{v}=(0,0)$: The identity element.  
\item $\mathbf{v'}=(-x,0)$: The inverse element of $\mathbf{v}=(x,0)$. 
\item $H_1((x_0,0),(x_1,0)) = ((x_0+x_1), 0)$: The addition operator. The result is a vector. 
\item $b (x,0) = (bx, 0)$: Multiplication by scalar, where $b$ is material.
\end{itemize}

It is straightforward to extend the standard model into this basic algebraic structure. All we need is to add negative numbers, which is already done in practice. However, it is critical not to add subtraction. It should be implemented as adding two numbers, one of them as a complement. Note that in practice it is acceptable to ignore $0$ in $(x,0)$, and use $x$. This is not really that unusual. Even for 3D vectors, we use $(x,y,z)$, not $(x,y,z,0)$. 

Since the operation $H_1$ does not cause one to obtain complex numbers, this particular structure is also consistent. The only issue with this particular structure, it does not allow wave-type of operations since for those we need phase shifts that can be created by complex numbers. 

\subsection{Extension of the Basic Model to Complex Numbers}

Extending the basic model to handle wave-type operations is straightforward. We can simply allow complex numbers with $H_1$ to be over vectors as follows. 

\begin{itemize}
\item $x$, $y$, $z$ or $a$, $b$, $c$: complex scalars.
\item $\mathbf{v}=(x,0)$: The Minkowski element.  
\item $\mathbf{v}=(0,0)$: The identity element.  
\item $\mathbf{v'}=(-x,0)$: The inverse element of $\mathbf{v}=(x,0)$. 
\item $H_1((x_0,0),(x_1,0)) = ((x_0+x_1), 0)$: The addition operator. The result is a vector. 
\item $b (x,0) = (bx, 0)$: Multiplication by scalar, where $b$ is material.
\end{itemize}

This extension is useful to obtain additional power. In this paper, we did not provide any example, but it can be possible to get
examples that are similar to those shown in Figures~\ref{fig_process1} and ~\ref{fig_process2} . 

\subsection{Extended Barycentric Models with H\"{o}lder elements}

For barycentric models, we need to extend elements to point (or H\"{o}lder) elements. We still need to work with $H_1$ to obtain barycentric operations as follows: 

\begin{itemize}
\item $x$, $y$, $z$ or $a$, $b$, $c$: complex scalars.
\item $\mathbf{v}=(x,a)$: The H\"{o}lder or Minkowski element ($a$ can also be 0).  
\item $\mathbf{v}=(0,0)$: The identity element.  
\item $\mathbf{v'}=(-x,-a)$: The inverse element of $\mathbf{v}=(x,a)$. 
\item $H_1((x_0,a_0),(x_1,a_1)) = ((x_0+x_1), a_0+a_1)$: The addition operator.
\item $b (x,a) = (bx,ba)$: Multiplication with scalar.
\end{itemize}

As we have proven in Theorem~\ref{them1}, if $a$'s are positive real numbers, $H_1((x_0,a_0),(x_1,a_1))$ is a Barycentric operator. Since the operator is associative, we can generalize for any barycentric case. For instance, let $B_i(t)$ denote always-positive real basic functions for $i=0,1,\ldots,n$. Examples of such functions are B-spline or B\'{e}zier basis functions. Now, the following function gives us a rational parametric curve in the color space where $w_i$ can be any positive real number as
$$H_1(w_0B_0(t)(x_0,1),w_1B_1(t)(x_1,1),\ldots w_nB_n(t)(x_n,1)) $$
Similarly, we can produce rational parametric surfaces. The advantage of Barycentric formulation is that we can design any desired shader easily. The advantage of even simpler interpolation is shown in Figures~\ref{fig_process1} and ~\ref{fig_process2}. If $x$'s are complex, we can obtain unusual results in compositing two images. Now, if we replace $b$'s with complex numbers, it is possible to extrapolate the colors. 

\subsection{Generalized Structures with H\"{o}lder elements}

\begin{figure}
    \centering  
        \begin{subfigure}[t]{0.32\textwidth}
        \includegraphics[width=1.0\textwidth]{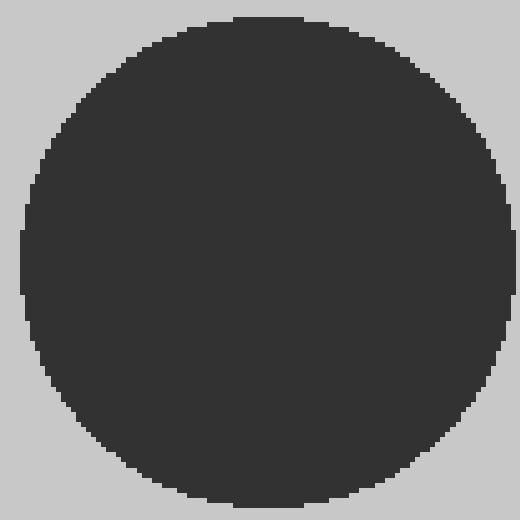}
         \includegraphics[width=1.0\textwidth]{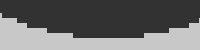}
        \caption{Original Images}
        \label{fig_circle0/circle}
    \end{subfigure}
    \hfill
    \begin{subfigure}[t]{0.32\textwidth}
        \includegraphics[width=1.0\textwidth]{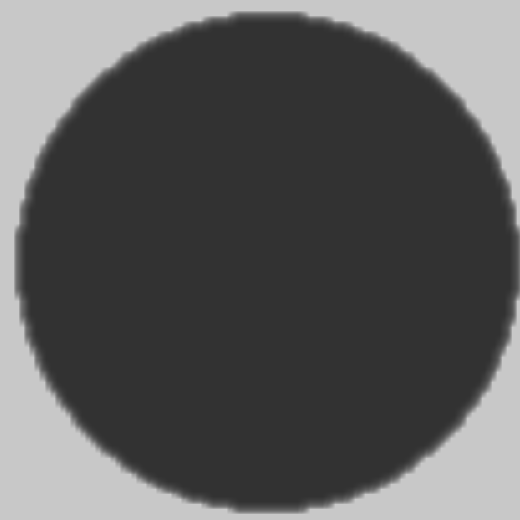}
        \includegraphics[width=1.0\textwidth]{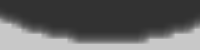}
        \caption{$p=-2$}
        \label{fig_images/filters/f-2}
    \end{subfigure}
          \hfill
    \begin{subfigure}[t]{0.32\textwidth}
        \includegraphics[width=1.0\textwidth]{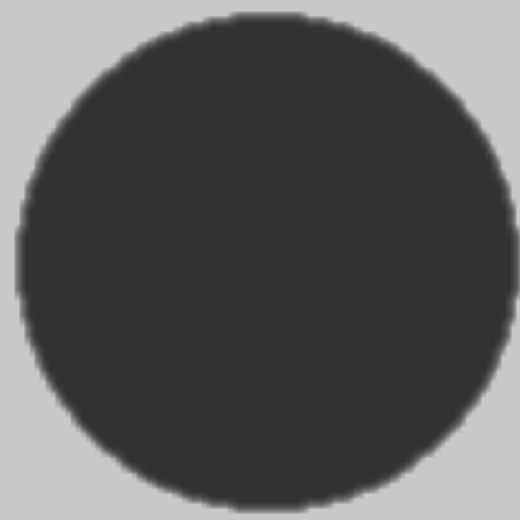}
        \includegraphics[width=1.0\textwidth]{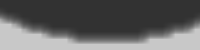}
        \caption{$p=-1$}
        \label{fig_images/filters/f-1}
    \end{subfigure}
          \hfill
        \begin{subfigure}[t]{0.32\textwidth}
        \includegraphics[width=1.0\textwidth]{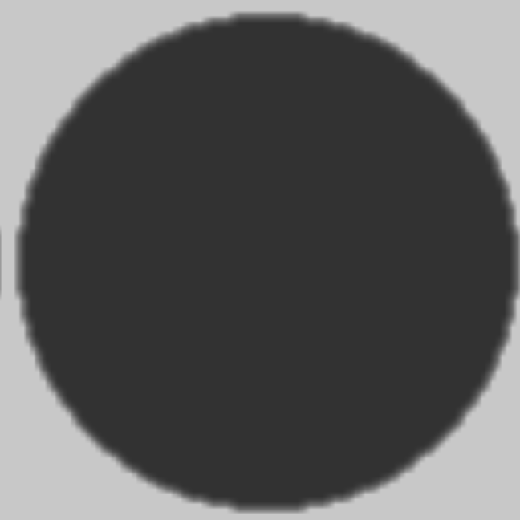}
        \includegraphics[width=1.0\textwidth]{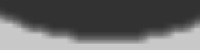}
        \caption{$p=-0.5$}
        \label{fig_images/filters/f-0.5}
    \end{subfigure}
          \hfill
        \begin{subfigure}[t]{0.32\textwidth}
        \includegraphics[width=1.0\textwidth]{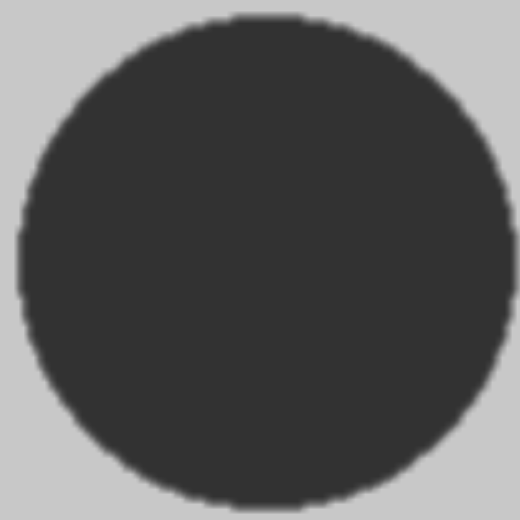}
        \includegraphics[width=1.0\textwidth]{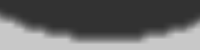}
        \caption{$p=0.5$}
        \label{fig_images/filters/f0.5}
    \end{subfigure}
          \hfill
    \begin{subfigure}[t]{0.32\textwidth}
        \includegraphics[width=1.0\textwidth]{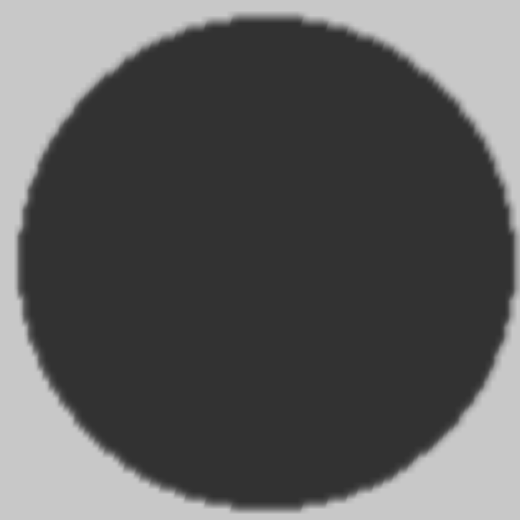}
        \includegraphics[width=1.0\textwidth]{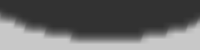}
        \caption{$p=1$}
        \label{fig_images/filters/f1}
    \end{subfigure}
\caption{An example of the effect of the value of $p$ on the aliased boundaries. Bottom images are details. In this case, we use the $15\times15$ box filter by changing the values to see the impact of the $p$ value on aliasing.  }
\label{fig_filterscircle}
\end{figure}

Generalized structures are simply extensions of barycentric forms to general $H_p$ operations. In this paper, we have provided two cases that can demonstrate the need for $H_p$ operations. The first is to remove derivative discontinuities at zero (see Figures~\ref{fig_teaser} and~\ref{fig_teaser2}) and  the second is to improve filtering operations (see Figure~\ref{fig_filterscircle}).

An interesting approach in these structures is to represent materials as H\"{o}lder elements in the form of $(0,a)$ instead of the multiplicative term $b$. This does not really fit the standard model. However, for some non-photorealistic rendering applications, it can be appropriate to use. 

We must point out that decoupling different frequencies is critical for turning general $H_p$ operations into commutative and associative operations. Because of these commutativity and associativity, it is possible to update a rendering when lights changed, moved, or removed, materials changed, or shapes moved, or removed. 

Using the same conceptual framework, it is also possible to couple different frequencies. However, this can only be done using $H_1$. It also comes with a price of losing commutativity. In the next subsection, we demonstrate the existence of matrix algebraic case that can allow coupling different frequencies. 

\subsection{On Matrix Models}
\label{subsec_matrices}

Now, assume that $b$'s are not just complex numbers, but periodic functions in the form of $\sum_{m}^{M_1}re^{iw_mt+\theta_m}$. By using these multiplication terms, we can shift frequencies to obtain other frequencies for the $H_1$ operation. Assume that we have three terms, i.e. $M=3$. Without loss of generality, we can call them red, green, and blue, which are given as three elements $(x_r,0)$, $(x_g,0)$, and $(x_b,0)$. Each term can be multiplied by periodic functions consisting of three terms as follows: 
$$(b_{r,r}+b_{r,g}+b{r,b})(x_r,0)$$
$$(b_{g,r}+b_{g,g}+b{g,b})(x_b,0)$$
$$(b_{b,r}+b_{b,g}+b{b,b})(x_g,0)$$
where $b_{g,r}$ moves green to red, for example. Note that in this case, we do not produce any frequency other than those related to red, green, and blue with the $H_1$ operation. The whole operation can be written in matrix form as follows: 

\begin{gather}
\begin{pmatrix}
(x'_r,0)\\
(x'_g,0)\\
(x'_b,0)
\end{pmatrix}
=
\begin{pmatrix}
b_{r,r} & b_{g,r} & b_{b,r}\\
b_{r,g} & b_{g,g} & b_{b,g}\\
b_{r,b} & b_{g,b} & b_{b,b}
\end{pmatrix}
\begin{pmatrix}
(x_r,0)\\
(x_g,0)\\
(x_b,0)
\end{pmatrix}
= M 
\begin{pmatrix}
(x_r,0)\\
(x_g,0)\\
(x_b,0)
\end{pmatrix}
\end{gather}
In this case, $M$ is material matrix it can also shift colors from red to green, green to blue, and so on. The number of elements can be more than three to have more control. We have simplified it to demonstrate the idea. Such an operation can particularly be useful non-photorealistic rendering application where we may want to change hues. This form is still useful since the matrices are associative. We can still concatenate material effects using the associative property. 

There are several limitations in this case. This operation is not necessarily commutative. Note that in commutative case we did not allow different channels to impact each other. Therefore, in this case, we cannot change one of the materials along a set of reflections. In our standard structures, we can change one the materials since material effects are simply multiplication with a complex or real scalar. That operation is commutative, and we can always invert it.  

Therefore, it is not surprising that we cannot use such functional multiplicative terms $b$ for the general $H_p$ case. This formulation is possible only for $H_1$, since the operation $H_1$ is linear. Another issue is that this formulation will not provide a mathematically consistent representation of projective alpha colors using $\alpha$ in terms of $a$ such as $a_r=a_b=a_g=\alpha$. It is also important to note that the operation involving the associative $\alpha$ is also different. Therefore, the only way to include $\alpha$ is by considering $\alpha$ as just another multiplicative term. We can then create a version of projective alpha colors using $4 \times 4$ matrices. 

\begin{gather}
\begin{pmatrix}
(x'_r,0)\\
(x'_g,0)\\
(x'_b,0)\\
\alpha'\\
\end{pmatrix}
=
\begin{pmatrix}
b_{r,r} & b_{g,r} & b_{b,r} & (x_{\alpha,r},0)\\
b_{r,g} & b_{g,g} & b_{b,g} & (x_{\alpha,b},0)\\
b_{r,b} & b_{g,b} & b_{b,b} & (x_{\alpha,b},0)\\
0 & 0 & 0 & b_{\alpha,\alpha}\\
\end{pmatrix}
\begin{pmatrix}
(x_r,0)\\
(x_g,0)\\
(x_b,0)\\
\alpha\\
\end{pmatrix}
\end{gather}

The first three terms in the last row must be zero, since we have not defined multiplicative inverses of Minkowski elements to obtain scalars. This also demonstrates that our projective space is conceptually different than projective alpha colors.

\begin{figure}
    \begin{subfigure}[t]{0.32\textwidth}
        \includegraphics[width=1.0\textwidth]{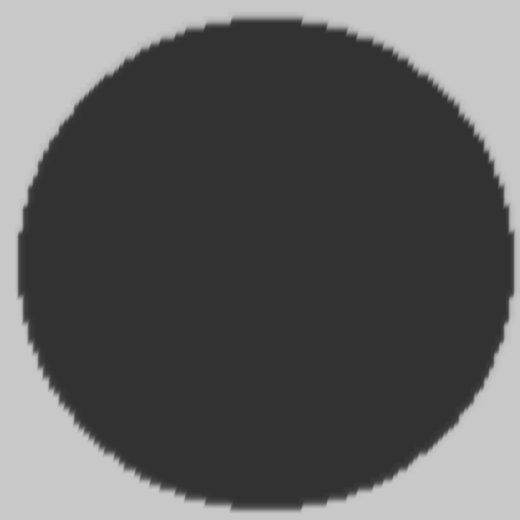} 
        \includegraphics[width=1.0\textwidth]{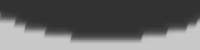}
        \caption{$p=2$}
        \label{fig_images/filters/f2}
    \end{subfigure}
        \hfill
            \begin{subfigure}[t]{0.32\textwidth}
        \includegraphics[width=1.0\textwidth]{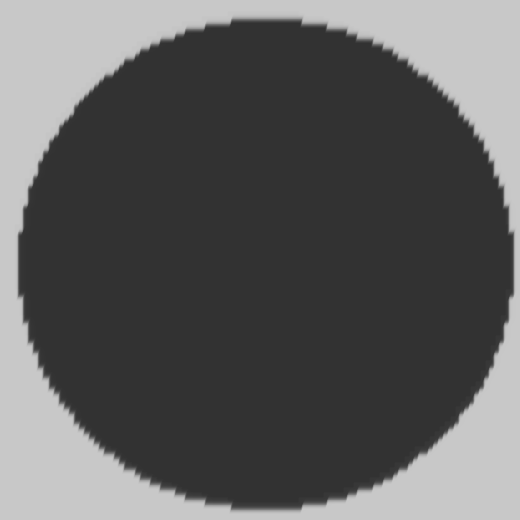}
        \includegraphics[width=1.0\textwidth]{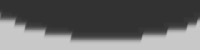}
        \caption{$p=3$}
        \label{fig_images/filters/f3}
    \end{subfigure}
          \hfill
              \begin{subfigure}[t]{0.32\textwidth}
        \includegraphics[width=1.0\textwidth]{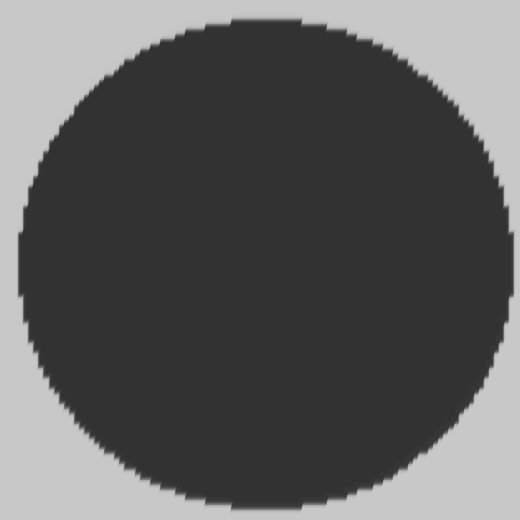}
        \includegraphics[width=1.0\textwidth]{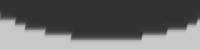}
        \caption{$p=4$}
        \label{fig_images/filters/f4}
    \end{subfigure}
          \hfill
    \begin{subfigure}[t]{0.32\textwidth}
        \includegraphics[width=1.0\textwidth]{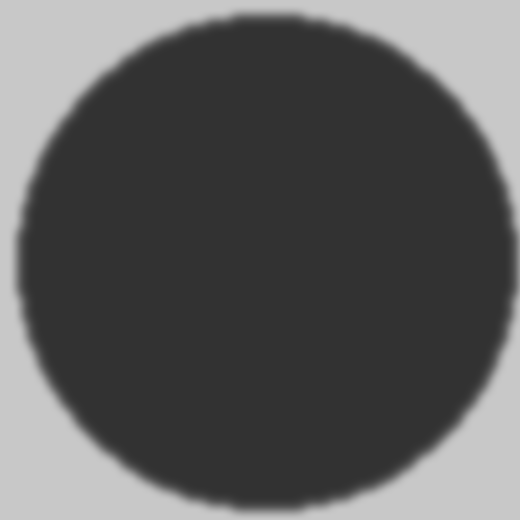}
        \includegraphics[width=1.0\textwidth]{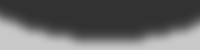}
        \caption{$p=0.5$ followed by $p=-0.5$}
        \label{fig_images/filters/f0.5-0.5}
    \end{subfigure}
    \hfill
    \begin{subfigure}[t]{0.32\textwidth}
        \includegraphics[width=1.0\textwidth]{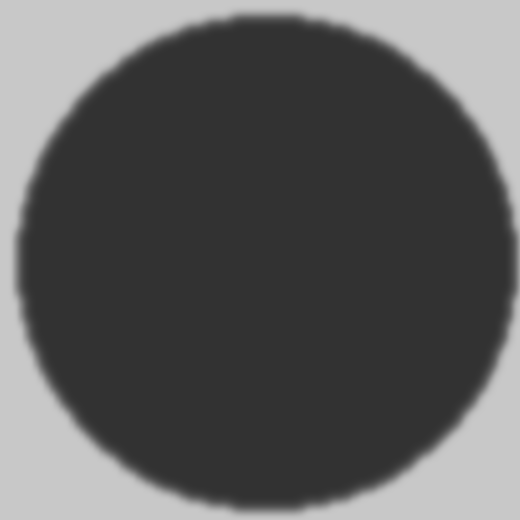}
        \includegraphics[width=1.0\textwidth]{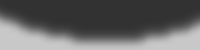}
        \caption{$p=1$ followed by $p=-1$}
        \label{fig_images/filters/f1-1}
    \end{subfigure}
    \hfill
    \begin{subfigure}[t]{0.32\textwidth}
        \includegraphics[width=1.0\textwidth]{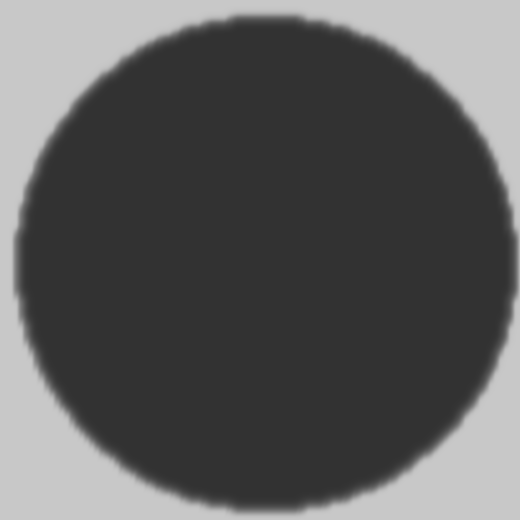}
        \includegraphics[width=1.0\textwidth]{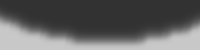}
        \caption{$p=2$ followed by $p=-2$}
        \label{fig_images/filters/f2-2}
    \end{subfigure}
    \hfill
\caption{Following up the Figure~\ref{fig_filterscircle} demonstrating more examples on the effect of the value of $p$ on the aliased boundaries. We observe that $p=1$ followed by $p-1$ gives the best results, which is similar to erosion followed by dilation.  }
\label{fig_filterscircle2}
\end{figure}

\section{Examples of Applications}
\label{sec_Applications}

There are three main applications of this set of algebraic structures: (1) Rendering, (2) Compositing, and (3) Filtering. We provided earlier a few examples of rendering and compositing that cannot be obtained without the H\"{o}lder-Minkowski framework.
In the section, we provide additional examples and discuss how to use the H\"{o}lder-Minkowski framework in filtering and composition. This section is just to demonstrate the existence of potential applications. The full treatment requires a deeper analysis. 

\subsection{Low-Pass Filtering}

Based on our generalized operations' associative property, we can reformulate
the generalized nonlinear filters as an extension of classical convolution filters as follows:

\begin{eqnarray*}
w' c'_{(x,y)} &=&  \left(\sum\limits_{i=0}^{2n}{}\sum\limits_{j=0}^{2m}{w_{(i,j)}^p c_{(x+i-n,y+j-m)}^p}\right)^{1/p}\\
w' &=&  \left(\sum\limits_{i=0}^{2n}{}\sum\limits_{j=0}^{2m}{w^p_{(i,j)}}\right)^{1/p}\\
 c'_{(x,y)} &=& \left( \frac{\sum\limits_{i=0}^{2n}{}\sum\limits_{j=0}^{2m}{w_{(i,j)}^p c_{(x+i-n,y+j-m)}^p}}{\sum\limits_{i=0}^{2n}{}\sum\limits_{j=0}^{2m}{w^p_{(i,j)}}}\right)^{1/p}
\end{eqnarray*}
where $c_{(x,y)}$ and $c'_{(x,y)}$ is the old and new colors of pixel $(x,y)$ and  $w_{(i,j)}$ are positive real numbers (kernel weights), $n$ and $m$ are two integers, and the filter size $(2n+1)\times(2m+1)$. As discussed earlier $p$ can be any real number except zero and the values of can create different filters.

\begin{itemize}
\item $p=1$: This is the classical convolution blur filter.
\item $p \rightarrow \infty$: The operation is maximum, and this is a dilation filter (or erosion depending on what is considered inside).
\item $p \rightarrow -\infty>$: The operation is minimum and this is an erosion filter.
\end{itemize}
Note that the use of very large numbers of $p$ does not make sense computationally. Instead, we can simply use the maximum and minimum. In those cases, the filters will be as follows:  

$$c'(x,y) =  \frac{max\left( w_{(i,j)} c_{(x+i-n ,y+j-m))}\right)}{ max \left( w_{(i,j)} \right)} \;\; \; \forall \; i,j $$

$$c'(x,y) =  \frac{min\left( w_{(i,j)} c_{(x+i-n ,y+j-m)}\right)}{ max \left( w_{(i,j)} \right)} \;\; \; \forall \; i,j $$

 \begin{figure*}
    \centering  
    \begin{subfigure}[t]{0.32\textwidth}
        \fbox{\includegraphics[width=1.0\textwidth]{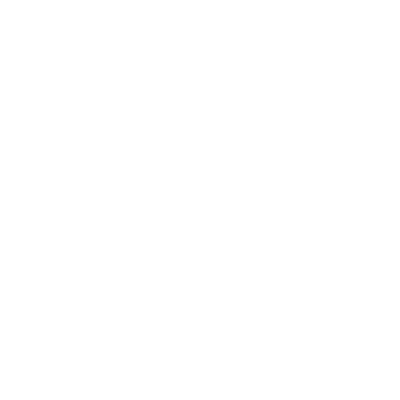}}
        \caption{$r_0=(1,1,1)$ white amplitude image. }
        \label{fig_process1/00}
    \end{subfigure}
          \hfill
    \begin{subfigure}[t]{0.32\textwidth}
        \fbox{\includegraphics[width=1.0\textwidth]{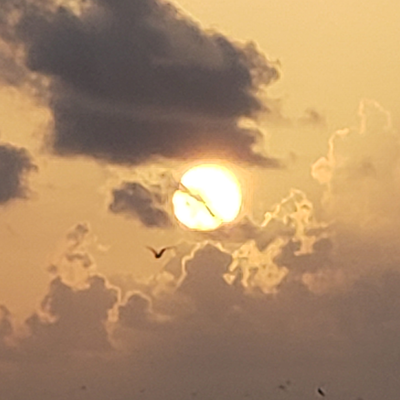}}
        \caption{$\theta_0$ phase image.}
        \label{fig_process1/01}
    \end{subfigure}
          \hfill
    \begin{subfigure}[t]{0.32\textwidth}
        \fbox{\includegraphics[width=1.0\textwidth]{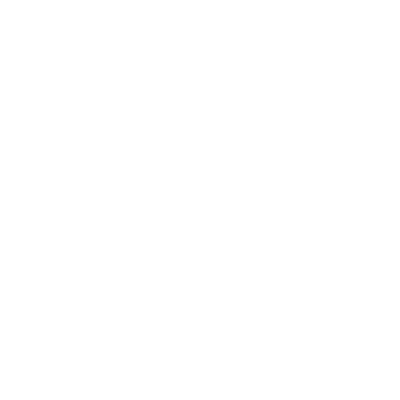}}
        \caption{$r_1=(1,1,1)$: another white amplitude image.}
        \label{fig_process1/10}
    \end{subfigure}
          \hfill
    \begin{subfigure}[t]{0.32\textwidth}
        \fbox{\includegraphics[width=1.0\textwidth]{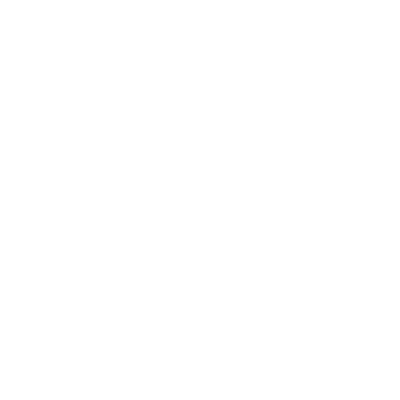}}
        \caption{$\theta_1=(1,1,1)$: white phase image. }
        \label{fig_process1/11}
    \end{subfigure}
      \hfill
    \begin{subfigure}[t]{0.32\textwidth}
        \fbox{\includegraphics[width=1.0\textwidth]{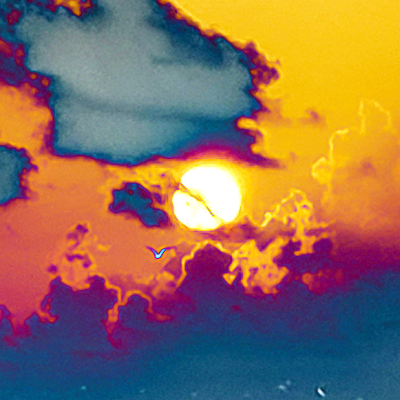}}
        \caption{Amplitude of \\$0.5(r_0e^{i\theta_0}+r_1e^{i\theta_1})$ }
        \label{fig_process1/final}
    \end{subfigure}
      \hfill
\caption{An example that demonstrates how to create complex representations using multiple images. This example also demonstrates the impact of the phase term. Note that each image consists of three channels, red, green, and blue. We assume that in phase images color value $c$ corresponds to $2c \pi$. For example, white-phase images correspond to $(2\pi, 2\pi, 2\pi )$. }
\label{fig_process2}
\end{figure*}

We do not know whether Minkowski average benefit from Fast Fourier transform except the obvious $p=1$ case. Since computation time is important in interactive applications, we have implemented H\"{o}lder-Minkowski filter in fragment shader using GLSL. We have observed that it is possible to obtain interactive real-time results. Therefore, we do not expect that the computation time poses any problem for H\"{o}lder-Minkowski filters. Therefore, H\"{o}lder-Minkowski formulations may also be useful in deep learning applications since they provide fast computation along with a single parameter to move one filter and another. In terms of visual quality, we observed that $p=1$ followed by $p-1$ provides the best results to remove aliasing in the boundaries (see Figure~\ref{fig_filterscircle}). In fact, values such as $p=0.5$ and $p=-0.5$ also give good results. This direction needs much more research. 

\subsection{Compositing}

Let $x$ and $\alpha$ denote the color and transparency of a given pixel in an image $I$. 
Then, the associative over operation\cite{porter1984} of foreground and background images, $I_1$ and $I_0$ is given as 
\begin{eqnarray*}
\alpha \; x &=& \alpha_1 \; x_1 + (1-\alpha_1) \; x_0\\
\alpha &=& \alpha_1 + (1-\alpha_1) \; \alpha_0 .
\end{eqnarray*}
where $\alpha x$ is usually called pre-multiplied color. Note that this operation is associative but not commutative. In operations we also want commutativity H\"{o}lder-Minkowski operations can provide an alternative. There can be two potential usages: (1) $\alpha$ can be used as a common $w$ for all channels; (2) an additional image can be used as color channels of $w$. 

H\"{o}lder-Minkowski operations can also be used in conjunction with generalized compositing operations where compositing is defined as  
\begin{eqnarray*}
\alpha \; x &=& \alpha_1 \; F + (1-\alpha_1) \; x_0\\
\alpha &=& \alpha_1 + (1-\alpha_1) \; \alpha_0 .
\end{eqnarray*}
where $F$ can be any function of $x_1$ and $x_0$. Some of these operations, such as maximum and minimum, produce a color, i.e., a number in a desired range. Some functions such as multiplication, i.e. $F = x_1  x_0$, produce numbers in the range of $[0,1]$ if the inputs also come from the same range $[0,1]$.  Some functions such as $F = x_1 + x_0$, on the other hand, can produce numbers outside of any range. This requires the cropping of numbers into the desired range. In applications where cropping is not desired, it makes sense to use an appropriate Minkowski average. For example, we can replace multiplication with $C(H_{-1})$ and addition with $C(H_{1})$. As we discuss next, it is also possible to extend Minkowski average to include operations such as negation with complex algebra.  Figures~\ref{fig_process1} and~\ref{fig_process2} demonstrate the effect of phase term. 

\section{Discussion}

In this paper, we introduced H\"{o}lder-Minkowski Colors that can be used in a broad range of applications. We have demonstrated that H\"{o}lder-Minkowski Colors have a wide range of properties. These operations are both commutative and associative in a projective space. They are scale-invariant and guarantee that the results of the operations are always between maximum and minimum numbers. They are also qualitatively similar to a weighted average. This theoretical framework is promising; however, there is still a need for a significant amount of additional work to identify its full power in practical applications. There are also many issues that need to be addressed. 

We have also demonstrated that H\"{o}lder-Minkowski Colors can be used with relative theoretical ease in practical applications such as rendering, compositing, and filtering. In the paper, we have devised extremely simple cases that strongly suggest that the new algebraic structures can provide additional power in rendering, compositing, and filtering. We purposely avoid complicated cases to avoid practical distractions in the establishment of theoretical foundations.   
For example, in the example provided in Figure~\ref{fig_filterscircle}, the weights and colors are real positive numbers. In other words, we did not attempt to try complex weights that can potentially provide new opportunities. However, there is a need for a theoretical understanding of how to produce these weights. Using complex coefficients, we no longer work with interpolating schemes, as visually demonstrated in Figures~\ref{fig_process1} and~\ref{fig_process2}. 

The need for an extension to complex numbers for H\"{o}lder-Minkowski inverse may initially be considered a disadvantage, since the need for H\"{o}lder-Minkowski inverse is not immediately obvious for most applications. We want to point out that the H\"{o}lder-Minkowski inverse is also useful as a representational power. For instance, even now some color operations, such as high-pass filtering, result in negative numbers. Having H\"{o}lder-Minkowski inverse makes representation complete and consistent. The discussion in section~\ref{sec_minkinverse}, provides a basic framework. However, additional theoretical development for dealing with H\"{o}lder-Minkowski inverses. One problem is that the inverse does not really exist for maximum and minimum operations. In these cases, $p$ corresponds to $\infty$ or $-\infty$, and the inverse obtained in Section~\ref{sec_minkinverse} is practically useless. Computation of the inverse will also be a problem for large values of $p$. On the other hand, having all range of values of $p$  is still useful to choose appropriate $p$ values that can work reasonably well. 

 \begin{figure*}
    \centering  
        \begin{subfigure}[t]{0.32\textwidth}
        \includegraphics[width=1.0\textwidth]{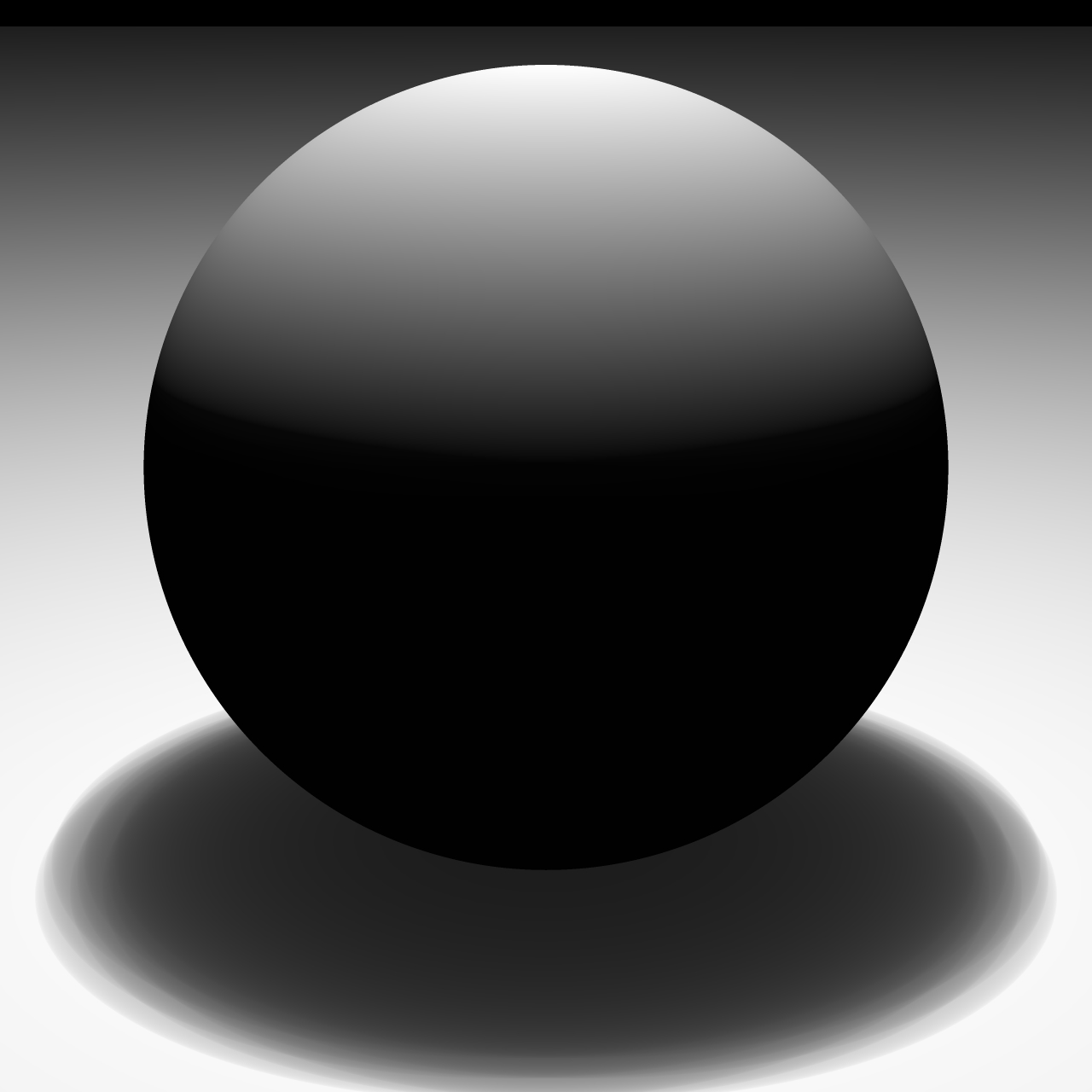}
        \includegraphics[width=1.0\textwidth]{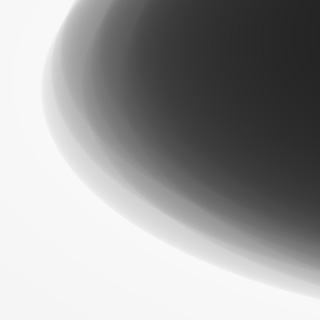}
        \caption{$p=-1.0$}
        \label{fig_sphere_single/2}
    \end{subfigure}
       \hfill 
         \begin{subfigure}[t]{0.32\textwidth}
        \includegraphics[width=1.0\textwidth]{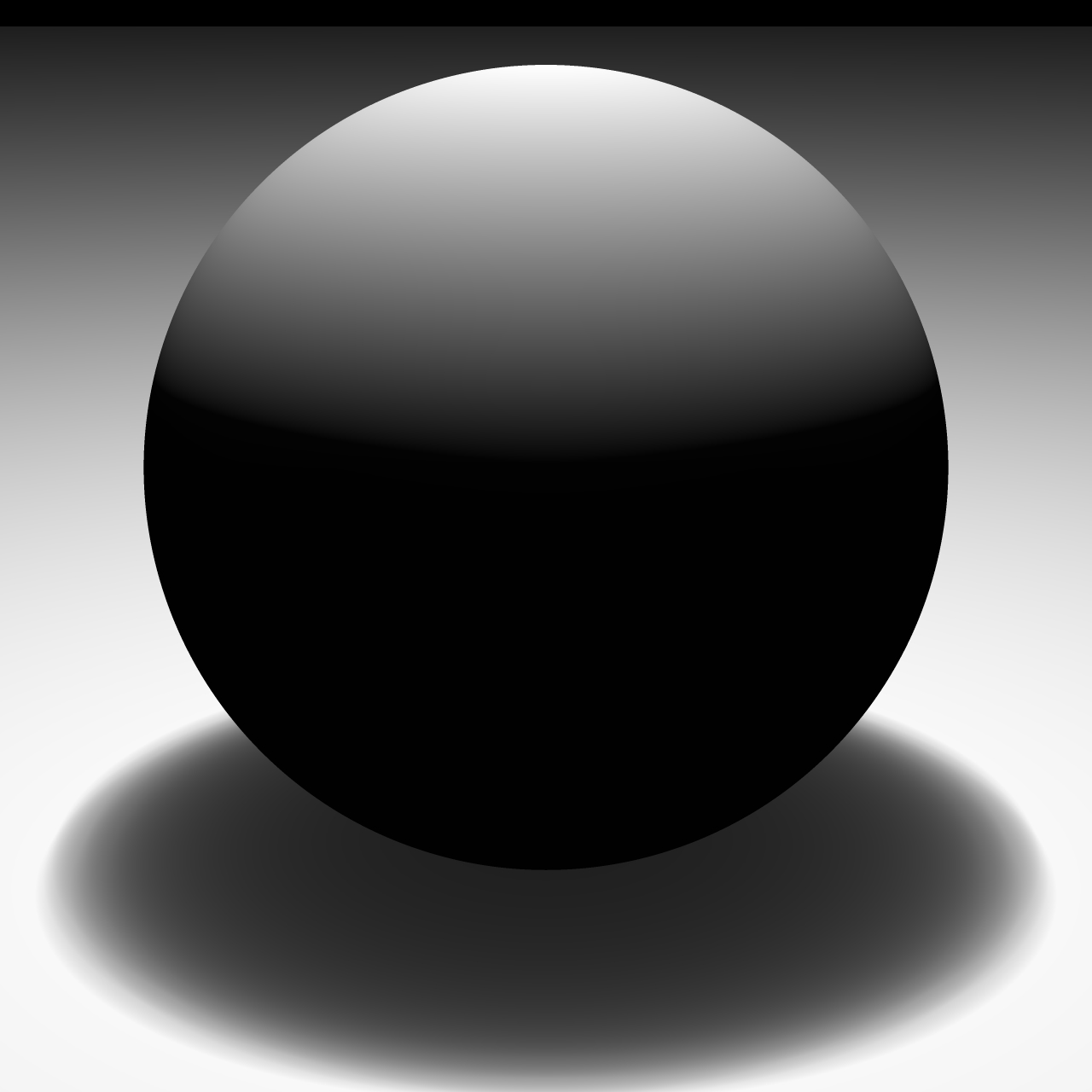}
        \includegraphics[width=1.0\textwidth]{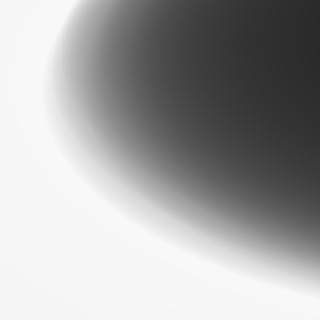}
        \caption{$p=-2.0$}
        \label{fig_sphere_single/4}
    \end{subfigure}
      \hfill 
           \begin{subfigure}[t]{0.32\textwidth}
        \includegraphics[width=1.0\textwidth]{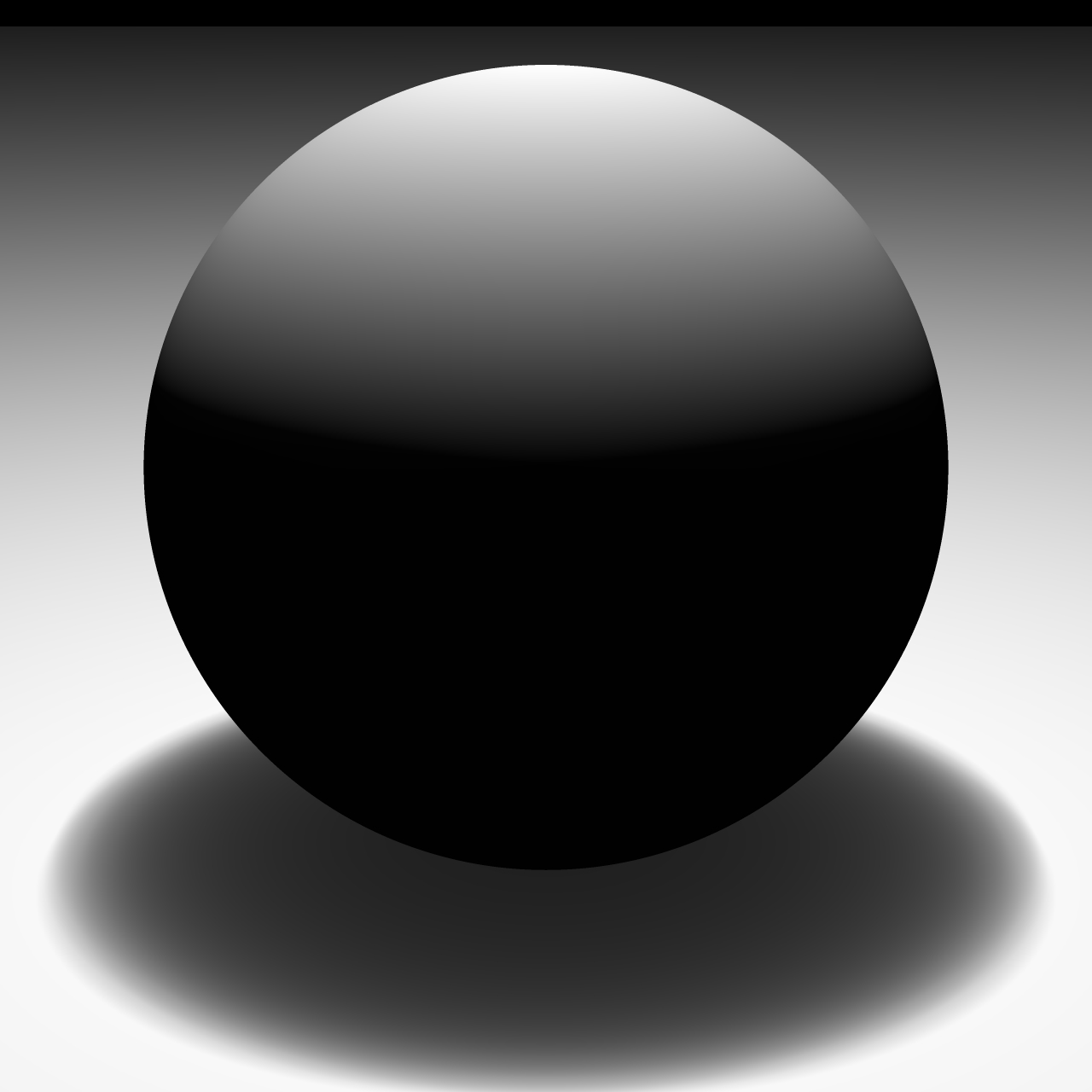}
        \includegraphics[width=1.0\textwidth]{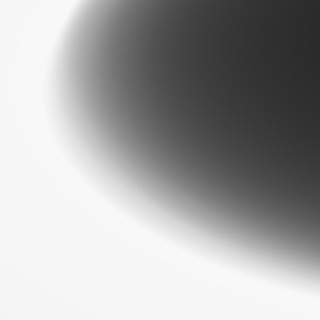}
        \caption{$p=-2.5$}
        \label{fig_sphere_single/5}
    \end{subfigure}
\caption{Another example of a sphere illuminated by nine point lights located in the vertices of a $3 \times 3$ square grid using a simplified form of equations. In this case, discontinuities exist in actual shadow functions.  The images show how the different H\"{o}lder-Minkowski cleans up discontinuities in the shadow functions. In this case, no value of $p$ corresponds exactly to the standard operations in computer graphics, since we simultaneously provide a compromise between addition and minimum. Parameter values smaller than $p=-1$ create smoother results by removing visual artifacts caused by derivative discontinuity. Very small values of $p$  approach the minimum operator and create a sharp shadow boundary}
\label{fig_teaser2a}
\end{figure*}

\begin{figure*}
    \centering  
     \begin{subfigure}[t]{0.32\textwidth}
        \includegraphics[width=1.0\textwidth]{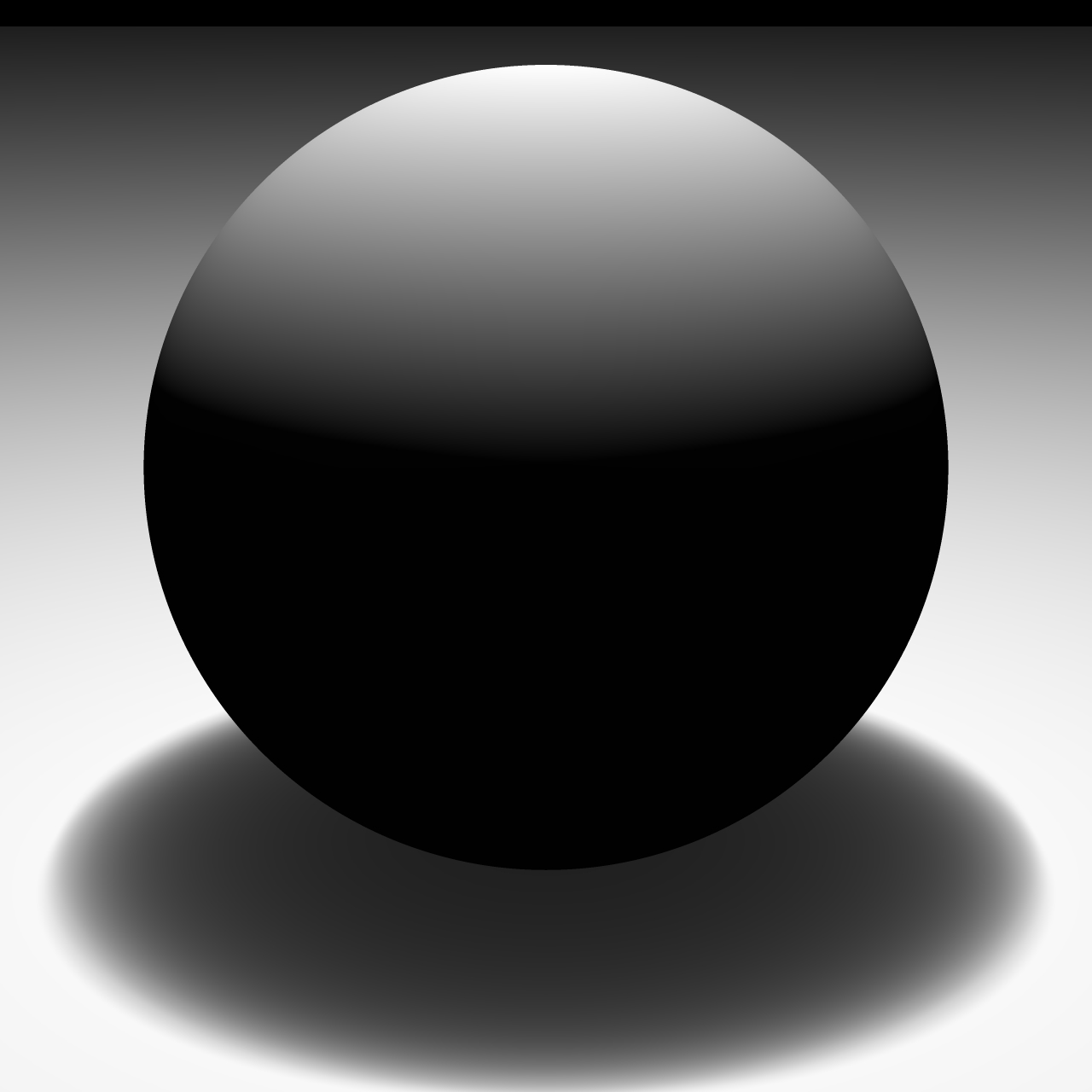}
        \includegraphics[width=1.0\textwidth]{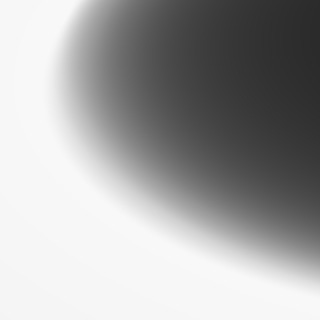}
        \caption{$p=-3.0$}
        \label{fig_sphere_single/6}
    \end{subfigure}
      \hfill 
    \begin{subfigure}[t]{0.32\textwidth}
        \includegraphics[width=1.0\textwidth]{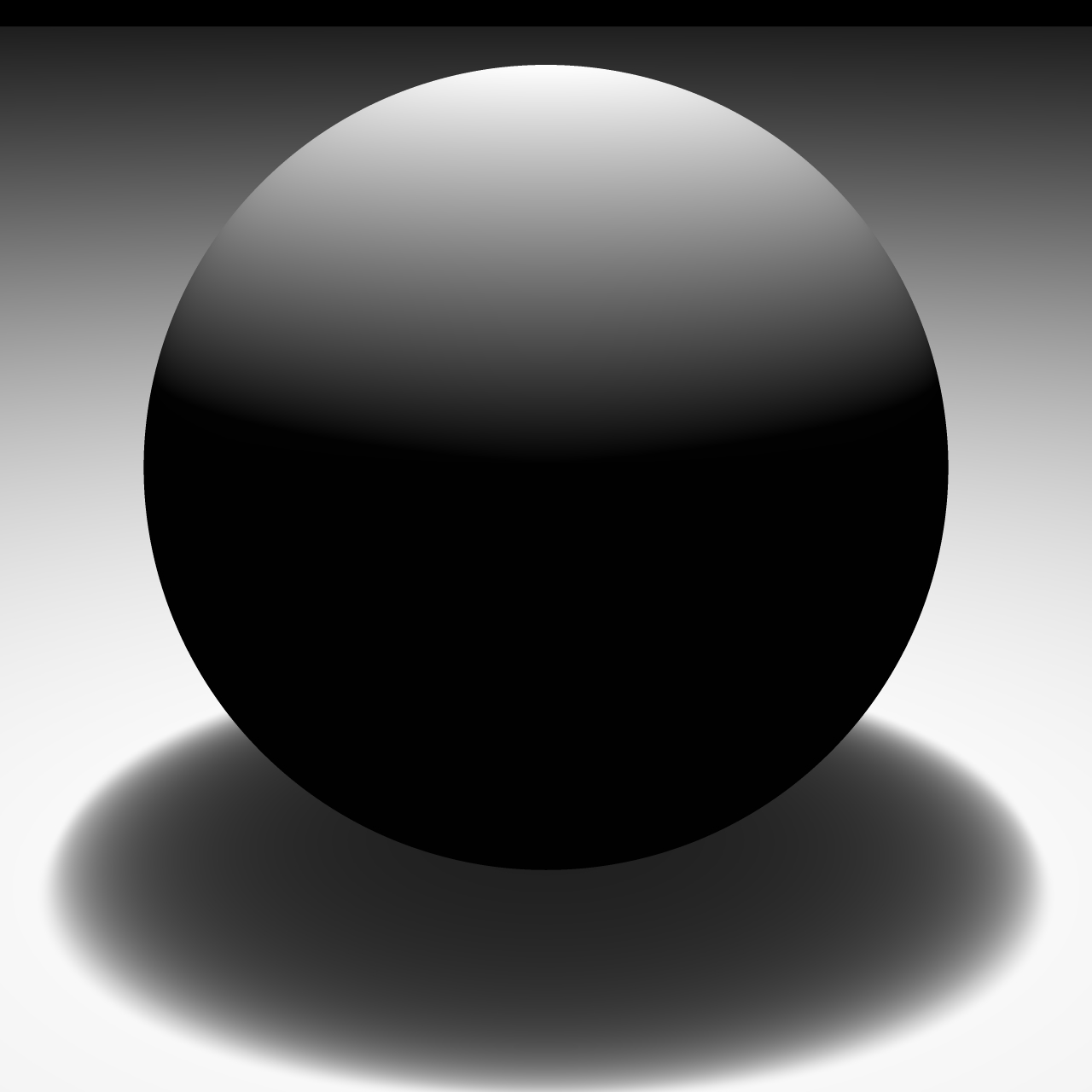}
        \includegraphics[width=1.0\textwidth]{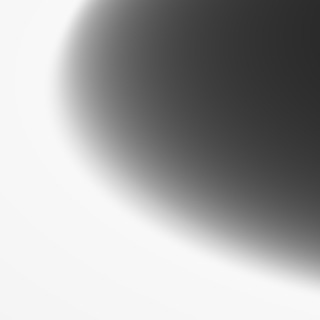}
        \caption{$p=-4.0$}
        \label{fig_sphere_single/8}
    \end{subfigure}
      \hfill 
          \begin{subfigure}[t]{0.32\textwidth}
        \includegraphics[width=1.0\textwidth]{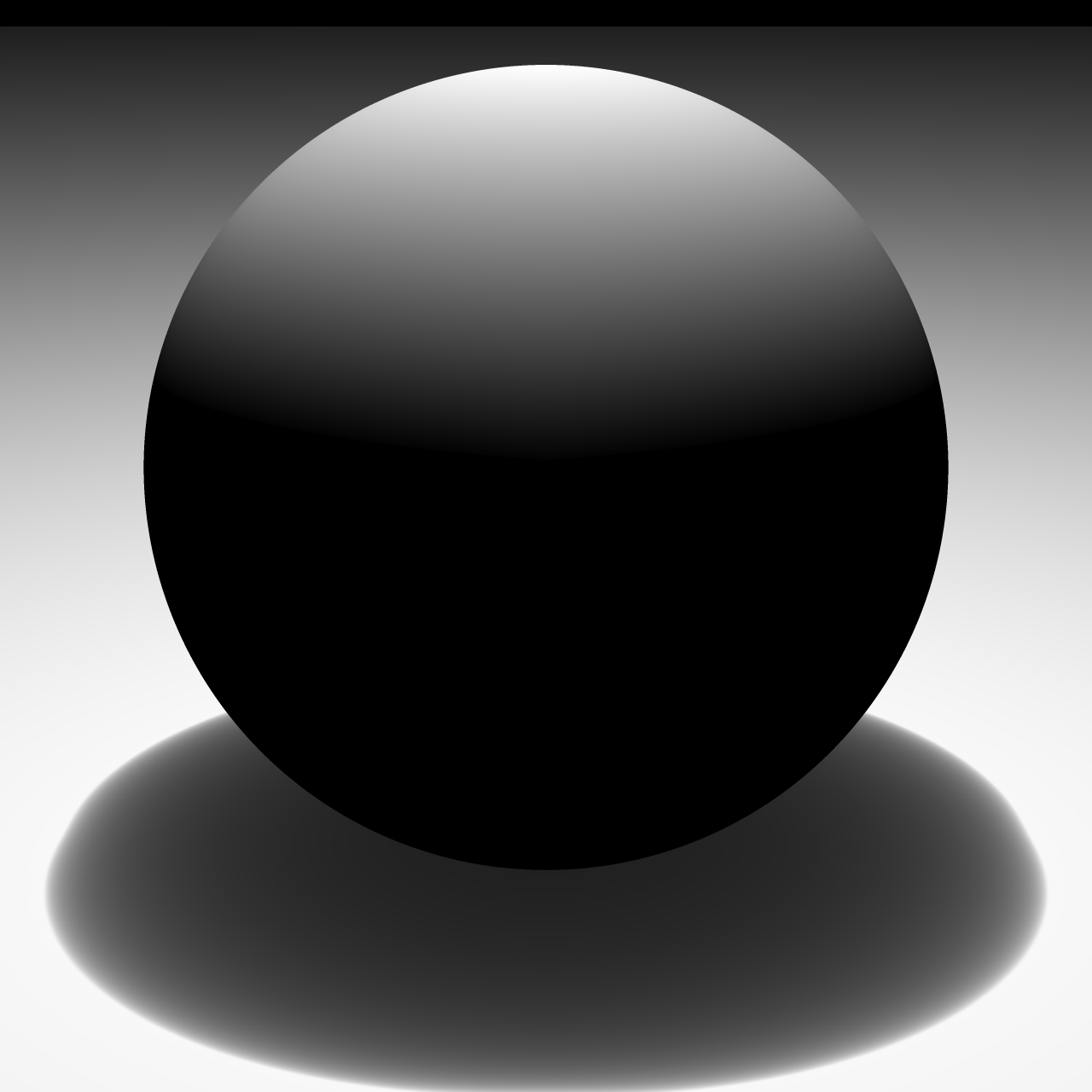}
        \includegraphics[width=1.0\textwidth]{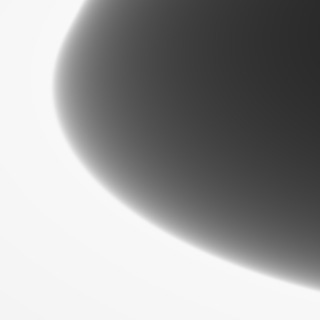}
        \caption{$p \rightarrow -\infty$}
        \label{fig_sphere_single/inf}
    \end{subfigure}
\caption{More $p$ values for example shown in Figure~\ref{fig_teaser2a}. Very small values of $p$  approach the minimum operator and create a sharp shadow boundary}
\label{fig_teaser2b}
\end{figure*}

Another potential approach is to replace the maximum of the minimum operators with some reasonably large values of $p$. One issue that needs attention is that averaging operators can only work with positive numbers or with some special cases, as we have discussed in the paper. For example, we cannot directly change $max(\cos \theta, 0)$ into a Minkowski term. Let $\cos \theta = -1$, then $((-1)^p + 0^p)=-1$ regardless of $p$ value (see subsection~\ref{subsec_remark} that explain the reason). This is definitely not an approximation of the maximum operator. However, there are ways to avoid this problem. For instance, the following formula 
$$C(H_{1}(H_{p}((\cos \theta + 1, 2), (1,1)),(-1,-1))),$$
is an approximation of $max(\cos \theta, 0)$ as
$$\lim_{p \rightarrow \infty}\frac{(\cos \theta +1)^p + 1)^{1/p} +(-1)}{ (2^p + 1)^{1/p} + (-1)} \longrightarrow max(\cos \theta, 0)$$
To exploit the associativity power of any H\"{o}lder-Minkowski Color, we need to use the same $p$ value for all operations. Therefore, it may be better to replace $H_{1}(.,.)$ with $H_{p}(.,.)$ to simplify the equation. We tried to do this by replacing $H_{1}(.,.)$ with $H_{p}(.,.)$ and replacing $(-1,-1)$ with the H\"{o}lder-Minkowski inverse of $(1,1)$. As a result, the term $(1,1)$ is canceled and we ended up with a simpler and nicer formula regardless of the value of $p$ as follows: 
$$ C((H_{p}((\cos \theta + 1, 2), (0,0))) = \frac{\cos \theta +1}{2}.$$
This formula, which was first introduced by Gooch \& Gooch \cite{gooch1998non}, is widely used for non-photorealistic rendering.
In this case, we also have to point out that this particular $\cos \theta$ term is used in part of the material. Therefore, it does not need to be formulated by H\"{o}lder-Minkowski. These approximations can be obtained using spline formulations \cite{bartels1995}. However, such examples demonstrate potential usage beyond color applications. 

An interesting avenue for future work would be to extend the H\"{o}lder-Minkowski average for geometry applications. The H\"{o}lder-Minkowski average can potentially be used in implicit shape modeling as an operation in CSG trees \cite{wyvill1999}. However, this suggests that we allow complex numbers as function values. This change also requires a major paradigm shift in implicit modeling. The relationship with parametric shape modeling is not even obvious. On the other hand, the partition of unity property suggests that there can be a relationship. 

We want to point out that one limitation of H\"{o}lder-Minkowski Colors is that they are not closed under translation. This limitation actually makes sense for dealing with colors. For instance, consider rendering: if we add $(x,0)$ to all lights by changing them exactly the same way, we cannot simply add $(x,0)$ to all rendering results. On the other hand, scale independence makes sense for colors. Consider rendering again: if all lights are multiplied by the same amount, we can scale the resulting computation by the same amount because of the scale Independence, which makes logical sense.

Although the premise of the H\"{o}lder-Minkowski Colors seems to be clear based on our initial results, it is necessary to take small steps using carefully designed case studies. In the future, if the images have also been produced using H\"{o}lder-Minkowski colors, phase terms can be used for further control of the results. For best practice, we must even develop image file formats that can allow complex numbers. We need to convert complex numbers into positive numbers only when we display colors as shown in Figures~\ref{fig_process1} and ~\ref{fig_process2}.

The major cause for aliasing in rendering is discontinuities. They introduce infinite frequencies, and it is not possible to resolve them just by increasing sampling frequencies. Therefore, adding noise to avoid aliasing is still the best method to avoid aliasing \cite{cook1984,cook1986}. On the other hand, if we turn discontinuities into derivative discontinuities at zero, H\"{o}lder-Minkowski operations can be used to remove artifacts. Turning discontinuities into derivative discontinuities at zero for shadows is potentially possible by changing the opacity of the object from one to zero near the boundary. We briefly experimented with this idea and obtained promising results (for a simple example that demonstrates the idea, see Figure~\ref{fig_teaser2a}). However, this direction needs further research.

\bibliographystyle{ACM-Reference-Format}
\bibliography{references,refakleman}

\appendix

\section{Algebraic Properties of Addition on Projective Space}
\label{section_APAPS}

In this appendix, we provide properties of addition. These results are straightforward extensions of our general theorems. We have provided them here for a special treatment of addition operation in case there is a need for some clarification. 

 \begin{figure*}
     \begin{subfigure}[t]{0.32\linewidth}
    \includegraphics[width=1.0\textwidth]{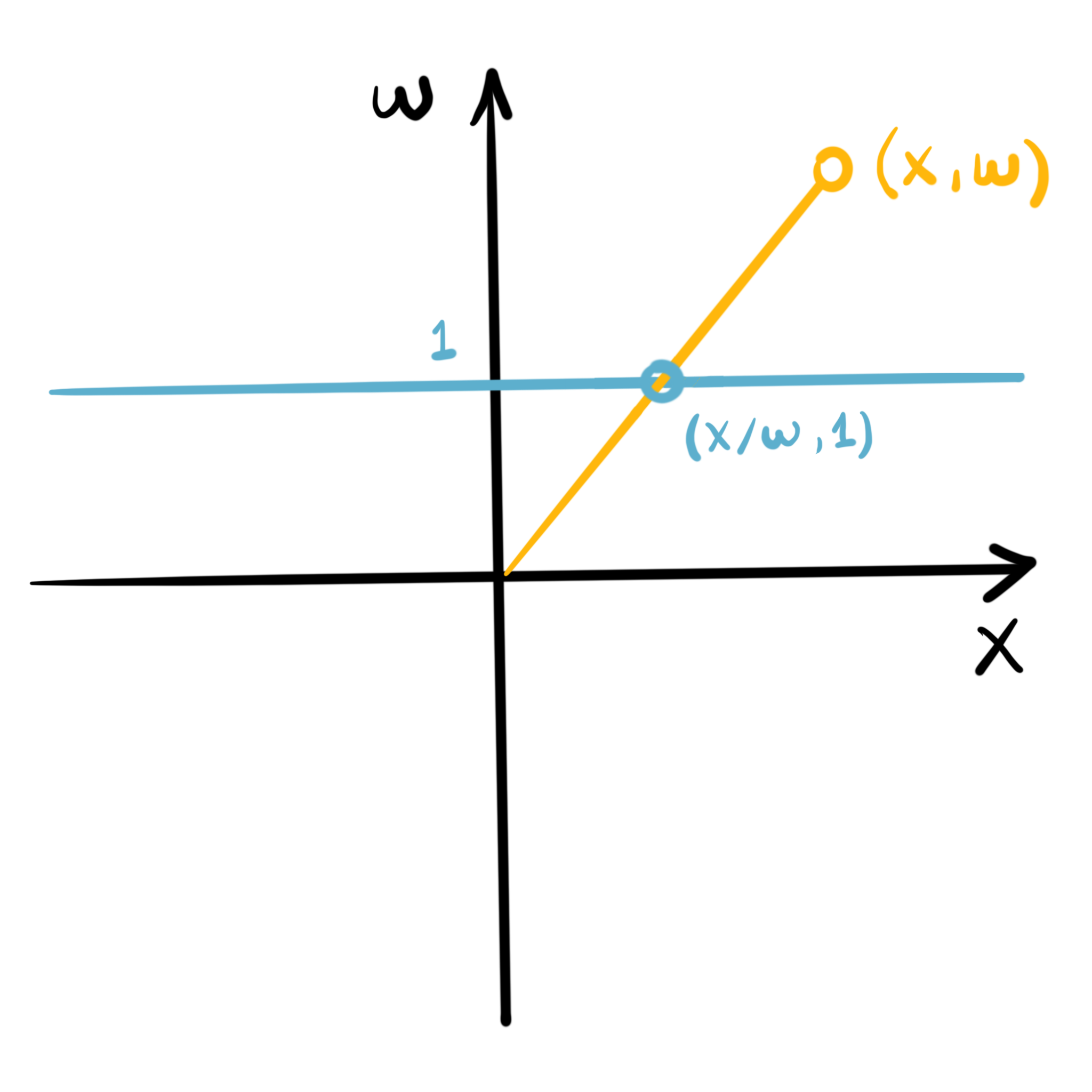}
    \caption{\footnotesize A point in projective space (yellow circle). It can be considered as a ray that starts at $(0,0)$, which intersects with $a =1$. }
    \label{fig_figures0}
 \end{subfigure}
 \hfill
     \begin{subfigure}[t]{0.32\linewidth}
    \includegraphics[width=1.0\textwidth]{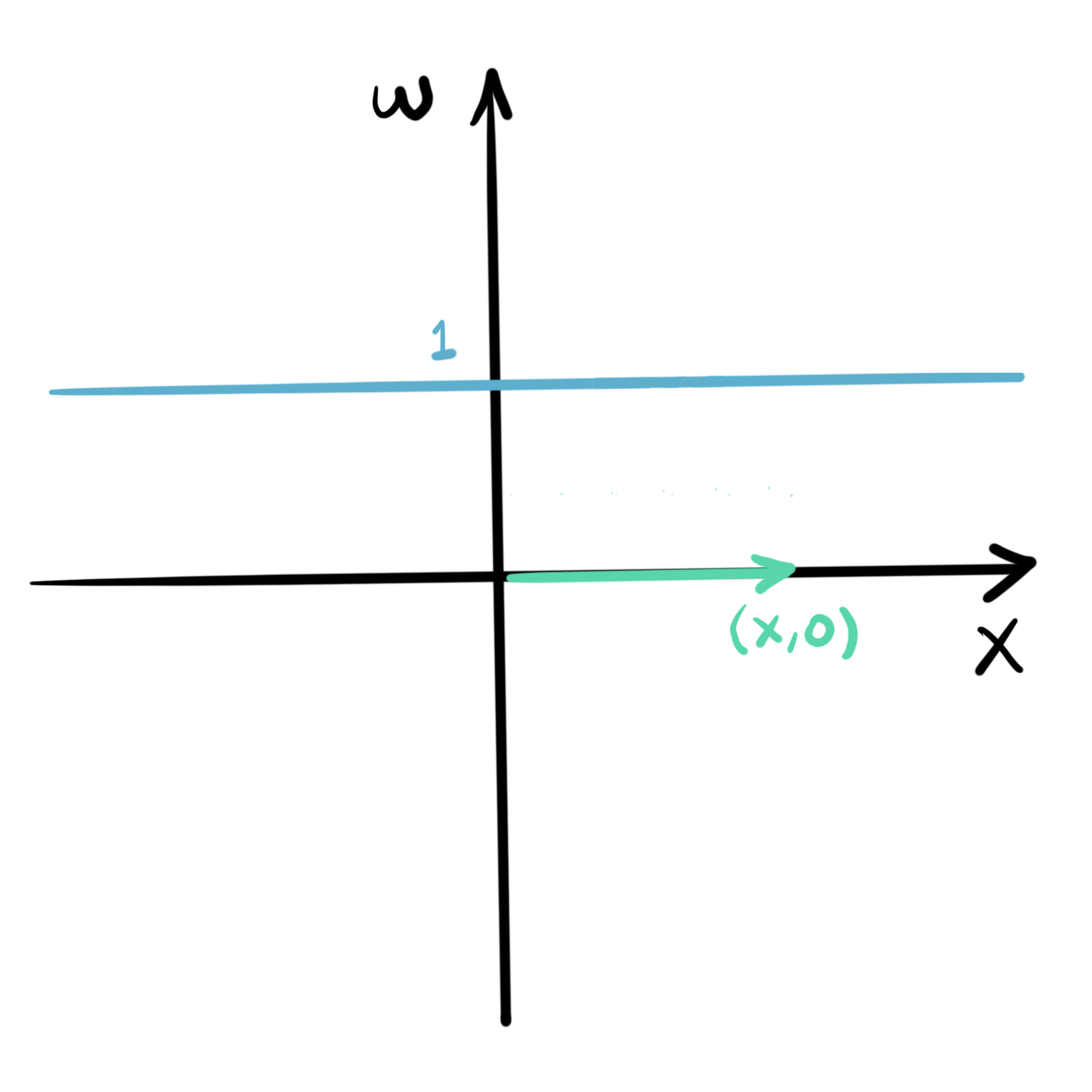}
    \caption{\footnotesize A vector in projective space (cyan arrrow). It can be considered as a ray starts at $(0,0)$ that is parallel to $a =1$.}
    \label{fig_figures1}
 \end{subfigure}
 \hfill
     \begin{subfigure}[t]{0.32\linewidth}
    \includegraphics[width=1.0\textwidth]{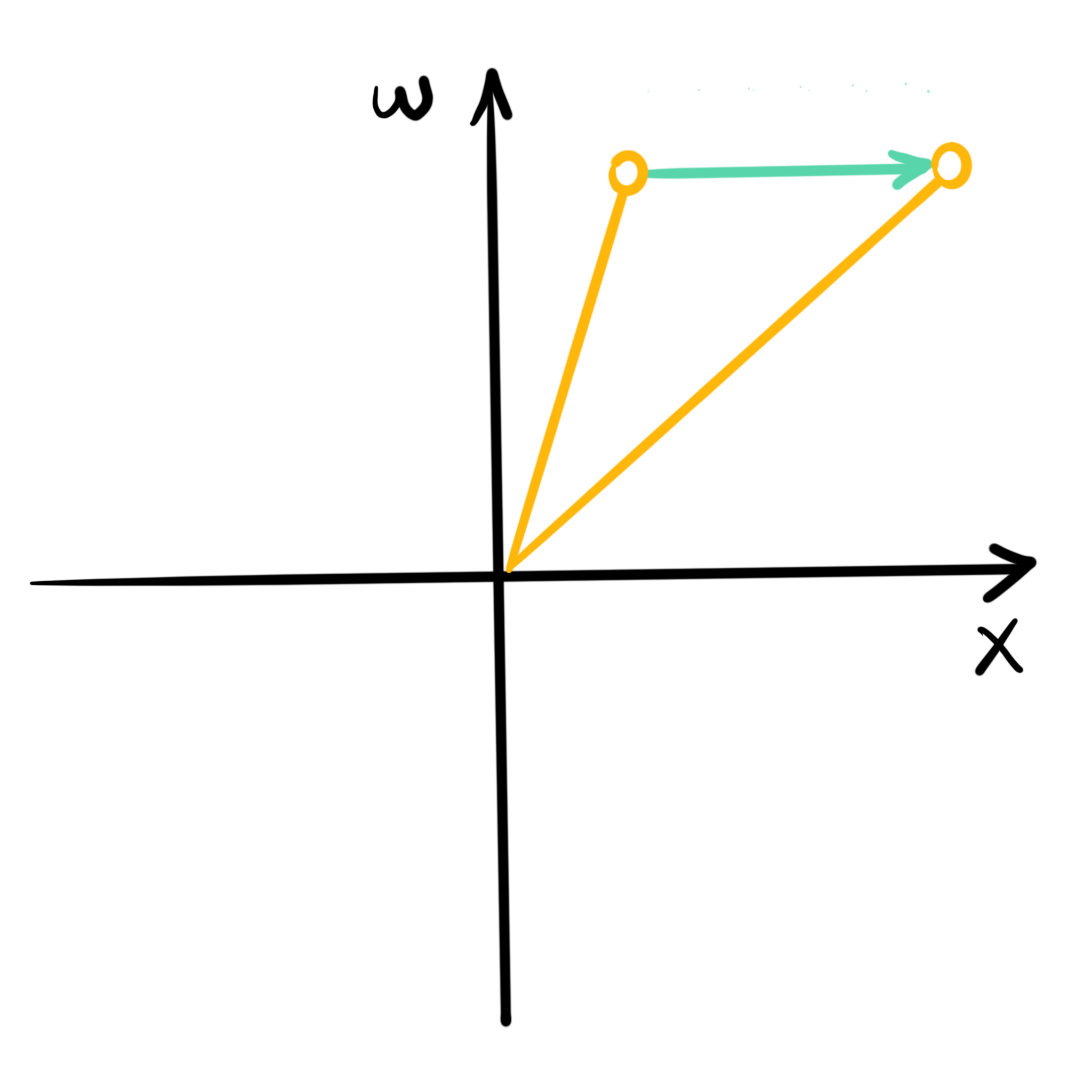}
    \caption{\footnotesize Adding a vector (blue arrow) to a point (yellow circle on the right) produces another point (yellow circle on the left).}
    \label{fig_figures3}
 \end{subfigure}
 \hfill
      \begin{subfigure}[t]{0.32\linewidth}
    \includegraphics[width=1.0\textwidth]{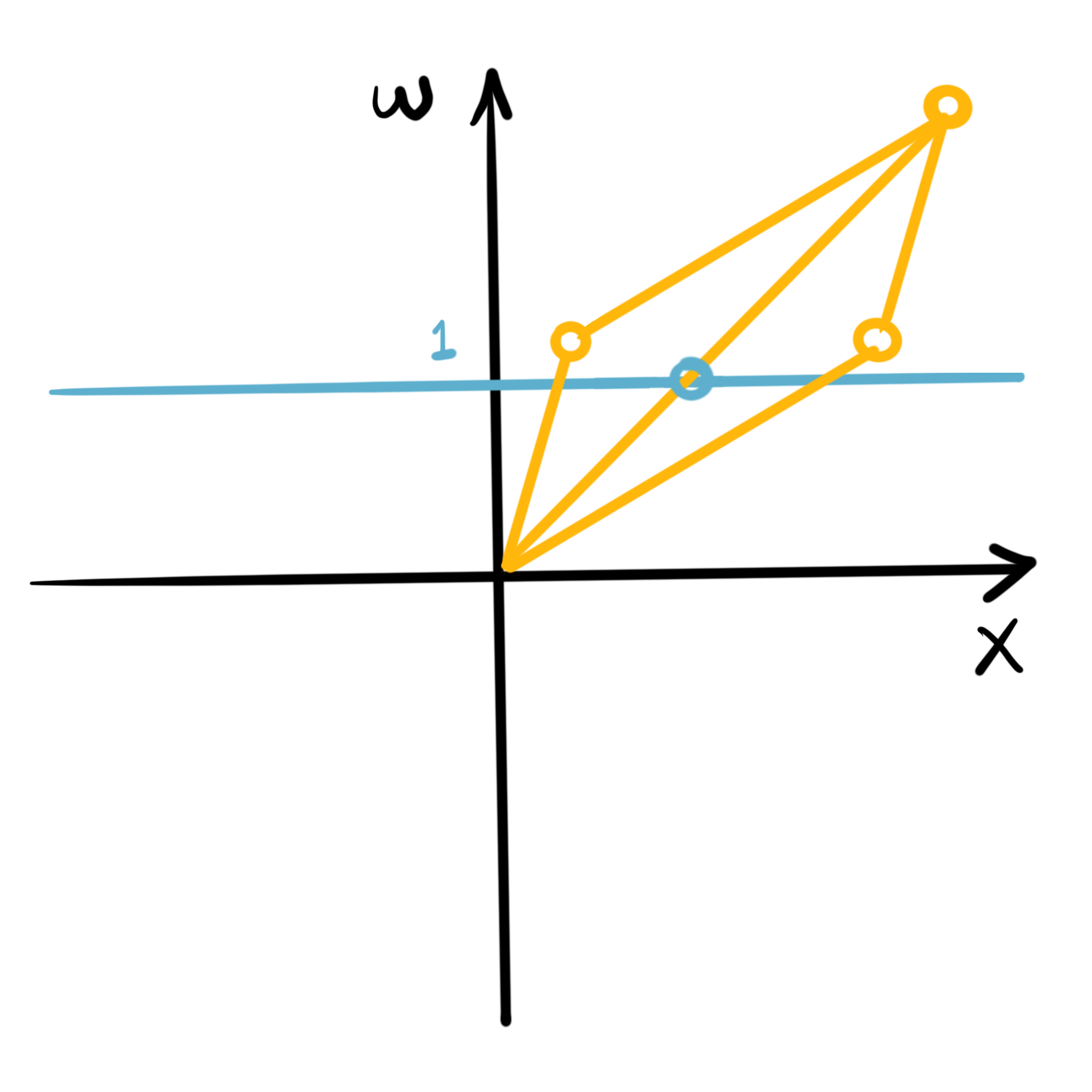}
    \caption{\footnotesize Projection of sum of two points (two yellow circles) creates a point that is the average of projections of these two points.}
    \label{fig_figures4}
 \end{subfigure}
 \hfill
      \begin{subfigure}[t]{0.32\linewidth}
    \includegraphics[width=1.0\textwidth]{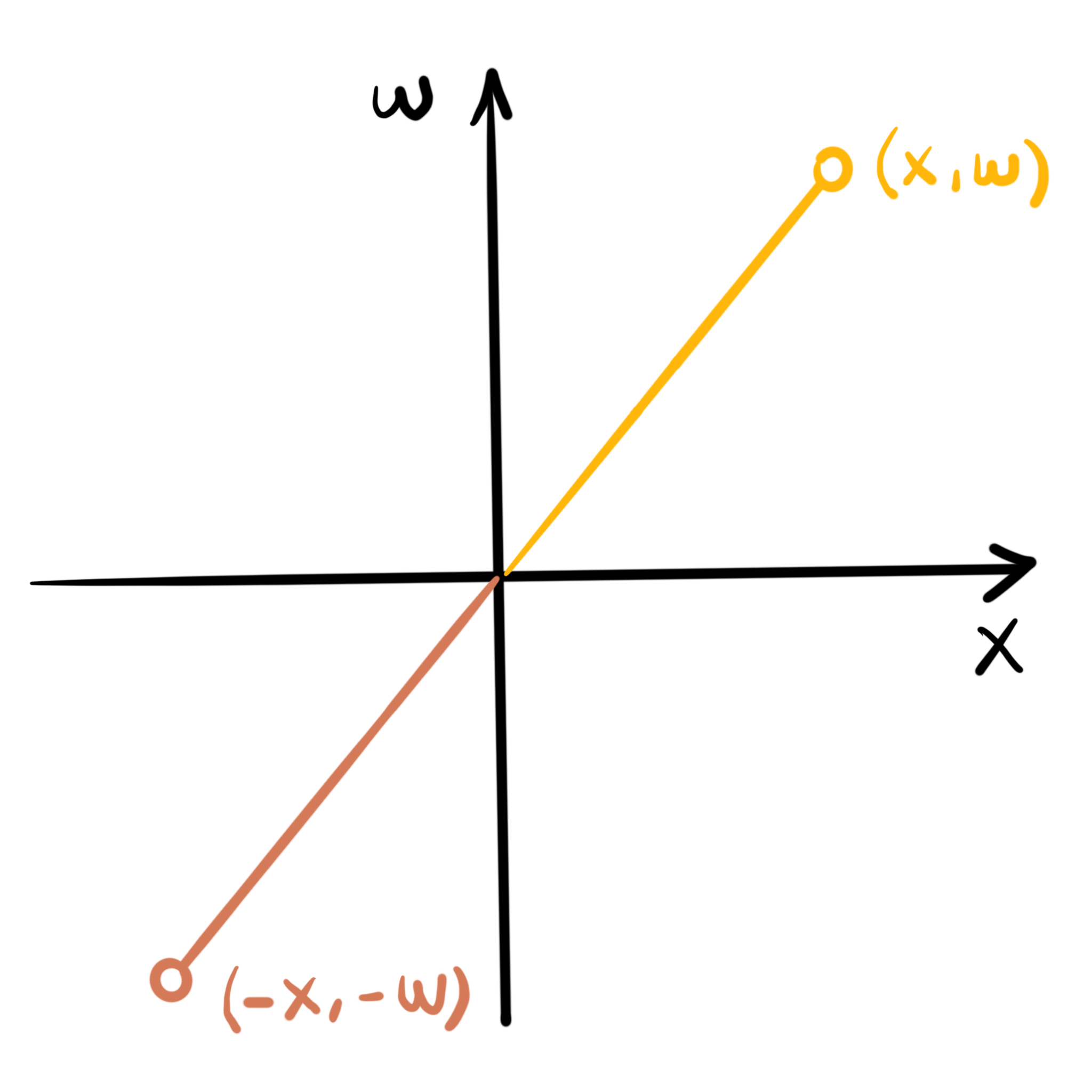}
    \caption{\footnotesize The vector $(0,0)$ is the identity element. And the inverse element for a point $(x,w)$ is $(-x,-w)$ for $H_1$  operation.}
    \label{fig_figures2}
 \end{subfigure}
    \caption{Elements of addition projective space: points \& vectors, operations over them, and identity elements.   
    }
\label{fig_projectivespace}
\end{figure*}

\begin{theorem}
\label{them00}
Let $C(ax,a)=(x,1)$ is called a projection operator and let a set of points be given as follows
$$(a_0 x_0,a_0); (a_1 x_1,a_1);\ldots,(a_{n-1}x_{n-1},a_{n-1}).$$
Then, the projection of their sum is a weighted average of their projected points $$(x_0,1); ( x_1,1);\ldots,(x_{n-1,1}).$$ (See Figure~\ref{fig_projectivespace}). The
process automatically satisfies the partition of unity property if and only if $s_i=a_{i}/a_{0}$ are all positive real numbers. 

\begin{proof}
Note that the projection of their sum can be rewritten as follows:
\begin{eqnarray*}
C\left( \sum_{i=0}^{n-1} a_i x_i,
\sum_{i=0}^{n-1} a_i \right) &=& 
\left( \frac{\sum_{i=0}^{n-1} a_i x_i}{\sum_{i=0}^{n-1} a_i}, 1 \right) \\
&=&
\left( \sum_{i=0}^{n-1} \frac{a_i}{
\sum_{i=0}^{n-1} a_i} x_i, 1 \right)\\
&=&
\left( \sum_{i=0}^{n-1} a_0 \frac{s_i}{ a_0
\sum_{i=0}^{n-1} s_i} x_i, 1 \right)\\
&=&
\left( \sum_{i=0}^{n-1} \frac{s_i}{ 
\sum_{i=0}^{n-1} s_i} x_i, 1 \right)\\
\end{eqnarray*}
Then, the weights satisfy the partition of unity since
$$ 0 \leq  \frac{s_i}{ 
\sum_{i=0}^{n-1} s_i} \leq 1 $$
and 
$$  \sum_{i=0}^{n-1} \left( \frac{s_i}{ 
\sum_{i=0}^{n-1} s_i} \right)= \frac{\sum_{i=0}^{n-1} s_i}{ 
\sum_{i=0}^{n-1} s_i}  =1 $$
This completes the proof. 
\end{proof}
\end{theorem}
In other words, all $a_i$ are on the same ray emanating from the origin, the projection operation still works as a weighted average of projected points even for complex numbers, and it still can produce a point inside the convex hull given by $(x_0,1); ( x_1,1);\ldots,(x_{n-1},1)$ by providing barycentric algebraic properties. The advantage of this formulation is that we can also extrapolate in a controlled manner. 

\begin{theorem}
\label{them01}
Let two points be given as follows
$$(a_0 x_0,a_0); (a_1 x_1,a_1)$$ and $s_1=a_{1}/a_{0}$ be any real number. Then the projection of their sum can be any point on the $(x_0,1); (x_1,1)$ line. 

\begin{proof}
Note that the projection of their sum can be rewritten as follows:
\begin{eqnarray*}
C\left( a_0 x_0+ a_1 x_1, a_0+a_1 \right) &=& 
\left( \frac{a_0 x_0+ a_1 x_1}{a_0+a_1}, 1 \right) \\
&=&
\left( \frac{a_0 (x_0+ s_1 x_1)}{a_0( 1+s_1)}, 1 \right) \\
&=&
\left( \frac{(x_0+ s_1 x_1)}{( 1+s_1)}, 1 \right) \\
&=&
\left( \frac{(s_1 (x_1-x_0)+ (1+s_1) x_1)}{( 1+s_1)}, 1 \right) \\
&=&
\left( x_0 + \frac{s_1}{( 1+s_1)} (x_1-x_0), 1 \right) \\
\end{eqnarray*}
This is an equation of a line segment that connects $(x_0,1) and (x_1,1)$ for positive real $s_1$ values. This is consistent with the previous theorem. For real values of $s_1$ this gives us an equation of an infinite line that passes from $(x_0,1) and (x_1,1)$. 
This completes the proof. 
\end{proof}
\end{theorem}
This property is also useful since we can cover the whole plane including $\infty$ without explicitly stating it since $\lim s_1 \rightarrow -1$, the line goes to $\infty$. If $a_0$ and $a_1$ are relatively complex (that is, $s_1$ is a complex number), then we can obtain any point in the complex plane. We will not prove the previous claim, since it is obvious. This property provides intuitive control as follows. Positive real $s_1$ values provide a line segment, real $s_1$ values provide the whole line, and non-zero imaginary components allow us to cover entire complex domain. 
This representation has additional power. It turns the weighted sum operation into an associative operator.

\begin{theorem}
\label{them02}
Addition is associative. 
\begin{proof}
Let three points be given as follows
$$(a_0 x_0,a_0); (a_1 x_1,a_1); (a_2 x_2,a_2).$$ Now, let us check the associative property. By adding first $(a_0 x_0,a_0)$ and $(a_1 x_1,a_1)$. Then adding $(a_2 x_2,a_2)$ to the sum we obtain the following: 
\begin{eqnarray*}
\left( (a_0 x_0,a_0) + (a_1 x_1,a_1) \right) + (a_2 x_2,a_2) = \\
\left( (a_0 x_0 + a_1 x_1,a_0+a_1) \right) + (a_2 x_2,a_2) =\\
\left( (a_0 x_0 + a_1 x_1 + a_2 x_2 ,a_0+a_1 + a_2) \right)  \\
\end{eqnarray*}
Next, first add $(a_1 x_1,a_1)$ and $(a_2 x_2,a_2)$. Then adding $(a_0 x_0,a_0)$ to the sum we obtain the following:
\begin{eqnarray*}
(a_0 x_0,a_0) +\left(  (a_1 x_1,a_1) + (a_2 x_2,a_2) \right)  = \\
(a_0 x_0,a_0) + \left( (a_1 x_1 + a_2 x_2,a_1+a_2) \right) =\\
\left( (a_0 x_0 + a_1 x_1 + a_2 x_2 ,a_0+a_1 + a_2) \right)  \\
\end{eqnarray*}
They are exactly the same, and this completes the proof. 
\end{proof}
\end{theorem}
This particular case is trivial; we demonstrate associativity later for general operations. The addition operation is also commutative. We do not provide it since it is even simpler than the previous proof. However, the impact of this property is that this formulation is different from an associative alpha operation, which is not commutative. Non-commutativity is a desired feature for some applications that involve pixel-level occlusion. On the other hand, for many other rendering operations, commutativity is useful when there is no order dependency. 

This introduction sets the stage for the next subject. In the next section, we demonstrate how to obtain inverse operations. 

 \subsection{Identity and Inverse Elements}
 
We now want to point out that the identity element for addition exists in homogeneous coordinates and is $(0,0)$. Note that these are actually null functions that map the entire domain to zero that can be written as $\mathbf{0}$. This is useful since we can now define inverse elements using an identify element. 
 
Now we can identify the inverse elements for any given point or vector. Let a point or vector be denoted as $(x_0,a_0)$, where $a_0$ can be zero to include vectors. Let the inverse of this point or vector denoted by again as $(x_1,a_1)$, because of the definition of identity and inverse elements, the following equality must hold: 
$$(x_0,a_0)+ (x_1,a_1) = (0,0)$$
This means that 
\begin{eqnarray*}
x_1 &=&  -x_0 \\
a_1 &=&  -a_0 \\
\end{eqnarray*}
The first term is not unexpected. We already let negative numbers be used to describe positions in the 3D space. On the other hand, the second term suggests that we need negative numbers for $a$. Even the first term is interesting in the sense that if we want to use the new representation for colors, we also need to allow negative colors. This observation is critical for problem-solving on color space, especially for compositing. In current practice, digital compositing is a set of hacks. This provides an extended algebra to handle digital compositing operations over a space of colors. This extension means that single channel scalars such as $z,a$ can be any real numbers. 

\begin{figure}[hbtp]
\centering
     \begin{subfigure}[t]{0.45\linewidth}
    \includegraphics[width=1.0\textwidth]{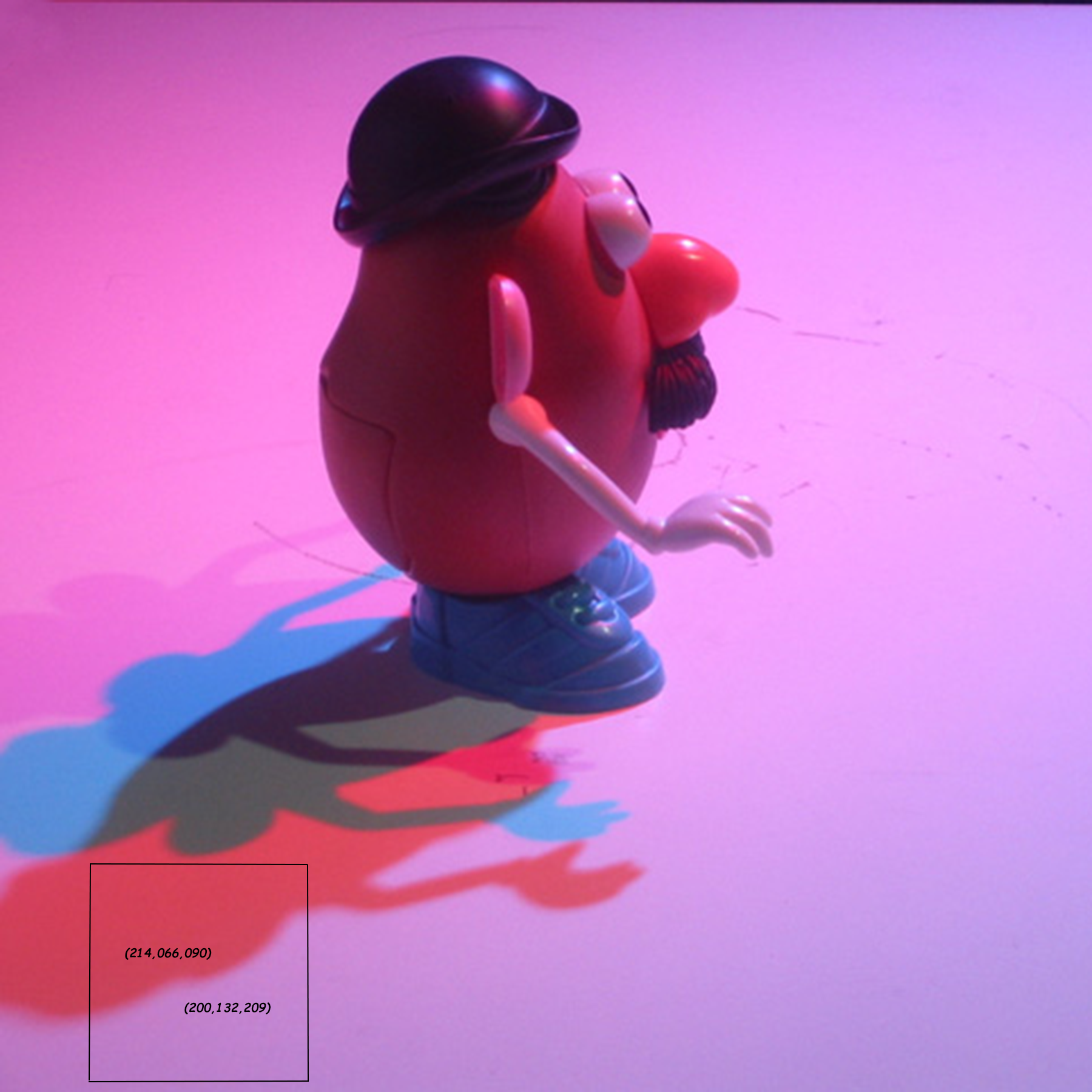}
    \caption{\footnotesize A photograph of Mr. Potato head illuminated by three colored lights. }
    \label{fig_photos/0}
 \end{subfigure}
 \hfill
     \begin{subfigure}[t]{0.45\linewidth}
    \includegraphics[width=1.0\textwidth]{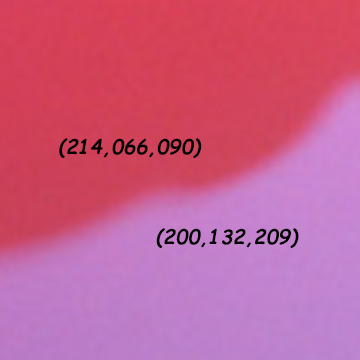}
    \caption{\footnotesize A detail of the photo that shows the red channel in the shadow region is brighter than in the non-shadow region. }
    \label{fig_photos/0a}
 \end{subfigure}
 \hfill
    \caption{An example that demonstrates we need negative numbers to solve digital compositing problems.   
    }
\label{fig_photos/00a}
\end{figure}

\begin{figure}[hbtp]
\centering
     \begin{subfigure}[t]{0.45\linewidth}
    \includegraphics[width=1.0\textwidth]{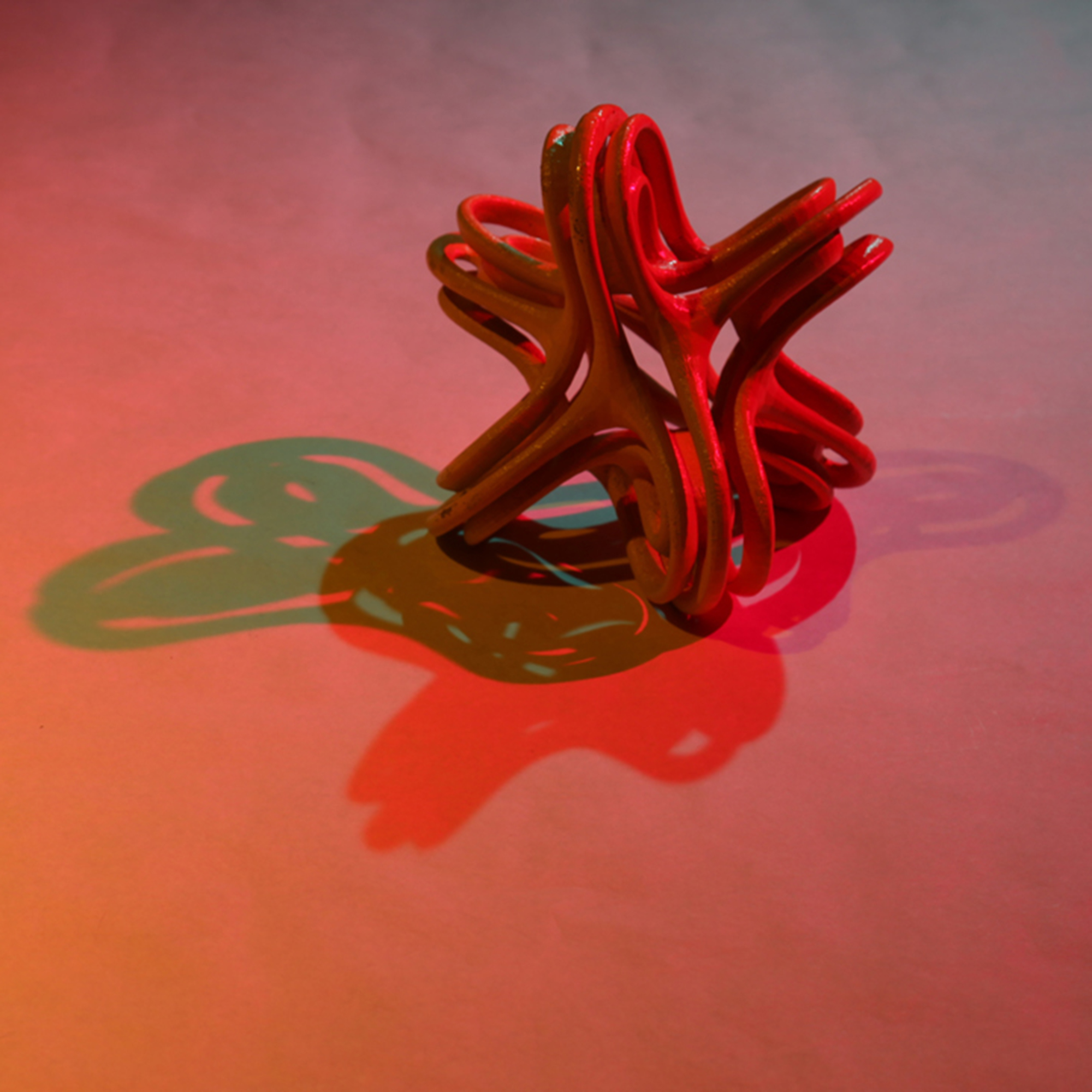}
    \caption{\footnotesize A photograph of a high-genus sculpture illuminated by three colored lights. }
    \label{fig_photos/3}
 \end{subfigure}
 \hfill
     \begin{subfigure}[t]{0.45\linewidth}
    \includegraphics[width=1.0\textwidth]{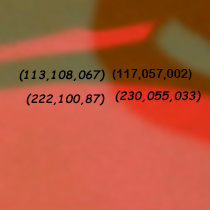}
    \caption{\footnotesize A detail of the photo that shows both red and green channels in the shadow regions are brighter than in the non-shadow region. }
    \label{fig_photos/3a}
 \end{subfigure}
 \hfill
    \caption{Another example that demonstrates we need negative numbers to solve digital compositing problems.   
    }
\label{fig_photos/33a}
\end{figure}

Now, let us demonstrate the effect of this change in real photographs.  Figure~\ref{fig_photos/00a} shows a photograph of a real object, Mr. Photo head, taken under three colored lights. Note that we gave details of an area that is the interface of shadow and illuminated regions. Since this area is small enough, we can assume that the change is linear and it is only effected by one of the lights. In this case, the colors in the shadow and illuminated regions are picked as
$$\mathbf{c_0}=(214,066,090)$$
$$\mathbf{c_1}=(200,132,209)$$
respectively. In the projective space, they are given as $(\mathbf(c_1),1)$ and $(\mathbf(c_0),1)$. 
Assume that light is a binary function $t$, which is $1$ and $0$. Our linear parametric function will be in the classical polynomial form as 
$$(\mathbf{c},1)=(\mathbf{a},1) + t (\mathbf{b},0)$$
if we compute two parameters of this linear function, we find 
$$(\mathbf{a},1)= (\mathbf(c_0),1) = (214,066,090,1) $$
$$(\mathbf{b},0) = (\mathbf(c_1),1) - (\mathbf(c_0),1) =  (\mathbf(c_1),1) + (-\mathbf(c_0),-1) = (-14, 66, 119, 0)$$
Note that the first term becomes negative. Note that the linear formula can also be written in barycentric form as an interpolation of the two colors. In that case, we do not need negative numbers as 
$$(\mathbf{c},1)=(\mathbf(c_1),1) t + (\mathbf(c_0),1) (1-t) $$
Now, let us check another real photograph case shown in Figure~\ref{fig_photos/33a}. In this case, we have another detail region that is at the intersection of two shadow regions. Assume that the colored lights are given $t_0$ and $t1$ and that they can only take Boolean values, i.e. $0$ or $1$. In this case, the colors in four regions are picked as
$$\mathbf{c_{00}}=(117,57,2)$$
$$\mathbf{c_{01}}=(230,55,33)$$
$$\mathbf{c_{10}}=(113,108,67)$$
$$\mathbf{c_{11}}=(222,100,87)$$
If we again assume that the region is small enough, we can interpolate these color points with a bilinear equation as follows: 
\begin{eqnarray*}
(\mathbf{c},1) &=& (1-t_0)(1-t_1)(\mathbf{c_{00}},1)+ (1-t_0)t_1 (\mathbf{c_{01}},1) +\\
&& t_0(1-t_1) \; (\mathbf{c_{10}},1) + t_0t_1 (\mathbf{c_{11}},1) \\
\end{eqnarray*}
In this barycentric form, all the numbers are positive. However, if we want to use the classical polynomial form, we need 
\begin{eqnarray*}
(\mathbf{c},1) &=& (\mathbf{c_{00}},1)+ \\
&&t_1 ((\mathbf{c_{01}},1) + (-\mathbf{c_{00}},-1))+\\
&&t_0 ((\mathbf{c_{10}},1) + (-\mathbf{c_{00}},-1))+\\ 
&&t_0 t_1((\mathbf{c_{11}},1) +(-\mathbf{c_{10}},-1) +(-\mathbf{c_{01}},-1) + (\mathbf{c_{00}},1))\\
&=& (117,57,2,1,1)+ \\
&&t_1 (113,-2,31,0)+\\
&&t_0 (-4,51,65, 0)+\\ 
&&t_0 t_1(-4, -6, -11,0 )\\
\end{eqnarray*}
Here, it is interesting to note that all values of the last vector is negative. Having auxiliary parameters with negative values helps the computation. If we did not allow negative values, we would not be able to form a model to obtain the actual values. In this formulation, it is also important to point out that negative values usually come in terms of vectors since vectors can be obtained as the difference of two points.  

Now, we have a commutative and associative algebra in our version of projective space. Next, we include Minkowski operations on this projective space, which will require complex number to develop a consistent algebra with commutativity and associativity property.  We will demonstrate commutativity and associativity properties in a more general setting.

\end{document}